\documentclass[10pt,epsfig]{article} 
\usepackage{amsbsy,amssymb,amscd,amsfonts,latexsym,amstext,delarray, 
amsmath,epsfig}  \headsep=15pt 
 
\font\cyr=wncyr10 
 
\input xypic 
 
\newtheorem{thm}{Theorem}[section] 
\newtheorem{prop}[thm]{Proposition} 
\newtheorem{cor}[thm]{Corollary} 
\newtheorem{lem}[thm]{Lemma} 
 
\newtheorem{defn}[thm]{Definition} 
\newtheorem{rem}[thm]{Remark} 
\newtheorem{ex}[thm]{Example} 
\newtheorem{claim}[thm]{Claim}

\numberwithin{equation}{section}

\def\R{{\mathbb R}}

\def\H{{\mathbb H}}

\def\mX{{\mathfrak X}} 
\def\Sp{{\rm{Spec}}}

\def\Q{{\mathbb Q}} 
\def\C{{\mathbb C}} 
\def\Z{{\mathbb Z}} 
\def\N{{\mathbb N}} 
\def\HH{{\mathbb H}} 
\def\P{{\mathbb P}} 
 
\def\PSL{{\rm PSL}} 
 
\def\SL{{\rm SL}} 
\def\Aut{{\rm Aut}} 
\def\Tr{{\rm Tr}} 
\def\w{{\text w}} 
\def\X{{{\tilde X^*}}} 
\def\XX{{{X^*}}} 
\def\K{{\mathbb K}} 

\def\cancel#1#2{\ooalign{$\hfil#1\mkern1mu/\hfil$\crcr$#1#2$}} 
\def\dirac{\mathpalette\cancel\partial}

\newcommand{\ie}{{\it i.e.\/}\ } 
\newcommand{\eg}{{\it e.g.\/}\ } 
\newcommand{\cf}{{\it cf.\/}\ } 
\newcommand{\resp}{{\it resp.\/}\ } 
\newcommand{\op}{{\it op.cit.\/}\ } 
 
\begin{document} 
\title{Non--commutative geometry, dynamics, and $\infty$--adic Arakelov 
geometry} 
 
\author{Caterina Consani\thanks{Partially supported by NSERC 
grant 72016789} and Matilde Marcolli 
\thanks{Partially supported by Humboldt Foundation 
Sofja Kovalevskaja Award}} 
\date{} 
\maketitle 
 
\bigskip 
 
\begin{center} 
{\em We dedicate this work to Yuri Manin, with admiration 
and gratitude} 
\end{center} 
 
\begin{abstract} 
In Arakelov theory a completion of an arithmetic surface is 
achieved by enlarging the group of divisors by formal linear 
combinations of the ``closed fibers at infinity''. Manin described 
the dual graph of any such closed fiber in terms of an infinite 
tangle of bounded geodesics in a hyperbolic handlebody endowed 
with a Schottky uniformization. In this paper we consider 
arithmetic surfaces over the ring of integers in a number field, 
with fibers of genus $g\geq 2$. We use Connes' theory of spectral 
triples to relate the hyperbolic geometry of the handlebody to 
Deninger's Archimedean  cohomology and the cohomology of the cone 
of the local monodromy $N$ at arithmetic infinity as introduced by 
the first author of this paper. First, we consider derived (cohomological) 
spectral data $({\rm A},H^\cdot(\XX),\Phi)$, where the algebra 
is obtained from the
$\SL(2,\R)$ action on the cohomology of the cone, induced by the 
presence of a polarized Lefschetz module structure, and its
restriction to the group ring of a Fuchsian Schottky group. In this 
setting we recover the alternating product of the Archimedean  
factors from a zeta function of a spectral triple. Then, we 
introduce a different construction, which is related to Manin's 
description of the dual graph of the fiber at infinity. We provide a
geometric model for the dual graph as the mapping torus of a dynamical
system $T$ on a Cantor set. We consider a noncommutative space which
describes the action of the Schottky group on its limit set and
parameterizes the ``components of the closed fiber at infinity''. 
This can be identified with a Cuntz--Krieger algebra ${\mathcal O}_A$
associated to a subshift of finite type. We construct a spectral
triple for this noncommutative space, via a representation on the 
cochains of a ``dynamical cohomology'', defined in terms 
of the tangle of bounded geodesics in the handlebody. In both
constructions presented in the paper, the Dirac operator agrees with
the grading operator $\Phi$, that represents the ``logarithm of a
Frobenius--type operator'' on the Archimedean  cohomology. In fact,
the Archimedean  cohomology embeds in the dynamical 
cohomology, compatibly with the action of a real Frobenius $\bar
F_\infty$, so that the local factor can again be recovered from these
data. The duality isomorphism on the cohomology of the cone of $N$
corresponds to the pairing of dynamical homology and cohomology. This
suggests the existence of a duality between the monodromy $N$ and the 
dynamical map $1-T$. Moreover, the ``reduction mod infinity'' is
described in terms of the homotopy quotient associated to the
noncommutative space ${\mathcal O}_A$ and the $\mu$-map of
Baum--Connes. The geometric model of the dual graph can also be
described as a homotopy quotient. 
\end{abstract} 
 
\vskip .5in 
 
\noindent {\bf Caterina Consani}, Department of Mathematics, University of 
Toronto, Canada. 
 
\noindent email: kc\@@math.toronto.edu 
 
\vskip .3in 
 
\noindent {\bf Matilde Marcolli}, Max--Planck--Institut f\"ur 
Mathematik, Bonn Germany. 
 
\noindent email: marcolli\@@mpim-bonn.mpg.de

\newpage 
 
\begin{verse} 
 
{\cyr Zolotoe runo, gde zhe t\char121{}, zolotoe runo? \\ 
Vsyu dorogu shumeli morskie tyazhel\char121{}e voln\char121{}, \\ 
i pokinuv korabl\char126{}, natrudivshi\char26{} v moryakh polotno, \\ 
Odisse\char26{} vozvratilsya, prostranstvom i vremenem 
poln\char121{}\char26{} \\ 
(Osip Mandel\char126{}shtam) } 
 
\end{verse} 
 
\tableofcontents 
 
\newpage 
 
\section{Introduction.}\label{0}

The aim of this paper is to show how non--commutative geometry 
provides a connection between two constructions in Arakelov theory 
concerning the Archimedean  fibers of a one-dimensional arithmetic fibration. 
On the one hand, we consider the cohomological construction introduced 
by the first author in \cite{KC}, that was inspired by the theory of 
limiting mixed Hodge structures on the limit fiber of a geometric 
degeneration over a disc and it is related to Deninger's 
Archimedean cohomology and regularized determinants 
(\cf \cite{Den}). On the other hand, we re--interpret Manin's geometric 
realization of the dual graph of the fiber at infinity as an infinite 
tangle of bounded geodesics inside a real 3--dimensional hyperbolic 
handlebody (\cf \cite{Man}) in the context of dynamical systems. 
The problem of relating the results of \cite{Man} with the 
cohomological constructions of Deninger was already addressed by Manin 
in \cite{Man}, but to our knowledge no further progress in this 
direction was made since then. 
 
Let $\K$ be a number field and let $O_\K$ be the ring of 
integers. The choice of a model $X_{O_\K}$ of a smooth, algebraic 
curve $X$ over $\K$ defines an arithmetic surface over 
$\Sp(O_\K)$. A closed vertical fiber of $X_{O_\K}$ over a prime $\wp$ 
in $O_\K$ is given by $X_{\wp}$: the reduction mod~$\wp$ of the 
model. It is well known that a completion of the fibered surface 
$X_{O_\K}$ is achieved by 
adding to $\Sp(O_\K)$ the {\em Archimedean  places} represented by the 
set of all embeddings $\alpha: \K \hookrightarrow \C$. The Arakelov 
divisors on the completion $\overline{X_{O_\K}}$  are 
defined by the divisors on $X_{O_\K}$ and by 
formal real combinations of the closed vertical fibers at 
infinity. Arakelov's geometry does not provide an explicit 
description of these fibers 
and it prescribes instead a Hermitian metric on each Riemann surface 
$X_{/\C}$, for each Archimedean prime $\alpha$. It is quite remarkable 
that the Hermitian geometry on each $X_{/\C}$ is sufficient to develop 
an intersection theory on the completed model, without 
an explicit knowledge of the closed fibers at infinity. 
For instance, Arakelov showed that intersection indices of divisors on the 
fibers at infinity are obtained via Green functions on the Riemann 
surfaces $X_{/\C}$. 
 
Inspired by Mumford's p-adic uniformization of algebraic 
curves \cite{Mum}, Manin realized that one could 
enrich Arakelov's metric structure by a choice of a Schottky 
uniformization. In this way, the Riemann surface $X_{/\C}$ is the 
boundary at infinity of a $3$-dimensional hyperbolic handlebody 
${\mathfrak X}_\Gamma$, 
described as the quotient of the real hyperbolic 3--space $\H^3$ by the 
action of the Schottky group $\Gamma$. The handlebody contains in its 
interior an infinite link of bounded geodesics, which are interpreted 
as the dual graph of the closed fiber at infinity, thus 
providing a first geometric realization of that space. 
 
A consequence of this innovative approach is a more concrete 
intuition of the idea that, in Arakelov geometry, 
the ``reduction modulo infinity'' of an arithmetic variety should 
be thought of as ``maximally degenerate'' (or totally split: all 
components are of genus zero). This is, in fact, 
the reduction type of the special fiber admitting a Schottky 
uniformization (\cf \cite{Mum}). 
 
In this paper we consider the case of an arithmetic surface over $\Sp 
(O_\K)$ where the fibers are of genus $g\ge 2$. 
The paper is divided into two parts. 
 
The first part consists of Sections \ref{3} and \ref{4}. 
Here we consider the formal construction of a cohomological theory for 
the ``maximally degenerate'' 
fiber at arithmetic infinity, developed in \cite{KC}. Namely, the 
Riemann surface $X_{/\C}$ supports a double complex $(K^{\cdot,\cdot}, 
d', d'')$ endowed with an endomorphism $N$. This complex is made of 
direct sums of vector spaces of real differential forms with certain 
``cutoff'' conditions on the indices, and was constructed as an 
Archimedean  analogue of the one 
defined by Steenbrink on the semi-stable 
fiber of a degeneration over a disc \cite{St}. 
The hyper-cohomologies of $(K^{\cdot}, d=d'+d'')$ and $({\rm Cone}(N)^\cdot, 
d)$ are infinite dimensional, graded real vector 
spaces. We show that their summands are isomorphic 
twisted copies of a same real de Rham cohomology group of $X$. The 
arithmetic meaning of $K^{\cdot}$ arises from the fact that the 
cohomology of $({\rm Coker}(N)^\cdot, d)$ computes the real Deligne 
cohomology of $X_{/\C}$, 
and the regularized determinant of an 
operator $\Phi$ on the subspace $\H^\cdot(K^\cdot,d)^{N=0}$ of the 
hyper-cohomology of the complex $(K^{\cdot}, d)$ recovers the 
Archimedean factors of \cite{Den}. 
The complex $K^{\cdot}$ carries an important structure of bigraded 
polarized Lefschetz module \`a la Deligne and Saito (\cf~\cite{Saito}). In 
particular one obtains an induced inner 
product on the hyper-cohomology and a representation of $\SL(2,\R)\times 
\SL(2,\R)$. 
 
The first part of this paper concentrates on the cohomology 
$H^\cdot(\XX)$ of $Cone(N)$. 
In the classical case of a semi--stable 
degeneration over a disc, the cohomology 
$H^\cdot(|{\mathcal G}|)$ of the dual graph of the special fiber can 
be described in terms of graded pieces, under the monodromy filtration 
on the cohomology of the geometric generic fiber, so that 
$H^\cdot(|{\mathcal G}|)$ 
provides at least a partial information on the mixed Hodge structure on 
$H^\cdot(\XX)$ (here $\XX$ denotes the complement of the special 
fiber in the model and the cohomology of this space has a 
second possible description as hypercohomology of the complex 
${\rm Cone}(N)^\cdot$). In the arithmetic case of a 
degeneration on the ring of integers of a local field, the cohomology 
group $H^\cdot(\XX)$ is still endowed with a graded structure, which is 
fundamental in arithmetic for determining the 
behavior of the local Euler factors at integer points on the left of the 
critical strip on the real line. In fact, the cohomology 
$H^\cdot(\XX)$ carries more arithmetical information than just the 
cohomology of the dual graph of 
the special fiber $H^\cdot(|{\mathcal G}|)$. \vspace{.05in} 
 
Using non--commutative geometry, we interpret the data 
of the cohomology $H^\cdot(\XX)$ at arithmetic infinity, 
with the operator $\Phi$ and the action of $\SL(2,\R)$ 
related to the Lefschetz operator, as a ``derived'' (cohomological)
version of a {\em spectral triple} \`a la Connes. 

More precisely, we prove that the bigraded polarized Lefschetz 
module structure on the complex $(K^\cdot, d=d'+d'')$ defines data
$({\rm A}, H^\cdot(\XX), \Phi)$, where the algebra ${\rm A}$ is obtained
from the action of the Lefschetz $\SL(2,\R)$ on the Hilbert space
completion of $H^\cdot(\XX)$ with respect to the  
inner product defined by the polarization on $K^\cdot$. The operator
$\Phi$ that determines the Archimedean factors of \cite{Den} satisfies
the properties of a Dirac operator. 

The data $({\rm A}, H^\cdot(\XX), \Phi)$ should be thought of as the
cohomological version of a more refined spectral triple, which encodes 
the full geometric data at arithmetic infinity in the structure of a
non-commutative manifold. The simplified cohomological information is
sufficient to the purpose of this paper, hence we leave the study of
the full structure to future work. 

The extra datum of the Schottky uniformization considered 
by Manin can be implemented in the data $({\rm A}, H^\cdot(\XX),
\Phi)$ by first associating to the Schottky group a pair of Fuchsian
Schottky groups in $\SL(2,\R)$ that correspond, via Bers' simultaneous 
uniformization, to a decomposition of $X_{\/\C}$ into two Riemann 
surfaces with boundary, and then making these groups act via the 
$\SL(2,\R)\times \SL(2,\R)$
representation of the Lefschetz module. One obtains this way a
non-commutative version of the handlebody ${\mathfrak X}_\Gamma$,
given by the group ring of $\Gamma$ acting via the representation of
$\SL(2,\R)$ associated to the Lefschetz operator.
The hyperbolic geometry is encoded in the Beltrami
differentials of Bers' simultaneous uniformization. In particular, we
show in \S \ref{unifsec} that, in the case of a real embedding
$\alpha:\K \hookrightarrow \C$, where the corresponding Riemann
surface is an orthosymmetric smooth real algebraic curve, the choice
of the Schottky group and of the quasi--circle giving the simultaneous
uniformization is determined canonically. 

This result allows us to reinterpret Deninger's regularized 
determinants describing the Archimedean factors in terms of an 
integration theory on the ``non-commutative manifold'' $({\rm 
A}, H^\cdot(\XX), \Phi)$. Theorem \ref{Gamma-zeta3} shows that 
the alternating product of the $\Gamma$-factors $L_\C 
(H^q(X_{/\C},\C),s)$ is recovered from a particular zeta function 
of the spectral triple. In \S \ref{ReiZsec} we interpret this 
alternating product as a Reidemeister torsion associated to the 
fiber at arithmetic infinity. 
 
The second part of the paper (Section \ref{6} and \ref{6bis}) 
concentrates on Manin's 
description of the dual graph ${\mathcal G}$ of the ``fiber at 
infinity'' of an arithmetic surface in terms of the infinite tangle of 
bounded geodesics in the hyperbolic handlebody ${\mathfrak 
X}_\Gamma$. 
 
More precisely, the suspension flow ${\mathcal S}_T$ of a dynamical 
system $T$ provides our model of the dual graph ${\mathcal G}$ of the 
fiber at infinity, which maps surjectively over the tangle of bounded 
geodesics considered in \cite{Man}. 
The map $T$ is a subshift of finite type which partially captures the 
dynamical properties of the action of the Schottky group on 
its limit set $\Lambda_\Gamma$. 
 
The first cohomology group of ${\mathcal S}_T$ is 
the ordered cohomology of the dynamical system $T$, in the sense of 
\cite{BoHa} \cite{PaTu} and it provides a model of the first 
cohomology of the dual graph of the fiber at infinity. The group 
$H^1({\mathcal S}_T)$ carries a natural filtration, 
which is related to the periodic orbits of the subshift of finite 
type. We give an explicit combinatorial description of homology and 
cohomology of ${\mathcal S}_T$ and of their pairing. 
 
We define a {\em dynamical cohomology} $H^1_{dyn}$ of the fiber at 
infinity as the graded space associated to the filtration of $H^1({\mathcal 
S}_T)$. Similarly, we introduce a {\em dynamical homology} $H_1^{dyn}$ as 
the sum of 
the spaces in the filtration of $H_1({\mathcal S}_T)$. These two 
graded spaces have an involution which plays a role analogous to the 
real Frobenius $\bar F_\infty$ on the cohomological theories of Section 
\ref{3}. 
 
Theorem \ref{map-ar-dyn} relates the dynamical cohomology to the 
Archimedean cohomology by 
showing that the Archimedean  cohomology sits as a particular subspace 
of the dynamical cohomology in a way that is compatible with the 
grading and with the action of the real Frobenius. The map that 
realizes this identification is obtained using the description of 
holomorphic differentials on the Riemann surface as Poincar\'e series 
over the Schottky group, which also plays a fundamental role in the 
description of the Green function in terms of geodesics in \cite{Man}. 
Similarly, in Theorem 
\ref{map-ar-dyn-2} we identify a 
subspace of the dynamical homology that is isomorphic to the image of 
the Archimedean  cohomology under the duality isomorphism acting on 
$\H^\cdot ({\rm Cone}^\cdot(N))$. This way we reinterpret this 
arithmetic duality as induced by the pairing of dynamical homology and 
cohomology. 
 
The Cuntz--Krieger algebra ${\mathcal O}_A$ associated to the subshift
of finite type $T$ (\cf \cite{Cu} \cite{CuKrie}) acts on the space
${\mathcal L}$  of cochains
defining the dynamical cohomology $H^1({\mathcal S}_T)$. This algebra
carries a refined information  
on the action of the Schottky group on its limit set. 

We introduce a Dirac operator $D$ on the Hilbert space of cochains
${\mathcal H}={\mathcal L} \oplus {\mathcal L}$, whose restriction to the 
subspaces isomorphic to the Archimedean  cohomology and its dual, 
recovers the Frobenius-type operator $\Phi$ of Section \ref{3}. 
We prove in Theorem \ref{dyn-SP3OA} that the data $({\mathcal O}_A,
{\mathcal H},D)$ define a spectral triple. In Proposition \ref{L-factor2} 
we show how to recover the local Euler factor from these 
data. 
 
In \S \ref{modinftysect} we describe the analog at arithmetic  
infinity of the $p$-adic reduction map considered in \cite{Man} and 
\cite{Mum} for Mumford curves, which is realized in terms of certain
finite graphs in a quotient of the Bruhat-Tits tree. The corresponding
object at arithmetic infinity is, as originally suggested in
\cite{Man}, constructed out of arcs geodesics in the handlebody
$\mX_\Gamma$ which have one end on the Riemann surface $X_{/\C}$ and
whose asymptotic behavior is prescribed by a limiting point on
$\Lambda_\Gamma$. The resulting space is a well known construction in
noncommutative geometry, namely the homotopy quotient $\Lambda_\Gamma
\times_\Gamma \H^3$ of the space ${\mathcal O}_A= {\rm
C}(\Lambda_\Gamma)\rtimes \Gamma$. Similarly, our geometric model
${\mathcal S}_T$ of the dual graph of the fiber at arithmetic infinity
is the homotopy quotient ${\mathcal S}\times_\Z \R$ of the
noncommutative space described by the corssed product algebra
${\rm C}({\mathcal S})\rtimes_T \Z$.

In the last section of the paper, we outline some possible further 
questions and directions for future investigations. 
 
\medskip 
 
Since this paper draws from the language and techniques of different 
fields (arithmetic geometry, non-commutative geometry, dynamical 
systems), we thought it necessary to include enough background 
material to make the paper sufficiently self contained and addressed 
to readers with different research interests.

\bigskip 
 
\noindent {\bf Acknowledgments.} In the course of this project we 
learned many things from different people to whom we are very grateful: 
to Paolo Aluffi for useful discussions and suggestions; to 
Alain Connes for beautiful lectures on spectral triples and 
noncommutative geometry that inspired the early development of this
work and for many discussions, comments and suggestions that greatly
improved the final version of the paper; to Curt McMullen for very 
enlightening conversations on the dynamics of the shift $T$ and on 
Schottky groups; to Victor Nistor for various remarks on cross product 
algebras; to Mika Sepp\"al\"a for useful information on real algebraic 
curves; and of course to Yuri Manin for sharing his vision and 
enlightenment on many aspects of this project. 
 
\subsection{Preliminary notions and notation} 
 
The three-dimensional real hyperbolic space $\H^3$ is the quotient 
\begin{equation}\label{1quotient} 
\H^3 = {\rm PGL}(2,\C) /{\rm SU}(2). 
\end{equation} 
It can also be described as the upper half space $\H^3 \simeq 
\C\times\R^+$ endowed with the hyperbolic metric. 
 
The group $\PSL(2,\C)$ is the group of orientation preserving 
isometries of $\H^3$. The action is given by 
\begin{equation}\label{action-H3} 
 \gamma : (z,y) \mapsto \left(\frac{(az+b)\overline{(cz+d)} + a\bar 
c y^2}{|cz+d|^2 + |c|^2 y^2}, \frac{y \, |ad-bc|}{|cz+d|^2 + |c|^2 
y^2}\right), 
\end{equation} 
for $(z,y)\in \C\times\R^+$ and 
$$ \gamma =\left(\begin{array}{cc} a& b \\ c& d \end{array}\right)\in 
\SL(2,\C). $$ 
 
The complex projective line $\P^1(\C)$ can be identified with the 
conformal boundary at infinity of $\H^3$. The action \eqref{action-H3} 
extends to an action on $\overline \H^3:=\H^3\cup \P^1(\C)$, where 
$\PSL(2,\C)$ acts on $\P^1(\C)$ by fractional linear transformations 
$$ \gamma : z \mapsto \frac{(az+b)}{(cz+d)}. $$ 
 
\medskip 
 
We begin by recalling some classical facts about Kleinian and Fuchsian 
groups (\cf \cite{Beardon} \cite{Bo} \cite{MaTa}). 
 
\smallskip 
 
A {\em Fuchsian group} $G$ is a discrete subgroup of $\PSL(2,\R)$, the group 
of orientation preserving isometries of the hyperbolic plane $\H^2$. A 
{\em Kleinian group} is a discrete subgroup of $\PSL(2,\C)$, the group of 
orientation preserving isometries of three-dimensional real hyperbolic 
space $\H^3$. 
 
\smallskip 
 
For $g\ge 1$, a {\em Schottky group} of rank $g$ is a discrete 
subgroup $\Gamma\subset \PSL (2,\C)$, which is purely loxodromic 
and isomorphic to a free group of rank $g$. Schottky groups are 
particular examples of Kleinian groups. 
 
\smallskip 
 
A Schottky group that is specified by real parameters 
so that it lies in $\PSL(2,\R)$ is called a {\em Fuchsian Schottky 
group}.  Viewed as a group of isometries of the hyperbolic plane 
$\H^2$, or equivalently of the Poincar\'e disk, a Fuchsian Schottky 
group $G$ produces a quotient $G\backslash \H^2$ which is 
topologically a Riemann surface with boundary. 
 
\medskip 
 
In the case $g=1$, the choice of a Schottky group $\Gamma\subset 
\PSL(2,\C)$ amounts to the 
choice of an element $q\in\C^*$, $|q|<1$. This acts on $\H^3$ by 
$$ \left(\begin{array}{cc} q^{1/2} & 0 \\ 0 & q^{-1/2} \end{array} 
\right) (z,y) = (qz, |q|y). $$ 
One sees that ${\mathfrak X} = \H^3/(q^\Z)$ is a solid torus with the 
elliptic curve $X_{/\C} = 
\C^*/(q^\Z)$ as its boundary at infinity. This space is known in the 
theory of quantum gravity as Euclidean BTZ black hole \cite{ManMar2}. 
 
In general, for $g\ge 1$, the quotient space 
\begin{equation}\label{handleb} 
{\mathfrak X}_\Gamma := \Gamma \backslash \H^3 
\end{equation} 
is topologically a handlebody of genus $g$. These also form an 
interesting class of Euclidean black holes (\cf \cite{ManMar2}). 
 
\medskip 
 
We denote by $\Lambda_\Gamma$, the {\em limit set} of the action 
of $\Gamma$. This is the smallest non--empty closed 
$\Gamma$--invariant subset of $\H^3\cup \P^1(\C)$. Since $\Gamma$ 
acts freely and properly discontinuously on $\H^3$, the set 
$\Lambda_\Gamma$ is contained in the sphere at infinity 
$\P^1(\C)$. This set can also be described as the 
closure of the set of the attractive and repelling fixed points 
$z^{\pm}(g)$ of the loxodromic elements $g\in \Gamma$. In the case 
$g=1$ the limit set consists of two points, but for $g\geq 2$ the 
limit set is usually a fractal of some Hausdorff dimension $0\leq \delta_H  
=\dim_H (\Lambda_\Gamma) < 2$. 
 
We denote by $\Omega_\Gamma$ the {\em domain of 
discontinuity} of $\Gamma$, that is, the complement of 
$\Lambda_\Gamma$ in $\P^1(\C)$. The quotient 
\begin{equation} 
X_{/\C} = \Gamma \backslash \Omega_\Gamma \label{RS} 
\end{equation} 
is a Riemann surface of genus $g$ and the covering 
$\Omega_\Gamma \to X_{/\C}$ is called a {\em Schottky uniformization} 
of $X_{/\C}$. 
Every complex Riemann surface $X_{/\C}$ admits a 
Schottky uniformization. 
 
\smallskip 
 
The handlebody \eqref{handleb} can be compactified by adding the 
conformal boundary at infinity $X_{/\C}$ to obtain 
\begin{equation}\label{handleb-comp} 
\overline{\mathfrak X}_\Gamma := {\mathfrak X}_\Gamma\cup X_{/\C}= 
\Gamma \backslash (\H^3 \cup \Omega_\Gamma). 
\end{equation} 
 
\medskip 
 
Let $\{ g_i \}_{i=1}^g$ be a set of 
generators of the Schottky group $\Gamma$. We write 
$g_{i+g}=g_i^{-1}$. There are $2g$ Jordan 
curves $\gamma_k$ on the sphere at 
infinity $\P^1(\C)$, with pairwise disjoint interiors $D_k$, 
such that the elements $g_k$ are given by fractional 
linear transformations that map the interior of $\gamma_k$ 
to the exterior of $\gamma_{j}$ with $|k-j|=g$. The curves $\gamma_k$ 
give a {\em marking} of the Schottky group. 
 
The choice of a Schottky uniformization for the Riemann surface 
$X_{/\C}$ provides a choice of a set of generators $a_i$, $i=1,\ldots 
g$, for $Ker(I_*)$, where $I_*: H_1(X_{/\C},\Z) \to H_1(\mX_\Gamma, 
\Z)$ is the map induce by the inclusion of $X_{/\C}$ in 
$\overline{\mX}_\Gamma$ as the conformal boundary at infinity. The 
$a_i$ are the images under the quotient map 
$\Omega_\Gamma \to X_{/\C}$ of the Jordan curves $\gamma_i$. 
 
\bigskip 
 
Recall that, if $\K$ is a number field with $n=[\K:\Q]$, there are $n$ 
Archimedean  primes which correspond to the embeddings $\alpha: \K 
\hookrightarrow \C$. Among these $n$ Archimedean  primes, there are $r$ 
embeddings into $\R$, and $s$ pairs of conjugate embeddings in $\C$ 
not contained in $\R$, so that $n=r+2s$. 
 
If $X$ is an arithmetic surface over ${\rm Spec}(O_\K)$, then at each 
Archimedean prime we obtain a Riemann surface $X_{/\C}$. If the 
Archimedean prime corresponds to a real embedding, the 
corresponding Riemann surface has a real structure, namely it is a 
smooth real algebraic curve $X_{/\R}$. 
 
\medskip 
 
A smooth real algebraic curve $X_{/\R}$ is a Riemann 
surface $X_{/\R}$ together with an involution $\iota: X_{/\R} \to 
X_{/\R}$ induced by complex conjugation $z\mapsto \bar z$. The fixed 
point set $X_{\iota}$ of the involution is the set of real points 
$X_{\iota}=X_{/\R}(\R)$ of $X_{/\R}$. If $X_{\iota} \neq \emptyset$, 
the components of $X_{\iota}$ are simple closed geodesics on $X_{/\R}$. 
A smooth real algebraic curve is called {\em orthosymmetric} if 
$X_{\iota} \neq \emptyset$ and the complement $X_{/\R} \backslash 
X_\iota$ consists of two connected components. If $X_{\iota} \neq 
\emptyset$, then $X_{/\R}$ can always be reduced to the orthosymmetric 
case upon passing to a double cover. 
 
\bigskip 
 
Even when not explicitly stated, all Hilbert spaces and algebras of 
operators we consider will be {\em separable}, \ie they admit a dense 
(in the norm topology) countable subset. 
 
\medskip 
 
An {\em involutive algebra} is an algebra over $\C$ with a 
conjugate linear involution $*$ (the adjoint) which is an 
anti-isomorphism. A {\em ${\rm C}^*$-algebra} is an involutive normed 
algebra, which is complete in the norm, and satisfies $\| a b \| \leq 
\| a \| \cdot \| b \|$ and $\| a^* a \| = \| a \|^2$. The analogous 
notions can be defined for algebras over $\R$. 
 
\medskip

\section{Cohomological Constructions.} \label{3}

In this chapter we give an explicit description of 
a cohomological theory for the Archimedean fiber of an Arakelov 
surface. The general theory, valid for any arithmetic variety, was 
defined in \cite{KC}. This construction 
provides an alternative definition and a refinement for the {\it Archimedean 
cohomology} $H^\ast_{\rm{ar}}$ introduced by Deninger in
\cite{Den}. The spaces  
$H^\cdot(\X)$ (\cf~Definition~\ref{gr-groups}) are 
{\it infinite dimensional} real vector spaces endowed with a 
monodromy operator $N$ and an endomorphism $\Phi$ 
(\cf~Section~\ref{2.5}). The groups $H^\ast_{\rm{ar}}$ 
can be identified with the subspace of the $N$-invariants (\ie 
${\rm Ker}(N)$) over which (the restriction of) $\Phi$ 
acts in the following way. The monodromy operator determines an 
integer, even 
graduation on $H^\cdot(\X) = \oplus_{p\in\Z}gr^\w_{2p}H^\cdot(\X)$ where 
each graded piece is still {\it infinite dimensional}. We 
will refer to it as to the {\it weight graduation}. This graduation 
induces a corresponding one on the subspace $H^\cdot(\X)^{N=0} := 
\oplus_{\cdot\ge 2p}gr^\w_{2p}H^\cdot(\XX)$. The summands 
$gr^\w_{2p}H^\cdot(\XX)$ are {\it finite dimensional} real vector spaces on 
which $\Phi$ acts as a multiplication by the weight $p$. 
 
When $X_{/\kappa}$ is a non-singular, projective curve defined over
$\kappa =  \C$ or $\R$, the 
description of $gr^\w_{2p}H^\cdot(\XX)$ ($\cdot\ge 2p$)  is particularly 
easy.  Proposition~\ref{grad} shows that $H^\cdot(\X)^{N=0}$ is 
isomorphic to an infinite direct sum of Hodge-Tate twisted copies of 
the same finite-dimensional vector space. For $\kappa = \C$, this 
space coincides with the de Rham 
cohomology $H^\ast_{DR}(X_{/\C},\R)$ of the Riemann surface $X_{/\C}$. 
 
For the reader acquainted with the classical theory of mixed 
Hodge structures 
for an algebraic degeneration over a disk (and its arithmetical 
counterpart theory of Frobenius weights), 
it will be immediately evident that the construction of the arithmetical 
cohomology defined in this chapter runs in parallel with the classical 
one defined by Steenbrink in \cite{St} and refined by M.~Saito in 
\cite{Saito}. The notation: $H^\cdot(\X)$, $H^\cdot(\XX)$, $H^\cdot(Y)$ 
followed in this section is 
purely formal. Namely, $\X$, $\XX$ and $Y$ are only symbols although 
this choice is motivated by the analogy with 
Steenbrink's construction in which $\X$, $\XX$ and $Y$ describe \resp the 
geometric generic fiber and the complement of the special fiber $Y$ in the 
model. The space $H^\cdot(\X)$ is the hypercohomology group of a
double complex $K^{\cdot,\cdot}$ of real, 
differential twisted forms (\cf Section~\ref{2.1}: \eqref{complex}) on 
which one defines a 
structure of polarized Lefschetz module that descends to its 
hypercohomology (\cf~Theorem~\ref{main thm} and 
Corollary~\ref{main cor}). 
 
The whole theory is inspired by the 
expectation that the fibers at infinity of an arithmetic variety should 
be thought to be semi-stable and more specifically to be `maximally 
degenerate or totally split'. 
We like to think that the construction of the 
complex $K^{\cdot,\cdot}$ on the Riemann surface $X_{/\kappa}$, whose 
structure and behavior 
gives the arithmetical information related to the `mysterious' fibers at 
infinity of an arithmetic surface, fits in with Arakelov's intuition 
that  Hermitian geometry on $X_{/\kappa}$ is enough 
to recover the intersection geometry on the fibers at infinity. 
 
\subsection{A bigraded complex with monodromy and Lefschetz 
operators}\label{2.1} 
 
Let $X_{/\kappa}$ be a smooth, projective curve defined over $\kappa = 
\C$ or $\R$. For 
$a,b\in\N$, we shall 
denote by $(A^{a,b}\oplus A^{b,a})_\R$ the abelian group of real 
differential forms (analytic or ${\rm C}^\infty$) on $X_{/\kappa}$ of type 
$(a,b) + (b,a)$. 
 
For $p\in\Z$, the expression $(A^{a,b}\oplus A^{b,a})_\R(p)$ means 
the $p$-th Hodge-Tate twist of $(A^{a,b}\oplus A^{b,a})_\R$, \ie 
\[ 
(A^{a,b}\oplus A^{b,a})_\R(p):= (2\pi\sqrt{-1})^p(A^{a,b}\oplus 
A^{b,a})_\R. 
\] 
 
Let $i,j,k\in\Z$. We consider the following complex 
(\cf\cite{KC},~\S 4 for the general construction) 
\begin{equation}\label{complex} 
K^{i,j,k} = 
\begin{cases} 
{\displaystyle \bigoplus_{\begin{subarray}{l} a+b=j+1\\|a-b|\le 
2k-i\end{subarray}} (A^{a,b}\oplus A^{b,a})_\R(\frac{1+j-i}{2})} 
&\text{if $1+j-i \equiv~0(2),~k\ge\text{max}(0,i)$} \\ 
0 &\text{otherwise.} 
\end{cases} 
\end{equation} 
 
On the complex $K^{i,j,k}$ one defines the following differentials 
\begin{align*} 
\,& d': K^{i,j,k} \to K^{i+1,j+1,k+1};\quad& \,& d'': K^{i,j,k} \to 
K^{i+1,j+1,k} \\ 
\,& d'= \partial+\overline\partial \quad& \,& d''= P^\perp
\sqrt{-1}(\overline\partial-\partial), 
\end{align*} 
with $P^\perp$ the orthogonal projection onto $K^{i+1,j+1,k}$.
These maps satisfy the property that ${d'}^2 = 0 = {d''}^2$ (\cf\cite{KC} 
Lemma~4.2). Since $X_{/\kappa}$ is a projective variety (hence 
K\"ahler) one uses the 
existence of the fundamental real (closed) $(1,1)$-form 
$\omega$ to define the 
following {\it Lefschetz map} $l$. The operator $N$ that is described 
in the next formula plays the role of the logarithm of the {\it local 
monodromy at infinity} 
\begin{equation}\label{monod} 
N: K^{i,j,k} \to K^{i+2,j,k+1},\quad N(f) = 
(2\pi\sqrt{-1})^{-1}f 
\end{equation} 
\begin{equation}\label{lef} 
l: K^{i,j,k} \to K^{i,j+2,k},\quad l(f) = (2\pi\sqrt{-1})f\wedge\omega 
\end{equation} 
 
These endomorphisms are known to commute with $d'$ and $d''$ and 
satisfy $[l,N] = 0$ (\cf\op 
Lemma~4.2). One sets $K^{i,j} = \oplus_k K^{i,j,k}$ and writes 
$K^\ast = \oplus_{i+j=\ast}K^{i,j}$ to denote the simple complex 
endowed with the total differential $d = d'+d''$ and with the action 
of the operators $N$ and $l$. 
 
\begin{rem}\label{rem1} {\rm In the complex \eqref{complex} 
the second index $j$ is subject to the constraint $a+b = j+1$ (where 
$a+b$ is the total degree of the differential forms). This implies 
that $j$ assumes only a finite number of values: $-1\le j \le 
1$, in fact $0\le a+b\le 2$ ($X$ is a Riemann surface). } 
\end{rem} 
 
\subsection{Polarized Hodge--Lefschetz structure}\label{2.2} 
 
In this paragraph we will review the theory of polarized bigraded 
Hodge-Lefschetz modules due to Deligne and Saito. The main result 
is Theorem~\ref{main thm} which states that the complex 
$K^{\cdot,\cdot}$ defined in \eqref{complex} together with the maps 
$N$ and $l$ as in \eqref{monod} and \eqref{lef} determine a 
Lefschetz module. A detailed 
description of the structure of polarized Hodge-Lefschetz modules is 
contained in \cite{Saito}; for a short and quite pleasant exposition 
we refer to  \cite{GNA}. 
 
\begin{defn}\label{def1} 
A bigraded Lefschetz module $(K^{\cdot,\cdot},L_1,L_2)$ is a bigraded 
real vector space 
$K=\oplus_{i,j} K^{i,j}$ with endomorphisms 
\begin{equation} \label{Li}  L_1: K^{i,j} \to K^{i+2,j}  \ \ \ \ \ \ 
L_2: K^{i,j} \to K^{i,j+2} 
\end{equation} 
satisfying $[ L_1,L_2 ]=0$. Furthermore, the operators $L_i$ are 
required to satisfy the following conditions 
\begin{enumerate} 
\item $L_1^i : K^{-i,j}\to K^{i,j}$ is an isomorphism for $i>0$ 
\item $L_2^j: K^{i,-j}\to K^{i,j}$ is an isomorphism for $j>0$. 
\end{enumerate} 
\end{defn} 
 
Bigraded Lefschetz modules correspond to 
representations of the Lie group $\SL(2,\R)\times \SL(2,\R)$ 
(\cf\cite{GNA} or \cite{Saito}). Therefore, given a 
bigraded Lefschetz module $(K^{\cdot,\cdot},L_1,L_2)$ this corresponds 
to the representation 
$$ \sigma: \SL(2,\R)\times \SL(2,\R) \to \Aut (K^{\cdot,\cdot}) $$ 
satisfying 
\begin{equation} \label{rep1} 
\sigma\left\{ \left(\begin{array}{cc} a & 0 \\ 0 & a^{-1} 
\end{array} \right) , \left(\begin{array}{cc} b & 0 \\ 0 & b^{-1} 
\end{array} \right) \right\} (x) = a^i b^j x \ \ \ \text{ for $x\in 
K^{i,j}$} \end{equation} 
\begin{equation} \label{rep2} 
d\sigma \left\{ \left( \begin{array}{cc} 0 & 1 \\ 0 & 0 \end{array} \right), 
0 \right\} = L_1 \end{equation} 
\begin{equation} \label{rep3} 
d\sigma \left\{ 0, \left( \begin{array}{cc} 0 & 1 \\ 0 & 0 \end{array} 
\right) \right\} = L_2. \end{equation} 
 
\medskip 
 
The Weyl reflection 
$$ w= \left( \begin{array}{cc} 0 & 1 \\ -1 & 0 \end{array} \right)\in 
\SL(2,\R)$$ 
defines the elements $\tilde w=\{w,w\},~w_1 = \{w,1\},~w_2 = 
\{1,w\} \in \SL(2,\R) \times \SL(2,\R)$. They 
determine isomorphisms $\sigma(\tilde w): K^{i,j}\to K^{-i,-j}$, 
$\sigma(w_1): K^{i,j}\to K^{-i,j}$ and $\sigma(w_2): K^{i,j}\to 
K^{i,-j}$, by taking $\sigma(w_1)=N^{-i}$ and
$\sigma(w_2)$ is the involution determined by the Hodge $*$ operator,
which induces the map $l^{-j}$ on the primitive cohomology (\cf
Definition \ref{def1} and \cite{Wells}, \S V.6).  
 
\medskip 
 
\begin{defn} A bigraded Lefschetz module is a Hodge--Lefschetz module if 
each $K^{i,j}$ carries a pure real Hodge structure and the $L_i$ (as 
in Definition~\ref{def1}) are morphisms of real Hodge structures. 
\end{defn} 
 
For convenience, we recall the definition of a {\it pure Hodge 
structure} over $\kappa = \C$ or $\R$. For a summary of mixed Hodge theory 
we refer to \cite{St1}. 
\begin{defn}\label{PHS} A pure Hodge structure over $\kappa$ is a finite 
dimensional $\C$-vector space $H = \oplus_{p,q}H^{p,q}$, together with 
a conjugate linear involution $c$ and in case $\kappa = \R$ a $\C$-linear 
involution $F_{\infty}$ such that 
\begin{enumerate} 
\item $c(H^{p,q}) = H^{q,p}$ 
\item the inclusion of $H_{\R} := H^{c=id}$ into $H$ induces an 
isomorphism $H = H_{\R}\otimes_{\R}\C$ 
\item in case $\kappa = \R$, $F_{\infty}$ commutes with $c$ and 
verifies $F_{\infty}(H^{p,q}) = H^{q,p}$. The action of 
$F_{\infty}$ on the space $H^{p,p}$ decomposes it as $H^{p,p} = 
H^{p,+}\oplus H^{p,-}$. We denote the dimensions of the eigenspaces by 
$h^{p,\pm} := \dim_{\C}H^{p,\pm(-1)^p}$. 
 
In the case $\kappa = \R$, $H$ is called a real, pure Hodge structure. 
\end{enumerate} 
\end{defn} 
 
\begin{ex}\label{ex} {\rm An example of pure Hodge structure 
    is given by the singular (Betti) cohomology $H^\ast_B(X_{/\C},\C)$ 
    on the Riemann surface $X_{/\C}$. The $\C$-linear involution 
    $F_{\infty}$ is induced by the complex conjugation on the Riemann 
    surface.}\end{ex} 
 
On a bigraded Lefschetz module $(K^{\cdot,\cdot},L_1,L_2)$ we consider 
the additional data of a differential $d$ and a pairing $\psi$: 
\[ 
d: K^{i,j}\to K^{i+1,j+1}, \ \ \ \ \psi : K^{-i,-j}\otimes K^{i,j} \to 
\R(1), 
\] 
satisfying the following properties: 
\begin{enumerate} 
\item $d^2 = 0 = [L_i,d]$ 
\item $\psi(x,y) = - \psi(y,x)$ 
\item $\psi(dx,y) = \psi(x,dy)$ 
\item $\psi(L_i x,y) + \psi(x,L_i y) = 0$ 
\item $\psi(\cdot, L_1^i L_2^j \cdot)$ is symmetric 
and positive 
definite on $K^{-i,-j} \cap {\rm Ker}(L_1^{i+1}) \cap 
{\rm Ker}(L_2^{j+1})$. 
\end{enumerate} 
 
\smallskip 
 
If $(K,L_1,L_2,\psi)$ is a polarized bigraded Lefschetz module 
(\ie $(K,L_1,L_2)$ is a bigraded Lefschetz module satisfying the 
properties {\rm{1.}-\rm{5.}}), then 
the bilinear form 
\begin{equation}\label{innerprod} 
\langle\cdot,\cdot\rangle: K\otimes K \to \R(1), \ \ \ \ 
\langle x, y \rangle : = \psi\left( x, \sigma(\tilde w) y \right) 
\end{equation} 
is symmetric and positive definite. 
 
\smallskip 
 
\begin{thm}\label{main thm} The differential complex $K^{\cdot,\cdot}$ 
defined in \eqref{complex} endowed with the operators $L_1 = N$ 
(\cf\eqref{monod}) and $L_2 = l$ (\cf\eqref{lef}) is a polarized 
bigraded Lefschetz module. 
The polarization is given by 
\[ 
\psi : K^{-i,-j,k}\otimes K^{i,j,k+i} \to \R(1) 
\] 
\[ 
\psi(x,y) := \left( \frac{1}{2\pi \sqrt{-1}} \right) 
\epsilon(1-j) (-1)^k \int_{X(\C)} x \wedge Cy. 
\] 
Here, for $m\in\Z$: $\epsilon(m) := (-1)^{\frac{m(m+1)}{2}}$ and 
$C(x):= (\sqrt{-1})^{a-b}  x$ is the Weil operator, for $x$ 
a differential form of type $(a,b)$ (\cf\cite{Wells} \S V.1). 
\end{thm} 
\noindent {\bf Proof.} We refer to \cite{KC} Lemmas~4.2,~4.5,~4.6 and 
Proposition~4.7). 
\noindent $\diamond$ 
\bigskip 
 
Such elaborate construction on the complex $K^{\cdot,\cdot}$ allows one 
to set up a harmonic theory as in \cite{KC} pp.~350-1, so that the 
polarized bigraded Lefschetz module structure passes to the 
hypercohomology $\HH^\ast(K^\cdot,d)$. More precisely, one defines a 
Laplace operator on $K^{\cdot,\cdot}$ as 
\[ 
\Box := d(^{t}d) + (^{t}d)d 
\] 
where $^{t}d$ is the transpose of $d$ relative to the bilinear form 
$\langle\cdot,\cdot\rangle$ defined in \eqref{innerprod}. Then, 
$\Box$ commutes with the action of $\SL(2,\R) \times \SL(2,\R)$ 
(\cf\cite{Wells} Lemma at p.~153). Using the properties of the 
bilinear form $\langle\cdot,\cdot\rangle$ one gets 
\[ 
\HH^\ast(K^\cdot,d) = {\rm Ker}(d) \cap {\rm Ker}(^{t}d) = 
{\rm Ker}(\Box) 
\] 
and $\Box$ is invariant for the action of $\SL(2,\R) \times 
\SL(2,\R)$. The following result holds 
\begin{cor}\label{main cor} The data $(\HH^\ast(K^\cdot,d),N,l,\psi)$ define 
a polarized, bigraded Hodge-Lefschetz module. 
\end{cor} 
\noindent {\bf Proof.} The statement follows from the isomorphism of 
complexes 
\[ 
K^\cdot \simeq {\rm Ker}(\Box) \oplus {\rm Image}(\Box) 
\] 
and from the facts that $d = 0$ on ${\rm Ker}(\Box)$ and that the 
complex ${\rm Image}(\Box)$ is $d$-acyclic. These three statements 
taken together imply the existence of an induced action of $\SL(2,\R) 
\times \SL(2,\R)$ on the hypercohomology of $K^{\cdot}$ (\cf\cite{GNA} 
for details). 
\noindent $\diamond$

\subsection{Cohomology groups}\label{2.3} 
 
It follows from the definition of the double complex $(K^{\cdot,\cdot}, 
d',d'')$ in \eqref{complex} that the total differential $d = d' + d''$ 
satisfies $d^2 = 0$ and commutes with the operator $N$. In particular, 
$d$ induces a differential on the graded groups 
\[ 
{\rm Ker}(N)^{\cdot,\cdot} = {\rm ker}(N:K^{\cdot,\cdot} \to K^{\cdot + 
2,\cdot}),\ \ \ \ {\rm Coker}(N)^{\cdot,\cdot} = {\rm coker}(N: 
K^{\cdot,\cdot} \to K^{\cdot +2,\cdot}) 
\] 
as well as on the mapping cone of $N$ 
\[ 
{\rm Cone}(N)^{\cdot,\cdot} = {\rm Cone}(N: K^{\cdot,\cdot}\to K^{\cdot 
+2, \cdot}) := K^{\cdot,\cdot}[1] \oplus K^{\cdot +2,\cdot}, \ \ \
D(a,b) = (-d(a),N(a)+d(b))  
\] 
 
\begin{defn}\label{gr-groups} For any non-negative integer $q$ and 
$p\in\Z$, define 
 
\begin{align}\label{X-co} 
gr^\w_{2p}H^q(\X) &= \frac{{\rm Ker}(d:K^{q-2p,q-1}\to 
K^{q-2p+1,q})}{{\rm Im}(d: K^{q-2p-1,q-2}\to 
K^{q-2p,q-1})}, \\[2mm] 
gr^\w_{2p}H^q(Y) &= \frac{{\rm Ker}(d:{\rm Ker}(N)^{q-2p,q-1}\to 
{\rm Ker}(N)^{q-2p+1,q})}{{\rm Im}(d: 
{\rm Ker}(N)^{q-2p-1,q-2}\to {\rm Ker}(N)^{q-2p,q-1})}, \label{Y-co} \\[2mm] 
gr^\w_{2p}H^q_Y(X) &= \frac{{\rm Ker}(d:{\rm Coker}(N)^{q-2p,q-3}\to 
{\rm Coker}(N)^{q-2p+1,q-2})}{{\rm Im}(d: 
{\rm Coker}(N)^{q-2p-1,q-4}\to {\rm Coker}(N)^{q-2p,q-3})},
\label{Y-X-co} \\[2mm]  
gr^\w_{2p}H^q(\XX) &=\frac{{\rm Ker}(d:{\rm Cone}(N)^{q-2p+1,q-2}\to 
{\rm Cone}(N)^{q-2p+2,q-1})}{{\rm Im}(d: 
{\rm Cone}(N)^{q-2p,q-3}\to {\rm Cone}(N)^{q-2p+1,q-2})}. \label{X-star-co} 
\end{align} 
We define: $H^q(\X) := \HH^q(K^{\cdot})$, $H^q(Y) := 
\HH^q({\rm Ker}(N)^{\cdot})$, $H^q(\XX) := \HH^q({\rm Cone}(N)^\cdot)$ 
and $H^q_Y(X) := \HH^q({\rm Coker}(N)^\cdot)$. These groups are 
identified with 
\begin{equation*} 
\left. 
\begin{array}{cc} 
H^q(\X) = \bigoplus_{p\in\Z}gr^\w_{2p}H^q(\X), & H^q(Y) = 
\bigoplus_{p\in\Z}gr^\w_{2p}H^q(Y), \\ \\ 
H^q(\XX) = \bigoplus_{p\in\Z}gr^\w_{2p}H^q(\XX),  & H^q_Y(X) = 
\bigoplus_{p\in\Z}gr^\w_{2p}H^q_Y(X). 
\end{array} 
\right. 
\end{equation*} 
\end{defn}\vspace{.1in} 
 
\begin{rem}\label{rem2} {\rm Note that the even graduation is a 
consequence of the 
parity condition $q+r \equiv~0~(2)$ imposed on the indices of the 
complex \eqref{complex}.}\end{rem} 
 
Because $\dim X_{/\kappa} = 1$, it is easy to verify from the definition of 
$K^{\cdot,\cdot}$ that $H^q(\X)$ and $H^q(Y)$ are $\neq 0$ only for 
$q = 0,1,2$. Furthermore, one easily finds that $H^q(\XX) \neq 0$ for $q 
= 0,1,2,3$ and $H^q_Y(X)\neq 0$ only when $q=2,3,4$.\vspace{.05in} 
 
The definition of these groups is inspired by the theory of 
degenerations of Hodge structures (\cf\cite{St}) where the symbols 
$\X$, $Y$ and $\XX$ have a precise geometric meaning: namely, they 
denote \resp the 
smooth fiber, the special fiber and the punctured space ${\mathcal X} - Y$, 
where ${\mathcal X}$ is the chosen model for a degeneration over a disk. 
In our set-up instead, $\X$, $Y$, and $\XX$ are 
only symbols but the general formalism associated to the 
hypercohomology of a double complex endowed with an operator commuting 
with the total differential can still be pursued and in fact it gives 
interesting arithmetical information. In the 
following we will show that the 
groups that we have just introduced enjoy similar properties as the 
graded quotients of the weight-filtration on the corresponding 
cohomology groups of \op.\medskip 
 
It is important to remark that the hypercohomology of the complex 
${\rm Cone}(N)^\cdot$ contains both the information coming from the 
cohomologies of ${\rm Ker}(N)^\cdot$ and ${\rm Coker}(N)^\cdot$, as the 
following proposition shows 
\begin{prop}\label{conedec} The following equality holds: 
\begin{multline*} 
H^q(\XX) = \bigoplus_{p\in\Z} gr^\w_{2p}H^q(\XX) = \\ 
= \bigoplus_{2p\le q-1} gr^\w_{2p}H^q(Y) \oplus gr^\w_{q}H^q(\XX) 
\oplus gr^\w_{q+1}H^q(\XX) 
\oplus\bigoplus_{2p > q+1}gr^\w_{2p}H^{q+1}_Y(X). 
\end{multline*} 
\end{prop}\medskip 
The proof of proposition~\ref{conedec} as well as an explicit 
description of each addendum in the sum is a consequence of the 
following lemmas 
\begin{lem}\label{ex-seqs} For all $p \in \Z$ and for $q \in \N$ there 
are exact sequences 
\begin{equation}\label{Wang} 
\cdots\to gr^w_{2p}H^q(\XX) \to gr^w_{2p}H^q(\X) 
\stackrel{\text{N}}{\to} gr^\w_{2p-2}H^q(\X) \to 
gr^\w_{2p}H^{q+1}(\XX) \to\cdots 
\end{equation} 
\begin{equation}\label{with-supp} 
\cdots\to gr^w_{2p}H^q(Y) \to gr^w_{2p}H^q(\XX) 
\to gr^\w_{2p}H^{q+1}_Y(X) \to 
gr^\w_{2p}H^{q+1}(Y) \to\cdots 
\end{equation} 
Furthermore, the maps $gr^\w_{2p}H^q_Y(X) \to gr^\w_{2p}H^q(Y)$ in
\eqref{with-supp} are zero unless $q = 2p$, in which case they  
coincide with the morphism 
\[ 
\frac{(A^{p-1,p-1})_\R(p-1)}{{\rm Im}(d'')} 
\stackrel{2\pi\sqrt{-1}d'd''}{\longrightarrow} 
{\rm Ker}\left( d': (A^{p,p})_\R(p) \to (
\bigoplus_{\begin{subarray}{l}a+b =  
2p+1 \\|a-b|\le 1\end{subarray}}A^{a,b})_\R (p)\right). 
\] 
\end{lem} 
\noindent {\bf Proof.} We refer to \cite{BGS}: Lemmas~3 and Lemma~4 
and to \cite{KC}: Lemma~4.3. 
\noindent $\diamond$\medskip 
\begin{lem}\label{descs} For all $p \in \Z$ and for $q \in 
\N$:\begin{enumerate} 
\item The group $gr^\w_{2p}H^q(Y)$ is zero unless $2p \le q$, in which 
case 
{\small \[
gr^w_{2p}H^q(Y) = \begin{cases} 
{\displaystyle {\rm Ker}(d': 
(A^{p,p})_\R(p)\to \bigoplus_{\begin{subarray}{l}a+b = 2p+1 \\|a-b|\le 
1\end{subarray}}(A^{a,b})_\R (p))} &\text{if $q=2p$} \\ \\ 
{\frac{\displaystyle 
{\rm Ker}(d':(\bigoplus_{\begin{subarray}{l}a+b = 
q \\|a-b|\le q-2p\end{subarray}}A^{a,b})_\R (p)  \to 
(\bigoplus_{\begin{subarray}{l}a+b = q+1 \\|a-b|\le 
q-2p+1\end{subarray}}A^{a,b})_\R (p))}{\displaystyle{\rm Im}(d')}} & 
\text{if $q \ge 2p+1.$} 
\end{cases} 
\] }
\item The group $gr^\w_{2p}H^q_Y(X)$ is zero unless $2p \ge q$, in 
which case 
{\small \[ 
gr^\w_{2p}H^q_Y(X) = \begin{cases} 
{\displaystyle {\rm Coker}(d'': 
(\bigoplus_{\begin{subarray}{l}a+b = 2p-3 \\|a-b|\le 
1\end{subarray}}A^{a,b})_\R (p-1)\to (A^{p-1,p-1})_\R (p-1))} 
&\text{if $q=2p$} \\ \\ 
{\frac{\displaystyle {\rm Ker}(d'':(\bigoplus_{\begin{subarray}{l}a+b = 
q-2 \\|a-b|\le 2p-q\end{subarray}}A^{a,b})_\R (p-1)  \to 
(\bigoplus_{\begin{subarray}{l}a+b = q-1 \\|a-b|\le 
2p-q-1\end{subarray}}A^{a,b})_\R (p-1))}{\displaystyle{\rm Im}(d'')}} & 
\text{if $q \le 2p-1.$} 
\end{cases} 
\] }
\end{enumerate} 
\end{lem} 
\noindent {\bf Proof.} We refer to \cite{KC} Lemmas~4.3. 
\noindent $\diamond$ \medskip 
 
\noindent {\bf Proof.}(of Prop.~\ref{conedec}) It is a straightforward 
consequence of Lemma~\ref{ex-seqs}. \noindent$\diamond$\medskip 
 
\begin{cor}\label{behv-N} The monodromy map 
\[ 
N: gr^\w_{2p}H^q(\X) \to gr^\w_{2(p-1)}H^q(\X)\ \ \ \text{is} \ \ 
\left\{\begin{array}{ccc} \text{injective} &\text{if $q<2p-1$} \\ 
\text{bijective} &\text{if $q = 2p-1$} \\ \text{surjective} &\text{if $q 
\ge 2p$.}\end{array} \right. 
\] 
Therefore, the following sequences are exact 
\begin{equation}\label{sequ1-ex} 
q\ge 2p:\quad 0 \to gr^\w_{2p}H^{q}(\XX) \to gr^\w_{2p}H^{q}(\X) 
\stackrel{N}{\to} gr^\w_{2(p-1)}H^{q}(\X) \to 0, 
\end{equation} 
\begin{equation}\label{sequ2-ex} 
q\le 2(p-1):\quad 0 \to gr^\w_{2p}H^{q}(\X) \stackrel{N}{\to} 
gr^\w_{2(p-1)}H^{q}(\X) \to gr^\w_{2p}H^{q+1}(\XX) \to 0. 
\end{equation} 
 
In particular, one obtains 
\[ 
{\rm Ker}(N) = \left\{\begin{array}{cc} gr^\w_{2p}H^q(\XX) &\text{if 
$q\ge 2p$} \\ \\ 
0 &\text{if $q \le 2p-1$.} \end{array} \right. 
\] 
\end{cor} 
\noindent {\bf Proof.} This follows from Lemma~\ref{ex-seqs}: 
\eqref{Wang} and from Lemma~\ref{descs}.\noindent $\diamond$ \medskip 
 
\begin{rem}\label{rem3}{\rm For future use, we explicitly remark 
    that when $q\ge 2p$ one has the following decomposition: 
\begin{align} 
gr^\w_{2p}H^q(\X) &= \operatorname*{\oplus}_{k\ge 
  q-2p}\frac{{\rm Ker}(d: K^{q-2p,q-1,k}\to\cdots)}{{\rm Im}(d)} = \\ 
  &= gr^\w_{2p}H^q(\XX) \operatorname*{\oplus}_{k\ge 
  q-2p+1}\frac{{\rm Ker}(d: K^{q-2p,q-1,k}\to\cdots)}{{\rm Im}(d)}.\nonumber 
\end{align} 
Hence, when $q = 0,1,2$, the group $gr^\w_{2p}H^q(\XX)$ coincides with 
the homology of the 
complex $$(\oplus_{k\ge q-2p}K^{q-2p,q-1,k}, d)$$ at $k = q-2p \ge 0$. 
}\end{rem}\medskip 
 
It is important to recall that the presence of a structure of 
polarized Lefschetz module on the hypercohomology $\HH^\ast(K^\cdot, 
d) = H^\ast(\X)$ allows one to state the following results 
 
\begin{prop}\label{N-isos} For $q,p\in\Z$ satisfying the conditions 
$q - 2p > 0$, $q\ge 0$, the operator $N$ induces 
isomorphisms 
\begin{equation}\label{Nisos} 
N^{q-2p}: gr^\w_{2(q-p)}H^q(\X) \stackrel{\simeq}{\to} 
gr^\w_{2p}H^q(\X). 
\end{equation} 
Furthermore, for $q \ge 2p$ the isomorphisms \eqref{Nisos} induce 
corresponding isomorphisms 
\begin{equation}\label{Hiso} 
(gr^\w_{2p}H^q(\X)^{N=0} \simeq)~gr^\w_{2p}H^{q}(\XX) 
\stackrel{N^{2p-q}}{\operatorname*{\to}_\simeq} gr^\w_{2(q-p+1)}H^{q+1}(\XX). 
\end{equation} 
\end{prop} 
\noindent {\bf Proof.} For the proof of \eqref{Nisos} we refer to 
\cite{KC}: Proposition~4.8. For a proof of the 
isomorphisms \eqref{Hiso} we refer to Corollary~\ref{behv-N} and 
either the proof of Proposition~\ref{isos-del} or to \op.\noindent 
$\diamond$ \medskip

\subsection{Relation with Deligne cohomology}\label{2.4} 
 
The main feature of the complex \eqref{complex}  is its 
relation with the {\it real Deligne cohomology} of $X_{/\kappa}$. This 
cohomology (\cf~Definition~\ref{del-coh}) measures how the natural 
real structure on 
the singular cohomology of a smooth projective variety behaves with 
respect to the de Rham 
filtration. One of the most interesting properties of Deligne cohomology 
is its connection with arithmetics. 
Proposition~\ref{pole-ord} describes a precise 
relation between the ranks of some real Deligne cohomology groups 
and the orders 
of pole, at non-positive integers, of the $\Gamma$-factors 
attached to a (real) Hodge structure $H = \oplus_{p,q}H^{p,q}$ over 
$\C$ (or $\R$). We recall that these factors are defined as (\cf\cite{Se}) 
\begin{align}\label{factors} 
L_{\C}(H,s) &= \prod_{p,q}\Gamma_{\C}(s-{\rm min}(p,q))^{h^{p,q}} 
\\ L_{\R}(H,s) &= 
\prod_{p<q}\Gamma_{\C}(s-p)^{h^{p,q}}\prod_p\Gamma_{\R}(s-p)^{h^{p,+}}
\Gamma_{\R}(s-p+1)^{h^{p,-}}
\nonumber  
\end{align} 
where $s\in\C$, $h^{p,q} = \dim_{\C}H^{p,q}$ and $h^{p,\pm}$ is the 
dimension of the $\pm(-1)^p$-eigenspace of the $\C$-linear involution 
$F_{\infty}$ on $H$ (\cf~Definition~\ref{PHS}). One sets 
\begin{equation}\label{gammas} 
\Gamma_{\C}(s) = (2\pi)^{-s}\Gamma(s), \ \ \ \ \Gamma_{\R}(s) = 
2^{-\frac{1}{2}}\pi^{-\frac{s}{2}}\Gamma(\frac{s}{2}) \ \ \ \
\text{where} \ \ \   
\ \Gamma(s) = \int_0^\infty e^{-t}t^s\frac{dt}{t}. 
\end{equation} 
These satisfy the Legendre--Gauss duplication formula
\begin{equation}\label{dupl} 
\Gamma_{\C}(s) = \Gamma_{\R}(s)\Gamma_{\R}(s+1). 
\end{equation} 
 
The $\Gamma$-function $\Gamma(s)$ is analytic in the whole complex
plane except  
for simple poles at the non-positive integer points on the real axis. 
 
\begin{defn}\label{del-coh} The real Deligne cohomology $H^\ast_{\mathcal 
D}(X_{/\C},\R(p))$ of a projective, smooth variety $X_{/\C}$ is the 
cohomology of the complex 
\[ 
\R(p)_{\mathcal D}:\ \ \R(p) \to \mathcal O_{X(\C)} \to \Omega^1 \to 
\cdots \to \Omega^{p-1} \to 0. 
\] 
The first arrow is the inclusion of the constant sheaf of 
Hodge-Tate twisted 
constants into the structural sheaf of the manifold: $\R(p) := 
(2\pi\sqrt{-1})^p\R \subseteq \C \subseteq 
\mathcal O_{X(\C)}$. $\Omega^i$ denotes the sheaf of holomorphic 
$i$-th differential forms on the manifold. 
 
The real Deligne cohomology of $X_{/\R}$ is defined as the subspace of 
elements invariant with respect to the de Rham conjugation 
\[ 
H^\ast_{\mathcal D}(X_{/\R},\R(p)) := H^\ast_{\mathcal 
D}(X_{/\C},\R(p))^{\bar F_{\infty}}. 
\] 
Precisely: $\bar F_{\infty}$ corresponds to the de Rham conjugation under 
the canonical identification between de Rham and Betti cohomology 
$H^\ast_{DR}(X_{/\C},\C) = H^\ast_B(X_{/\C},\C)$ on the manifold $X_{/\C}$. 
\end{defn} 
 
\begin{rem}\label{rem4} {\rm To understand correctly the meaning of 
    $\bar F_\infty$ in the definition, it is worth recalling that both 
    de Rham and singular cohomology of the manifold $X_{/\C}$ have 
    real structures. The real structure on the singular cohomology 
    $H^\ast_B(X_{/\C},\R) \otimes_\R \C = H^\ast_B(X_{/\C},\C)$ is given 
    by the $\R$-linear involution:~~$\bar{}$~~(on the right hand side) 
    which is induced by complex conjugation on the 
    coefficients. On the other hand, by GAGA the algebraic de Rham 
    cohomology of $X_{/\R}$ defines a real structure 
    $H^\ast_{DR}(X_{/\R},\R) \otimes_\R \C = H^\ast_{DR}(X_{/\C},\C)$ on 
    the analytic de Rham cohomology. The complex conjugation on the 
    pair $(X_{/\C},\Omega^\cdot)$ induces a $\R$-linear involution on 
    the right hand side called de Rham conjugation. Such de Rham 
    conjugation corresponds to the operator $\bar F_\infty$ on the 
    Betti cohomology (\cf~example~\ref{ex}) through the identification 
    $H^\ast_{DR}(X_{/\C},\C) = H^\ast_B(X_{/\C},\C)$ on the manifold 
    $X_{/\C}$. 
 
For example, assume that $X_{/\R}$ is a smooth, real algebraic 
    curve, that is, a symmetric Riemann surface with an involution 
    $\iota: X_{/\R} \to X_{/\R}$ induced by the complex 
    conjugation. Then, $\iota$ determines a corresponding involution 
    $\iota^\ast$ in cohomology. On $H^0(X_{/\R},\R)$ this is the 
    identity, whereas on $H^2(X_{/\R},\R)$ it reverses the orientation, 
    hence is identified with $-{\rm id}$. On $H^1(X_{/\R},\R)$, 
    $\iota^\ast$ acts as a non-trivial involution that exchanges 
    $H^{1,0}_\R$ and $H^{0,1}_\R$, hence it separates the cohomology group 
    into two eigenspaces corresponding to the eigenvalues of $\pm 1$: 
    $H^1(X_{/\R},\R) = E_1 \oplus E_{-1}$, with $\dim E_1 = \dim 
    E_{-1} = g =$ genus of $X$. On twisted cohomology groups such as 
    $H^\ast(X_{/\R},\R(p))$, $\bar F_\infty$ acts as the composition of 
    $\iota^\ast$ with the involution that acts on the real Hodge 
    structure $\R(p) = (2\pi i)^p\R \subset \C$ as $(-1)^p$.}\end{rem} 
 
{}From the short exact sequence of complexes 
\[ 
0 \to \Omega^\cdot_{< p}[-1] \to \R(p)_{\mathcal D} \to \R(p) \to 0 
\] 
one gets the following long exact cohomology sequence that explains the 
statement we made at the beginning of this chapter 
($F^pH^\ast_{DR}(X_{/\C},\C) =$ the Hodge filtration on de Rham 
cohomology of the complex manifold) 
\[ 
\cdots\to H^\ast(X_{/\C},\R(p)) \to H^\ast_{DR}(X_{/\C},\C)/F^p \to H^{\ast 
+1}_{\mathcal D}(X_{/\C},\R(p)) \to H^{\ast +1}(X_{/\C},\R(p)) \to\cdots 
\] 
\begin{prop}\label{pole-ord} Let $X_{/\kappa}$ be a smooth, projective 
variety over $\kappa = \C$ or $\R$. Let $i,~m$ be two integers 
satisfying the conditions: $i \ge 0$,  $2m \le i$. Set $n = 
i+1-m$. Then 
\[ 
-\dim_{\R}H^{i+1}_{\mathcal D}(X_{/\kappa},\R(n)) = 
 \operatorname*{\mathrm{ord}}_{s=m} L_{\kappa}(H^i_B(X_{/\C},\C),s) 
\] 
Here $H^i_B(X_{/\C},\C)$ denotes Betti's cohomology of the manifold 
$X_{/\C}$ endowed with its pure Hodge structure over $\kappa$. 
\end{prop} 
\noindent {\bf Proof.} We refer to \cite{Sch} Section~2.\noindent 
$\diamond$ \medskip 
 
The definition of the double complex $K^{\cdot,\cdot}$ in 
\eqref{complex} was motivated by the expectation 
that the geometry on 
$X_{/\kappa}$ supports interesting arithmetical information on 
the Archimedean fiber(s) of an arithmetic 
variety. Proposition~\ref{pole-ord} shows that the real Deligne 
cohomology of $X_{/\kappa}$ carries such information. The main 
goal was to construct a complex 
and an operator $N$ acting on it which carries interesting arithmetic 
information so that the hypercohomology of ${\rm Cone}(N)^\cdot$ becomes 
isomorphic to the real Deligne cohomology of $X_{/\kappa}$. 
Proposition~4.1 of \cite{KC}  shows that the complex 
$K^{\cdot,\cdot}$ has such property. When $X_{/\kappa}$ is a projective 
non-singular 
curve defined over $\kappa = \C$ or $\R$, this can be stated as follows 
\begin{prop}\label{cone-deligne} Let $p,~q$ be two non-negative
integer. Then,  for each 
fixed valued of $p$, the 
complex ${\rm Cone}(N: K^{q-2p,q-1} \to K^{q-2p+2,q-1})$ is 
quasi-isomorphic to the complex 
\begin{multline}\label{deligne} 
\ldots\stackrel{d''}{\to}\underbrace{(A^{p-2,p-1}(p-1)\oplus 
A^{p-1,p-2}(p-1))_\R}_{q=2p-2}\stackrel{d''}{\to} 
\underbrace{A^{p-1,p-1}_\R 
(p-1)}_{q=2p-1}\stackrel{2\pi\sqrt{-1}d'd''}{\longrightarrow} 
\underbrace{A^{p,p}_\R(p)}_{q=2p} \stackrel{d'}{\to}\\ 
\stackrel{d'}{\to}\underbrace{(A^{p+1,p}(p)\oplus 
A^{p,p+1}(p))_\R}_{q=2p+1} \stackrel{d'}{\to} \ldots 
\end{multline} 
whose homology, in each degree $q$, is isomorphic to the real Deligne 
cohomology of $X_{/\C}$. 
 
Furthermore, taking the $\bar 
F_{\infty}$-invariants of such homology yields a description of 
$H^q_{\mathcal D}(X_{/\R},\R(p))$. 
\end{prop} 
\noindent {\bf Proof.} We refer to \cite{KC} Section~4, 
Proposition~4.1.\noindent $\diamond$ \medskip 
 
\begin{rem}\label{rem5} {\rm Notice that the left-side of
\eqref{deligne}, with respect to  
the central map $2\pi\sqrt{-1}d'd''$ is quasi-isomorphic to the
complex obtained by  
exchanging the map $d'$ with $d''$.} 
\end{rem} 
 
\begin{prop}\label{isos-del} For $q \ge 2p \ge 0$, the following 
isomorphisms hold 
\begin{align} 
H^q_{\mathcal D}(X_{/\C},\R(p)) &\simeq H^{q+1}_{\mathcal 
D}(X_{/\C},\R(q+1-p))\label{iso-1} \\ 
H^q_{\mathcal D}(X_{/\R},\R(p)) &\simeq 
H^{q+1}_{\mathcal D}(X_{/\C},\R(q+1-p))^{(-1)^q\bar
F_{\infty}=id}\label{iso-2}  
\end{align} 
\end{prop} 
\noindent {\bf Proof.} We consider in detail the case $q=2p$. From 
Proposition~\ref{N-isos}, the following composite map ($N^2$) is an 
isomorphism 
\[ 
gr^\w_{2(p+1)}H^{2p}(\X) \stackrel{N_{2(p+1)}}{\hookrightarrow} 
gr^\w_{2p}H^{2p}(\X) \stackrel{N_{2p}}{\twoheadrightarrow} 
gr^\w_{2(p-1)}H^{2p}(\X), \] 
where we use the notation $N_{2p}=N|_{gr^\w_{2p}H^{2p}(\X)}$.
This implies, using the results of Corollary~\ref{behv-N},  that 
$gr^\w_{2p}H^{2p}(\XX) = {\rm Ker}{N_{2p}}$ is mapped isomorphically to 
the group $gr^\w_{2(p+1)}H^{2p+1}(\XX) ={\rm Coker}(N_{2(p+1)})$: this 
isomorphism is induced by the sequence of maps \eqref{Wang} in 
Lemma~\ref{ex-seqs} (case $q = 2p$). 
It follows from Proposition~\ref{cone-deligne} that 
$gr^\w_{2p}H^{2p}(\XX) \simeq H^{2p}_{\mathcal D}(X_{/\C},\R(p))$ 
and that $$ gr^\w_{2p}H^{2p}(\XX)^{\bar F_{\infty}=id} \simeq 
H^{2p}_{\mathcal D}(X_{/\R},\R(p)).$$ Similarly, one gets from the same 
proposition that $$gr^\w_{2(p+1)}H^{2p+1}(\XX) \simeq 
H^{2p+1}_{\mathcal D}(X_{/\C},\R(p+1)),$$ hence we obtain 
\eqref{iso-1}. Taking the invariants for the action of 
$(-1)^q\bar F_\infty$  yields \eqref{iso-2}. 
The proof in the case $q\ge 2p+1$ is a generalization of the one just 
finished. For details on this part we refer to 
\cite{KC}: pp~352-3. \noindent $\diamond$ \medskip 
 
It is well known that the algebraic de Rham cohomology 
$H^\ast_{DR}(X_{/\C},\R(p))$, ($p\in\Z$) is the homology of the 
complex 
\[ 
0 \to \R(p) \to A^{0,0}_\R(p)\stackrel{d'}{\to}(A^{1,0}\oplus 
A^{0,1})_\R(p)\stackrel{d'}{\to}A^{1,1}_\R(p) \to 0;\ \ \ \  d' = \partial + 
\bar\partial. 
\] 
 
Using Lemma~\ref{descs} together with 
Proposition~\ref{cone-deligne} and Remark~\ref{rem5} we obtain the 
following description 
\begin{prop}\label{descr} Let $X = X_{/\kappa}$ be a smooth, projective 
curve over $\kappa = 
\C$ or $\R$. 
 
For $\kappa = \C$, the following description holds: 
\begin{equation}\label{HY} \begin{aligned} 
H^0(Y) &= \bigoplus_{p\le 0}gr^\w_{2p} H^0(Y) = \bigoplus_{p\le 0} 
H^0(X_{/\C},\R(p))\\ 
H^1(Y) &= \bigoplus_{p\le 0}gr^\w_{2p}H^1(Y) = \bigoplus_{p\le 0} 
H^1(X_{/\C},\R(p)) \\ 
H^2(Y) &= gr^\w_2H^2(Y)\oplus\bigoplus_{p\le 0}gr^\w_{2p}H^2(Y) = 
A^{1,1}_\R(1)\oplus\bigoplus_{p\le 0} H^2(X_{/\C},\R(p)). 
\end{aligned} 
\end{equation} 
\begin{equation}\label{HYX} 
\begin{aligned} 
H^2_Y(X) =& gr^\w_2H^2_Y(X)\oplus\bigoplus_{p\ge 
2}gr^\w_{2p}H^2_Y(X) 
\simeq A^{0,0}_\R \oplus\bigoplus_{p\ge 2}H^1_{\mathcal 
D}(X_{/\C},\R(p)) \\ 
H^3_Y(X) =& \bigoplus_{p\ge 2}gr^\w_{2p}H^3_Y(X) \simeq 
\bigoplus_{p\ge 2} H^2_{\mathcal D}(X_{/\C},\R(p)) \\ 
H^4_Y(X) =& \bigoplus_{p\ge 2}gr^\w_{2p}H^4_Y(X) 
\simeq \bigoplus_{p\ge 2}H^3_{\mathcal D}(X_{/\C},\R(p)). 
\end{aligned} 
\end{equation} 
For $X_{/\R}$ similar results hold by taking $\bar 
F_{\infty}$-invariants on both sides of the equalities. 
\end{prop} 
 
Using Proposition~\ref{descr}, the description of $H^q(\XX)$ given in 
Proposition~\ref{conedec} can be made more explicit. 
 
\begin{prop}\label{grad} Let $X_{/\kappa}$ be a smooth, projective curve 
over $\kappa = \C$ or $\R$. 
\begin{enumerate} 
\item For $\kappa = \C$ and $\forall~q \ge 0$ one has 
\[ 
H^q(\X)^{N=0} = \operatorname*{\oplus}_{p\in\Z}gr^\w_{2p}H^q(\X)^{N=0} = 
\operatorname*{\oplus}_{q\ge 2p}gr^\w_{2p}H^q(\XX) 
\] 
\item In particular: $H^q(\XX) = 0$ for $q\notin [0,3]$ and the 
following description holds 
\begin{align*} 
H^0(\XX) &= H^0(\X)^{N=0} = \bigoplus_{p\le -1}gr^\w_{2p}H^0(Y)\oplus 
gr^\w_0H^0(\XX) 
= \bigoplus_{p\le 0}H^0(X_{/\C},\R(p)) \\ 
H^1(\XX) &= H^1(\X)^{N=0} \oplus 
(gr^\w_2H^1(\XX)\oplus\bigoplus_{p\ge 2}gr^\w_{2p}H^2_Y(X)) \simeq \\ 
&\simeq\bigoplus_{p\le 0}H^1(X_{/\C},\R(p))\oplus 
\bigoplus_{p\ge 1}H^0(X_{/\C},\R(p-1))\\ 
H^{2}(\XX) &= H^2(\X)^{N=0} \oplus\bigoplus_{p\ge 2}gr^\w_{2p}H^3_Y(X) 
\simeq\\ 
&\simeq \bigoplus_{p\le 1}H^2(X_{/\C},\R(p))\oplus \bigoplus_{p\ge 
2}H^1(X_{/\C},\R(p-1))\\ 
H^3(\XX) &= (gr^\w_4H^3(\XX)\oplus\bigoplus_{p\ge 3}gr^\w_{2p}H^4_Y(X)) 
\simeq \bigoplus_{p\ge 2}H^2(X_{/\C},\R(p-1)). 
\end{align*} 
\end{enumerate} 
When $\kappa = \R$ similar results hold by taking $\bar 
 F_{\infty}$-invariants on both sides. 
\end{prop} 
 
\noindent {\bf Proof.} {\it 1.} follows from 
Corollary~\ref{behv-N}. The first statement in {\it 2.} is a 
consequence of $\dim X_{/\kappa} = 1$. For $q \in [0,3]$, the 
description of the 
graded groups $H^q(\XX)$ follows from Proposition~\ref{conedec}, 
Proposition~\ref{cone-deligne}, Proposition~\ref{isos-del}. 
For $p,q \ge 2$, the isomorphisms $H^q_{\mathcal D}(X_{/\C},\R(p)) 
\simeq H^{q-1}(X_{/\C},\R(p-1))$ are a consequence of Remark~\ref{rem5}, 
whereas the isomorphisms:  $gr^\w_2H^1(\XX) \simeq H^1_{\mathcal 
D}(X_{/\C},\R(1)) \simeq H^0_{\mathcal D}(X_{/\C},\R) \simeq 
H^0(X_{/\C},\R)$ and 
$gr^\w_2H^2(\XX) \simeq H^2_{\mathcal D}(X_{/\C},\R(1)) \simeq 
H^3_{\mathcal D}(X_{/\C},\R(2)) \simeq 
H^2(X_{/\C},\R(1))$ follow from Proposition~\ref{isos-del}. In 
particular the last isomorphism holds because $\dim X = 1$. Finally, 
the case $\kappa = \R$ is a consequence of the fact that the $\bar 
F_{\infty}$-invariants of the homology of the complex 
\eqref{deligne} give $H^\ast_{\mathcal D}(X_{/\R},\R(p))$. \noindent 
$\diamond$ 
 
\subsection{Archimedean Frobenius and regularized determinants}\label{2.5} 
 
On the {\it infinite dimensional} real vector space 
$gr^\w_{2p}H^\ast(\X)$ (\cf\eqref{X-co}) one defines a linear operator 
\begin{equation} \label{Phi} \Phi:  gr^\w_{2p} H^\ast(\X) \to gr^\w_{2p} 
H^\ast(\X), \ \ \ \ \ \  \Phi(x) = p\cdot x \end{equation} 
and then extend this definition to the whole group $H^\ast(\X)$ 
according to the decomposition $H^q(\X) = \oplus_p gr^\w_{2p}H^q(\X)$. 
 
In this section we will consider the operator $\Phi$ restricted to the 
subspace $H^\ast (\X)^{N=0}.$ Following the description of this 
space given in Proposition~\ref{grad}, we write $\Phi = \oplus_{q=0}^2 
\Phi_q$, where 
\[ 
 \Phi_q : H^q (\X)^{N=0} \to H^q (\X)^{N=0}. 
\] 
 
Given a self-adjoint operator $T$ with pure point spectrum, the 
{\it zeta-regularized determinant} is defined by  
\begin{equation}\label{reg-det0} 
\det_\infty ( s - T ) = \exp \left( 
-\frac{d}{dz}~\zeta_{T}(s,z)_{|{z=0}} \right), 
\end{equation} 
where 
\begin{equation}\label{reg-zeta} 
\zeta_{T}(s,z) = \sum_{\lambda \in {\rm Spec}(T)} 
m_\lambda (s-\lambda)^{-z}. 
\end{equation} 
Here, ${\rm Spec}(T)$ denotes the spectrum of $T$ and $m_\lambda =\dim
E_\lambda(T)$ is the multiplicity of the eigenvalue $\lambda$ with
eigenspace $E_\lambda(T)$. 
 
$\Phi_q$ is a self-adjoint operator with respect to the inner product 
induced by \eqref{innerprod}, with spectrum 
$$ {\rm Spec}(\Phi_q) = \left\{ \begin{array}{lr} \{ n\in \Z, n\leq 0 \} & 
q=0,1 \\ \{ n \in \Z, n \leq 1 \} & q =2. \end{array} \right. $$ 
 
The eigenspaces $E_{n}(\Phi_q)=gr_{2n}^\w H^q (\X)^{N=0}$ have 
dimensions $\dim E_{n}(\Phi_q) = b_q$, the $q$-th Betti number of 
$X_{/\C}$.

\begin{prop}\label{infdet} 
The zeta regularized determinant of $\Phi_q$ is given by 
\begin{equation}\label{detPhiq} 
\det_\infty \left( \frac{s}{2\pi} - \frac{\Phi_q}{2\pi}\right) = \Gamma_\C 
(s)^{-b_q}, 
\end{equation} 
for $q=0,1$ and 
\begin{equation}\label{detPhi2} 
\det_\infty \left( \frac{s}{2\pi} - \frac{\Phi_2}{2\pi}\right) = \Gamma_\C 
(s-1)^{-b_2}, 
\end{equation} 
with $\Gamma_\C(s)$ and $\Gamma_\R(s)$ as in \eqref{gammas}. 
\label{detC} 
\end{prop} 
 
\noindent {\bf Proof.} We write explicitly the zeta function for the operator 
$\Phi_q/(2\pi)$. 
When $q=0,1$, this has spectrum $\{ n/(2\pi) \}_{n\leq 0}$, 
hence we have 
$$ \zeta_{\Phi_q/(2\pi)} (s/(2\pi),z)= \sum_{n\leq 0} b_q (s/(2\pi) - 
n/(2\pi))^{-z}= b_q (2\pi)^z \zeta(s,z). $$ 
 
$\zeta(s,z)$ is the Hurwitz zeta function 
$$ \zeta(s,z) = \sum_{n\geq 0} \frac{1}{(s+n)^z}. $$ 
For $q=2$, similarly we have 
$$ \zeta_{\Phi_2/(2\pi)} (s/(2\pi),z)= b_2 (2\pi)^z (\zeta(s,z) + 
(s-1)^{-z}). $$ 
 
It is well known that the Hurwitz zeta function satisfies the 
following properties: 
\begin{equation}\label{H-zeta} 
 \zeta(s,0) = \frac{1}{2} -s, \ \ \ \ \ 
\frac{d}{dz}~\zeta(s,z)_{|_{z=0}} = \log \Gamma(s) - \frac{1}{2} \log (2\pi). 
\end{equation} 
When $q = 0,1$, the computation of 
$\frac{d}{dz}~\zeta_{\Phi_q/(2\pi)}(s/(2\pi),z)_{|{z=0}}$ yields 
 
$$ \frac{d}{dz}~\zeta_{\Phi_q/(2\pi)}(s/(2\pi),z) 
 = b_q \left(\log (2\pi) (2\pi)^{z} \zeta(s,z) + (2\pi)^{z} 
\frac{d}{dz}\zeta(s,z) \right). $$ 
At $z=0$, this gives 
$$ \frac{d}{dz}~\zeta_{\Phi_q/(2\pi)}(s/(2\pi),z)_{|_{z=0}} = 
b_q \left( \log (2\pi) \zeta(s,0) + \frac{d}{dz}\zeta(s,z)_{|_{z=0}} 
\right) $$ 
$$ = b_q \left( \log (2\pi) (\frac{1}{2} -s) + \log 
\Gamma(s) - \frac{1}{2} \log (2\pi) \right) 
 = b_q ( -s\log (2\pi) + \log \Gamma(s) ). $$ 
Taking the exponential we get 
\begin{align*} 
 \exp \left( -  \frac{d}{dz}~\zeta_{\Phi_q/(2\pi)}(s/(2\pi),z)_{|_{z=0}} 
\right) &=\exp(-b_q ( -s\log (2\pi) + \log \Gamma(s) )) = \\ 
&= \left( (2\pi)^{-s} \Gamma(s) \right)^{-b_q} 
= \Gamma_{\C}(s)^{-b_q}. 
\end{align*} 
 
When $q=2$, one has similarly 
$$ \frac{d}{dz}~\zeta_{\Phi_2/(2\pi)}(s/(2\pi),z) = 
b_2 \left( (\log (2\pi) (2\pi)^{z} \zeta(s,z) + (2\pi)^{z} 
\frac{d}{dz}\zeta(s,z) \right) 
+ b_2 \frac{d}{dz}~\frac{(2\pi)^z}{(s-1)^z}. $$ 
At $z=0$ this gives 
$$ b_2 \left( -s\log (2\pi) + \log \Gamma(s)  + \log(2\pi) -\log(s-1) 
\right). $$ 
Thus, we have 
\begin{align*} 
 \exp \left( -  \frac{d}{dz}~\zeta_{\Phi_2/(2\pi)}(s/(2\pi),z)_{ 
|_{z=0}} \right) &= \left( (2\pi)^{-s+1} \frac{\Gamma(s)}{(s-1)} 
\right)^{-b_2} = \left( (2\pi)^{-(s-1)}\Gamma(s-1) \right)^{-b_2} =\\\\ 
&= \Gamma_{\C}(s-1)^{-b_2}. 
\end{align*} 
\noindent $\diamond$ 
 
\begin{rem} {\rm When $X_{/\C}$ is a smooth complex 
algebraic curve, that is, when $X_{/\C} = X_{\alpha(K)}$ for an Archimedean 
prime that corresponds to a complex (non-real) embedding $\alpha: K 
\hookrightarrow \C$, the description of the complex Euler factor 
is given by (\cf\eqref{factors}) 
$$ L_\C (H^q(X_{/\C},\C),s) = \left\{ \begin{array}{lr} \Gamma_\C 
(s)^{b_q} & q=0,1 \\ 
\Gamma_\C (s-1)^{b_2} & q=2, \end{array} \right. $$ 
where $H^q(X_{/\C},\C)$ is the Betti 
cohomology. The relation to the determinants \eqref{detPhiq} 
\eqref{detPhi2} is then 
\begin{equation}\label{det=L} 
\det_\infty \left( \frac{s}{2\pi} - \frac{\Phi_q}{2\pi}\right)^{-1} 
=L_\C (H^q(X_{/\C},\C),s). 
\end{equation} 
This result was proved in \cite{KC}: \S 5, via comparison to Deninger's 
pair $(H^\ast_{ar},\Theta)$.} 
\end{rem} 
 
Assume now that $X_{/\R}$ is a smooth real algebraic 
curve; that is $X_{/\R} = X_{\alpha(K)}$ for an Archimedean  prime 
that corresponds to a real embedding $\alpha: K \hookrightarrow 
\R$. In this case $X_{/\R}$ is a symmetric Riemann surface, namely a 
compact Riemann surface with an involution $\iota: X_{/\R} \to X_{/\R}$ 
induced by complex conjugation. 
Such involution on the manifold induces an action of the real Frobenius 
$\bar F_\infty$ on $H^q(\X)^{N=0}$: we refer to Remark~\ref{rem4} for 
the description of this operator. 
 
For instance, following the decomposition 
given in Proposition~\ref{grad}, $$H^1 (\X)^{N=0} = \oplus_{p\le 
0}gr^\w_{2p}H^1(\XX)$$ splits as the sum 
of two eigenspaces for $\bar F_\infty$ with eigenvalues $\pm 1$: 
$$ H^1 (\X)^{N=0} = E^+ \oplus E^-, $$ 
where 
\begin{align}\label{deco} 
E^+ &:= H^1 (\X)^{N=0, \bar F_\infty=id} = \bigoplus_{p\leq 0} E_1(2p) 
\oplus \bigoplus_{p\leq -1} E_{-1}(2p+1), \\ 
E^- &:= H^1 (\X)^{N=0, \bar F_\infty= - id} = \bigoplus_{p\leq -1} 
E_1(2p+1)\oplus \bigoplus_{p\leq 0} E_{-1}(2p). \nonumber 
\end{align} 
 
\medskip 
 
We consider once more the operator $\Phi$ acting on $H^\ast (\X)^{N=0}$, 
and we denote by $\hat\Phi_q$ the restriction of this operator to the 
subspace $H^q (\X)^{N=0, \bar F_\infty =id}$. 
 
\begin{prop} The regularized determinant for the operator 
$ \hat\Phi_q = \Phi_{|_{H^q (\X)^{N=0, \bar F_\infty =id}}}$ 
is given by ($g =$ genus of $X_{/\kappa}$) 
\begin{equation}\label{detPhi0R} 
\det_\infty \left( \frac{s}{2\pi} - \frac{\hat\Phi_0}{2\pi}\right) = 
\Gamma_\R(s)^{-b_0} 
\end{equation} 
\begin{equation}\label{detPhi1R} 
\det_\infty \left( \frac{s}{2\pi} - \frac{\hat\Phi_1}{2\pi}\right) = 
\Gamma_\R(s)^{-b_1/2} \Gamma_\R(s+1)^{-b_1/2} = \Gamma_\C (s)^{-g} 
\end{equation} 
\begin{equation}\label{detPhi2R} 
\det_\infty \left( \frac{s}{2\pi} - \frac{\hat\Phi_2}{2\pi}\right) = 
\Gamma_\R(s-1)^{-b_2}. 
\end{equation} 
\label{detR} 
\end{prop} 
 
\noindent {\bf Proof.} We write explicitly the zeta function for the 
operators $\hat\Phi_q$ on $H^q (\X)^{N=0, \bar F_\infty=id}$. The 
spectrum of $\hat\Phi_q$ is given by $\{ n\in \Z, n\leq 0 \}$ for 
$q=0,1$ and $\{ n\in \Z, n\leq 1 \}$ for $q=2$. Because complex 
conjugation is the identity on $H^0(X_{/\R},\R(2n))$, the eigenspaces 
where $\bar F_\infty$ acts as identity are $E_n(\hat\Phi_0) = 
gr_{4n}^\w H^0(\XX)$. 
On the other hand, $\bar F_{\infty}$ acts as the identity on 
$H^2(X_{/\R},\R(2n+1))$, hence $E_n(\hat\Phi_2) = 
gr_{2(2n+1)}^\w H^2(\XX)$. The action of $\bar F_\infty$ on 
$H^1(X_{/\R},\R(n))$ is the identity precisely on the eigenspaces 
$E_n(\hat\Phi_1) = E_1(2n) \oplus E_{-1}(2n+1)$ as in \eqref{deco}. 
 
The zeta function of $\hat\Phi_0/(2\pi)$ is therefore of the form 
$$ \zeta_{\hat\Phi_0/(2\pi)}(s/(2\pi),z) =\sum_{n\geq 0} b_0 \left( 
\frac{s+2n}{2\pi} \right)^{-z}= b_0 (2\pi)^z \sum_{n\geq 0} 
\frac{1}{(s+2n)^z} = b_0 (\pi)^z \zeta(s/2,z), $$ 
where $\zeta(s,z)$ is the Hurwitz zeta function. Using the identities 
\eqref{H-zeta}, we obtain 
$$ \frac{d}{dz}~\zeta_{\hat\Phi_0/(2\pi)}(s/(2\pi),z)_{|_{z=0}} = 
b_0 ( \log(\pi) (1/2 - s/2) + \log \Gamma(s/2) -1/2 \, \log(2\pi)  ) $$ 
$$ = b_0 ( -s/2 \log(\pi) -1/2 \log(2) +\log  \Gamma(s/2) ).$$ 
Hence, using the equalities \eqref{gammas}, we obtain 
$$ \exp\left( -\frac{d}{dz}~\zeta_{\hat\Phi_0/(2\pi)}(s/(2\pi),z)_{|_{z=0}} 
\right) = \exp( -b_0 ( -s/2 \log(\pi) -1/2 \log(2) +\log  \Gamma(s/2) ) 
) $$ 
$$ = \left( 2^{-1/2} \pi^{-s/2} \Gamma(s/2) 
\right)^{-b_0}=\Gamma_\R(s)^{-b_0}. $$ 
 
The determinant for $\hat\Phi_1/(2\pi)$ is given by the product 
$$ \det_\infty \left( \frac{s}{2\pi} - \frac{\hat\Phi_1}{2\pi}\right) 
= $$ 
$$ \det_\infty \left( \frac{s}{2\pi} - 
\frac{\hat\Phi_1}{2\pi}|_{\oplus_n E_1(2n)} \right) \cdot 
\det_\infty \left( \frac{s}{2\pi} - 
\frac{\hat\Phi_1}{2\pi}|_{\oplus_n E_{-1}(2n+1)}\right). $$ 
The zeta function for the 
first operator is given by 
$$ \zeta_{\frac{\hat\Phi_1}{2\pi}|_{\oplus_n E_1(2n)}}(s/(2\pi),z) = 
\frac{b_1}{2} \pi^z \zeta(s/2,z) $$ 
while the for the second operator is 
$$ \zeta_{\frac{\hat\Phi_1}{2\pi}|_{\oplus_n E_{-1}(2n+1)}}(s/(2\pi),z) = 
 \frac{b_1}{2}(2\pi)^z \sum_{n\geq 0} \frac{1}{(s+1+2n)^z} = 
\frac{b_1}{2}\pi^z \zeta((s+1)/2,z). $$ 
Thus, we obtain 
$$ \det_\infty \left( \frac{s}{2\pi} - 
\frac{\hat\Phi_1}{2\pi}|_{\oplus_n E_1(2n)} \right) = \Gamma_\R 
(s)^{-b_1/2} $$ 
and 
$$ \det_\infty \left( \frac{s}{2\pi} - 
\frac{\hat\Phi_1}{2\pi}|_{\oplus_n E_{-1}(2n+1)}\right) = \Gamma_\R 
(s+1)^{-b_1/2}. $$ 
 
Then, \eqref{detPhi1R} follows using the equality $\Gamma_\R(s) 
\Gamma_\R(s+1) = \Gamma_\C(s)$. 
 
Finally, for $\hat\Phi_2/(2\pi)$, we have 
$$ \zeta_{\frac{\hat\Phi_2}{2\pi}}\left(\frac{s}{2\pi},z \right) = 
b_2 (2\pi)^z \sum_{n\geq 0} \frac{1}{(s-1+2n)^z} =  b_2  \pi^z 
\zeta((s-1)/2,z), $$ 
hence 
$$ \frac{d}{dz}~\zeta_{\frac{\hat\Phi_2}{2\pi}}\left(\frac{s}{2\pi},z 
\right)_{|_{z=0}} = b_2 \left( -(s-1)/2 \log \pi - 1/2 \, \log 2 + \log 
\Gamma((s-1)/2) \right).$$ 
 
Therefore 
$$ \det_\infty \left( \frac{s}{2\pi} - 
\frac{\hat\Phi_2}{2\pi} \right) = \left( 2^{-1/2} \pi^{-(s-1)/2} 
\Gamma((s-1)/2) \right)^{-b_2} = \Gamma_\R (s-1)^{-1}. $$ 
 
\noindent $\diamond$ 
 
\begin{rem} {\rm When $X_{/\R}$ is a smooth, real 
algebraic curve of genus $g$, that is, when $X_{/\R} = X_{\alpha(K)}$ 
for an Archimedean  
prime that corresponds to a real embedding $\alpha: K 
\hookrightarrow \R$, the description of the real Euler factor 
is given by (\cf\eqref{factors}) 
$$ L_\R (H^q(X_{/\R},\R),s) = \left\{ \begin{array}{lr} \Gamma_\R 
(s) & q=0 \\ 
\Gamma_\C (s)^g & q=1 \\ 
\Gamma_\R (s-1) & q=2, \end{array} \right. $$ 
As for the complex case, this result was proved in \cite{KC}: \S 5, 
via comparison to Deninger's pair $(H^\ast_{ar},\Theta)$.} 
\end{rem}

\section{Arithmetic spectral triple.} \label{4} 
 
In this Section we show that the polarized Lefschetz module structure
of Theorem~\ref{main thm} together with the operator $\Phi$ define a
``cohomological'' version of the structure of a {\em spectral 
triple} in the sense of Connes (\cf \cite{Co94} \S VI).

\medskip 
 
In this Section, we will use {\em real coefficients}. In 
fact, in order to introduce spectral data compatible with the 
arithmetic construction of Section \ref{3}, we need to preserve 
the structure of real vector spaces. For this reason, the algebras  we
consider in this construction will be real group rings. 
 
\medskip 

Let $(H^\cdot (\XX),\Phi)$ be the cohomological theory of the fiber at 
the Archimedean  prime introduced in Section \ref{3}, endowed with the 
structure of polarized Lefschetz module. 

In Theorem \ref{repres} we show that the Lefschetz representation of
$\SL(2,\R)$ given by the Lefschetz module structure on $K^\cdot$
induces a representation  
\begin{equation}\label{rhoS3} 
 \rho : \SL(2,\R) \to {\mathcal B}(H^\cdot 
(\X)), 
\end{equation} 
where ${\mathcal B}(H^\cdot (\X))$ is the algebra of 
bounded operators on a real Hilbert space completion of $H^\cdot 
(\X)$ (in the inner product determined by the polarization of 
$K^{\cdot,\cdot}$). The representation $\rho$ extends to the real
group ring compatibly with the Lefschetz module structure on $H^\cdot
(\X)=\H^\cdot (K,d)$. We work with the group ring, since for the
purpose of this paper we are interested in considering the restriction
of \eqref{rhoS3} to certain discrete subgroups of
$\SL(2,\R)$. A formulation in terms of the Lie algebra and its
universal enveloping algebra will be considered elsewhere.
 
The main result of this section is Theorem \ref{spectral3}, where we 
prove that the inner product on $H^\cdot (\X)$ and the representation 
\eqref{rhoS3} induce an inner product on $H^\cdot (\XX) = 
\H^\cdot({\rm Cone}^\cdot(N))$ and a corresponding representation
$\rho_N$ in ${\mathcal B}(H^\cdot 
(\XX))$. We then consider the spectral data $({\rm A},H^\cdot 
(\XX),\Phi)$, where ${\rm A}$ is the image under $\rho_N$ of the real
group ring, and show that the operator $\Phi$ satisfies the properties
of a Dirac operator (in the sense of Connes' theory of spectral
triples), which has bounded commutators with the elements of ${\rm A}$. 

\medskip 
 
In non-commutative geometry, the notion of a spectral triple 
provides the correct generalization of the classical structure of 
a Riemannian manifold. The two notions agree on a commutative 
space. In the usual context of Riemannian geometry, the definition 
of the infinitesimal element $ds$ on a smooth spin manifold can be 
expressed in terms of the inverse of the classical Dirac operator 
$D$. This is the key remark that motivates the theory of spectral 
triples. In particular, the geodesic distance between two points 
on the manifold is defined in terms of $D^{-1}$ (\cf \cite{Co94} 
\S VI). The spectral triple that describes a classical Riemannian 
spin manifold is $(A,H,D)$, where $A$ is the algebra of complex 
valued smooth functions on the manifold, $H$ is the Hilbert space 
of square integrable spinor sections, and $D$ is the classical 
Dirac operator (a square root of the Laplacian). These data 
determine completely and uniquely the Riemannian geometry on the 
manifold. It turns out that, when expressed in this form, the 
notion of spectral triple extends to more general non-commutative 
spaces, where the data $(A,H,D)$ consist of a ${\rm C}^*$-algebra 
$A$ (or more generally of a smooth subalgebra of a ${\rm 
C}^*$-algebra) with a representation as bounded operators on a 
Hilbert space $H$, and an operator $D$ on $H$ that verifies the 
main properties of a Dirac operator. The notion of smoothness 
is determined by $D$: the smooth elements of $A$ are defined 
by the intersection of domains of powers of the derivation given 
by commutator with $|D|$.

\medskip

The basic geometric structure encoded by the theory of spectral
triples is Riemannian geometry, but in more refined cases, such as
K\"ahler geometry, the additional structure can be easily encoded as
additional symmetries. In our case, for instance, the algebra ${\rm
A}$ corresponds to the action of the Lefschetz operator, hence it
carries the information (at the cohomological level) on the K\"ahler
form. 

\medskip 

In the theory of specral triples, in general, the Hilbert space $H$ is
a space of {\em cochains} on which the natural algebra of the geometry
is acting. Here we are considering a simplified triple of spectral
data defined on the {\em cohomology}, hence we do not expect the full
algebra describing the geometry at arithmetic infinity to
act. We show in Theorem \ref{Gamma-zeta3} that the spectral data
$({\rm A},H^\cdot (\XX),\Phi)$ are sufficient to recover the
alternating product of the local factor. In fact, 
the theory of spectral triples encodes important arithmetic 
information on the underlying non-commutative space, expressed via 
an associated family of {\em zeta functions}. By studying the zeta 
functions attached to the data $({\rm A},H^\cdot 
(\XX),\Phi)$, we find a natural one whose associated Ray-Singer 
determinant is the alternating product of the $\Gamma$-factors for 
the real Hodge structure over $\C$ given by the Betti cohomology 
$H^q(X_{/\C}, \C)$. A more refined 
construction of a spectral triple associated to the Archimedean  places
of an arithmetic surface (using the full complex $K^\cdot$ instead of
its cohomology) will be considered elsewhere.    

\medskip 

Moreover, we show that, in the case of a  Riemann surface $X_{/\C}$ of
genus $g\geq 2$, one can enrich the cohomological spectral data 
$({\rm A},H^\cdot (\XX),\Phi)$ by the additional datum of a Schottky
uniformization. Given the group $\Gamma\subset 
\PSL(2,\C)$, which gives a Schottky uniformization of the Riemann 
surface $X_{/\C}$ and of the hyperbolic handlebody $\mX_\Gamma \cup 
X_{/\C} =\Gamma \backslash (\H^3 \cup \Omega_\Gamma)$,  
Bers simultaneous uniformization (\cf \cite{Bers} 
\cite{Bo}) determines a pair of Fuchsian Schottky groups $G_1,G_2\subset 
\SL(2,\R)$, which correspond geometrically to a decomposition of 
the Riemann surface $X_{/\C}$ as the union of two Riemann surfaces 
with boundary. We let the Fuchsian Schottky groups act on the complex 
$K^{\cdot,\cdot}$ and on the cohomology via the restriction of the 
representation of $\SL(2,\R)\times \SL(2,\R)$ of the 
Lefschetz module structure to a normal subgroup $\tilde\Gamma \subset
\Gamma$ determined by the simultaneous uniformization, with 
$G_1 \simeq \tilde\Gamma \simeq G_2$. Geometrically,
this group corresponds to the choice of a covering $\mX_{\tilde\Gamma} \to 
\mX_{\Gamma}$ of $\mX_{\Gamma}$ by a handlebody 
$\mX_{\tilde\Gamma}$. The image ${\rm A}(\tilde\Gamma)$ of
the group ring of $\tilde\Gamma$ under the representation $\rho_N$ 
encodes in the spectral data $({\rm A},H^\cdot (\XX),\Phi)$ the
information on the topology of $\mX_{\tilde\Gamma}$. 
 
In the interpretation of the tangle of bounded geodesics in the handlebody 
$\mX_{\tilde\Gamma}$ as the dual graph of the closed fiber at 
arithmetic infinity, the covering 
$\mX_{\tilde\Gamma} \to \mX_{\Gamma}$ produces a corresponding 
covering of the dual graph by geodesics in $\mX_{\tilde\Gamma}$. 
Passing to the covering $\mX_{\tilde\Gamma}$ may be regarded as an 
analog, at the Archimedean  primes, of the 
refinement of the dual graph of a Mumford curve that corresponds to 
a minimal resolution (\cf \cite{Mum} \S 3). 
 
\smallskip 
 
When the Archimedean  prime corresponds to a real embedding 
$K\hookrightarrow \R$, so that the corresponding Riemann surface 
$X_{/\R}$ acquires a real structure, Proposition \ref{real-rep} 
shows that if $X_{/\R}$ is a smooth orthosymmetric real algebraic
curve (in particular, the set of real points $X_{/\R}(\R)$ is
non-empty), then there is a preferred choice of a Fuchsian Schottky
group $\Gamma$ determined by the real structure,
for which the simultaneous uniformization consists of
cutting the Riemann surface along $X_{/\R}(\R)$.

\subsection{Spectral triples} 
 
We recall the basic setting of Connes theory of spectral triples. 
For a more complete treatment we refer to \cite{Connes}, 
\cite{Co94}, \cite{ConnesMosc}. 
 
\smallskip 
 
\begin{defn}\label{specDef} 
a spectral triple $({\mathcal A}, {\mathcal H}, D)$ consists of 
an involutive algebra ${\mathcal A}$ with a representation 
$$ \rho : {\mathcal A} \to {\mathcal B}({\mathcal H}) $$ 
as bounded operators on a Hilbert space ${\mathcal H}$, and an 
operator  $D$ (called the Dirac operator) on ${\mathcal H}$, which 
satisfies the following properties: 
\begin{enumerate} 
\item $D$ is self--adjoint. 
\item For all $\lambda\notin \R$, the resolvent $(D-\lambda)^{-1}$ is 
a compact operator on ${\mathcal H}$. 
\item For all $a\in {\mathcal A}$, the commutator $[D,a]$ is a bounded 
operator on ${\mathcal H}$. 
\end{enumerate} 
\end{defn} 
 
\begin{rem}\label{defSP3} {\em 
The property {\em 2.} of Definition \ref{specDef} generalizes 
ellipticity of the standard Dirac operator on a compact manifold. 
Usually, the involutive algebra ${\mathcal A}$ satisfying property 
{\em 3.} can be chosen to be a dense subalgebra of a ${\rm 
C}^*$--algebra. This is the case, for instance, when we consider 
smooth functions on a manifold as a subalgebra of the commutative 
${\rm C}^*$-algebra of continuous functions. In the classical case of 
Riemannian manifolds, property {\em 3.} is equivalent the Lipschitz 
condition, hence it is satisfied by a larger class than that of 
smooth functions. In {\em 3.} we write $[D,a]$ as shorthand for 
the extension to all of ${\mathcal H}$ of the operator $[D,\rho(a)]$ 
defined on the domain ${\rm Dom}(D) \cap \rho(a)^{-1}({\rm Dom}(D))$, 
where ${\rm Dom}(D)$ is the domain of the unbounded operator $D$.} 
\end{rem} 
 
We review those aspects of the theory of spectral triples which are of 
direct interest to us. For a more general treatment we refer to 
\cite{Connes}, \cite{Co94}, \cite{ConnesMosc}. 
 
\smallskip 
 
\noindent {\bf Volume form.} A spectral triple $({\mathcal A}, 
{\mathcal H}, D)$ is said to be 
of dimension $n$, or $n$--{\it summable} if the operator 
$|D|^{-n}$ is an infinitesimal of 
order one, which means that the eigenvalues $\lambda_k(|D|^{-n})$ 
satisfy the estimate $\lambda_k(|D|^{-n})=O(k^{-1})$. 
 
For a positive compact operator $T$ such that 
$$ \sum_{j=0}^{k-1} \lambda_j(T) = O(\log k), $$ 
the Dixmier trace $\Tr_\omega(T)$ is the coefficient of this 
logarithmic divergence, namely 
\begin{equation}\label{DixTr} 
\Tr_\omega(T) = \lim_\omega \frac{1}{\log k} 
\sum_{j=1}^k \lambda_j(T). \end{equation} 
Here the notation $\lim_\omega$ takes into account the fact that 
the sequence 
$$ S(k,T):=\frac{1}{\log k} \sum_{j=1}^k \lambda_j(T) $$ 
is bounded though possibly non-convergent. For this reason, the usual 
notion of limit is replaced by a 
choice of a linear form $\lim_\omega$ on the set of bounded sequences 
satisfying suitable conditions that extend analogous properties 
of the limit. When the sequence $S(k,T)$ converges \eqref{DixTr} is 
just the ordinary limit $\Tr_\omega(T)= \lim_{k\to \infty} S(k,T)$. 
So defined, the Dixmier trace \eqref{DixTr} extends to any 
compact operator that is an infinitesimal of order one, since any such 
operator is a combination $T=T_1-T_2+i(T_3-T_4)$ of positive ones $T_i$. 
The operators for which the Dixmier trace does not depend on 
the choice of the linear form $\lim_\omega$ are called {\em measurable 
operators}. 
 
On a non-commutative space the operator $|D|^{-n}$ generalizes the 
notion of a volume form. The volume is defined as 
\begin{equation} \label{VolDix} V = \Tr_\omega 
(|D|^{-n}). \end{equation} 
More generally, consider the algebra $\tilde {\mathcal A}$ generated 
by ${\mathcal A}$ and $[D,{\mathcal A}]$. Then, for $a\in \tilde 
{\mathcal A}$, integration with respect to the volume form $|D|^{-n}$ 
is defined as 
\begin{equation} \label{intDix} \int a := \frac{1}{V} \Tr_\omega (a 
|D|^{-n}). \end{equation} 
 
The usual notion of integration on a Riemannian spin manifold $M$ can be 
recovered in this context (\cf \cite{Co94}) through the 
formula ($n$ even): 
$$ \int_M f dv = \left( 2^{n-[n/2]-1} \pi^{n/2} n \Gamma(n/2) \right) 
\Tr_\omega (f |D|^{-n}). $$ 
Here $D$ is the classical Dirac operator on $M$ associated to the 
metric that determines the volume form $dv$, and $f$ in the right hand 
side is regarded as the multiplication operator acting on the Hilbert 
space of square integrable spinors on $M$. 
 
\medskip 
 
\noindent {\bf Zeta functions.} An important function associated to 
the Dirac operator $D$ of a spectral triple 
$({\mathcal A}, {\mathcal H}, D)$ is its zeta function 
\begin{equation}\label{zetaD} 
 \zeta_D (z) := \Tr(|D|^{-z}) = \sum_{\lambda} \Tr(\Pi(\lambda,|D|)) 
\lambda^{-z}, \end{equation} 
where $\Pi(\lambda,|D|)$ denotes the orthogonal projection on the 
eigenspace $E(\lambda,|D|)$. 
 
\smallskip 
 
An important result in the theory of spectral triples (\cite{Co94} \S 
IV Proposition 4) relates the volume \eqref{VolDix} with the residue 
of the zeta function \eqref{zetaD} at $s=1$ through the formula 
\begin{equation} \label{VolRes} 
V = \lim_{s\to 1+} (s-1)\zeta_D(s) = Res_{s=1} \Tr 
(|D|^{-s}). \end{equation}

\smallskip 
 
There is a family of zeta functions associated to a spectral 
triple $({\mathcal A}, {\mathcal H}, D)$, to which \eqref{zetaD} 
belongs. For an operator $a\in \tilde {\mathcal A}$, we can define the 
zeta functions 
\begin{equation} \label{aDzeta} 
\zeta_{a,D}(z) := \Tr (a |D|^{-z}) = \sum_{\lambda} 
\Tr(a\, \Pi(\lambda,|D|)) \lambda^{-z} 
\end{equation} 
and 
\begin{equation} \label{aDzetas} 
\zeta_{a,D}(s,z) := \sum_{\lambda} \Tr(a\, \Pi(\lambda,|D|)) 
(s-\lambda)^{-z}. 
\end{equation} 
These zeta functions are related to the heat kernel $e^{-t|D|}$ by 
Mellin transform 
\begin{equation} \label{Mellin} 
\zeta_{a,D}(z) = \frac{1}{\Gamma(z)} \int_0^\infty t^{z-1} \Tr( a\, 
e^{-t|D|} )\, dt 
\end{equation} 
where 
\begin{equation} \label{theta} 
\Tr( a\, e^{-t|D|} ) = \sum_\lambda \Tr(a\, \Pi(\lambda,|D|)) 
e^{-t\lambda} =: \theta_{a,D}(t). 
\end{equation} 
Similarly, 
\begin{equation} \label{Mellins} 
\zeta_{a,D}(s,z) = \frac{1}{\Gamma(z)} \int_0^\infty \theta_{a,D,s}(t) 
\, t^{z-1} \, dt 
\end{equation} 
with 
\begin{equation} \label{thetas} 
\theta_{a,D,s}(t) := \sum_\lambda \Tr(a\, \Pi(\lambda,|D|)) 
e^{(s-\lambda)t}. 
\end{equation} 
Under suitable hypothesis on the asymptotic expansion of 
\eqref{thetas} (\cf Theorem 2.7-2.8 of \cite{Man4} \S 2), the functions 
\eqref{aDzeta} and \eqref{aDzetas} admit a unique analytic 
continuation (\cf \cite{ConnesMosc}) and there is an  associated 
regularized determinant in the sense of Ray--Singer (\cf \cite{RS}): 
\begin{equation} \label{zetadet} 
 {\det_\infty}_{a,D}(s) := \exp \left( -\frac{d}{dz} \zeta_{a,D}(s,z) 
|_{z=0} \right) 
\end{equation} 
 
\smallskip 
 
The family of zeta functions \eqref{aDzeta} also provides a refined 
notion of dimension for a spectral triple $({\mathcal A}, {\mathcal 
H}, D)$, called the {\it dimension spectrum}. This is a subset 
$\Sigma=\Sigma({\mathcal A}, {\mathcal H}, D)$ in $\C$ with the 
property that all the zeta functions \eqref{aDzeta}, 
as $a$ varies in $\tilde {\mathcal A}$, extend holomorphically to $\C 
\setminus \Sigma$. 
 
\medskip 
 
\subsection{Lefschetz modules and cohomological spectral
data}\label{cohsp3} 

We consider the polarized bigraded Lefschetz module 
$(K^{\cdot,\cdot}, N, \ell, \psi)$ associated to the Riemann surface 
$X_{/\C}$ at an Archimedean  prime, as described in Section \ref{3}. 

We set 
\begin{equation}\label{PhiK} 
\tilde\Phi: K^{i,j,k} \to  K^{i,j,k}   \ \ \ \ \ \   \tilde\Phi(x) = 
\frac{(1+j-i)}{2} x. 
\end{equation} 
The operator $\tilde\Phi$ induces the operator $\Phi$ of \eqref{Phi} 
on the cohomology $H^\cdot(\X)^{N=0}$.

We have the following result. 
 
\begin{thm} \label{repres} Let $(K^{\cdot,\cdot}, d, N, \ell, \psi)$  be 
the polarized bigraded Lefschetz module associated to a Riemann 
surface $X_{/\C}$. Then the following holds. 
\begin{enumerate} 
\item The group $\SL(2,\R)$ acts, via the representation 
$\sigma_2$ of Lemma \ref{repGamma}, by bounded operators on the 
Hilbert completion of $\H^\cdot (K,d)$ in the inner product 
defined by the polarization $\psi$. This defines a representation 
\begin{equation}\label{rho-preS3} \rho: \SL(2,\R)  \to 
{\mathcal B}(\H^\cdot (K,d)) 
\end{equation} 
\item Let ${\rm A}$ be the image of the group ring in
${\mathcal B}(\H^\cdot (K,d))$, obtained by extending
\eqref{rho-preS3}. Then the operator  
$\tilde \Phi$ defined in \eqref{PhiK} has bounded commutators with 
all the elements in ${\rm A}$. 
\end{enumerate} 
\end{thm} 
 
\noindent {\bf Proof.} {\em 1.} The representations $\sigma_1$ and 
$\sigma_2$ of Lemma \ref{repGamma} extend by linearity to 
representations $\sigma_i$ of the real group ring in ${\rm 
Aut}(K)$. 
By Theorem \ref{main thm} and Corollary \ref{main cor}, the 
cohomology $\H^\cdot (K,d)$ has an induced Lefschetz module 
structure, thus we obtain induced actions of the real group 
ring on $\H^\cdot (K,d)$. 
We complete $H^\cdot (\X)=\H^\cdot (K,d)$ to a real Hilbert space with 
respect to the inner product induced by the polarization $\psi$. 
Consider operators of the form \eqref{rep1} 
with $b=1$, 
$$ U_a (x) := \sigma\left\{ \left(\begin{array}{cc} a & 0 \\ 0 & a^{-1} 
\end{array} \right) , \left(\begin{array}{cc} 1 & 0 \\ 0 & 1 
\end{array} \right) \right\} (x) = a^i x \ \ \ \text{ for $x\in 
K^{i,j}$}. $$ 
A direct calculation shows that the $U_a$ are in general {\em 
unbounded operators}: since the index $i$ varies over a 
countable set, it is not hard to construct examples of infinite sums 
$x= \sum_i x_i$ that are in the Hilbert space completion of $H^\cdot 
(\X)$ but such that $U_a(x)$ is no longer contained in this space. 
 
On the other hand, the 
index $j$ in the complex $K^{i,j}$ varies subject to the constraint 
$j+1=q$, where $q$ is the degree of the differential forms (\cf Remark 
\ref{rem1}). Thus, expressions of the form 
\begin{equation} \label{sigma2b} 
 \sigma\left\{ \left(\begin{array}{cc} 1 & 0 \\ 0 & 1 
\end{array} \right), \left(\begin{array}{cc} b & 0 \\ 0 & b^{-1} 
\end{array} \right) \right\} (x) = b^j x \ \ \ \text{ for $x\in 
K^{i,j}$} 
\end{equation} 
give rise to {\em bounded operators}. Thus the representation 
$\sigma_2$ of $\SL(2,\R)$ determines an action of the real group ring
by bounded operators in ${\mathcal B}(\H^\cdot (K,d))$. 
 
(2) It is sufficient to compute explicitly the following commutators with
the operator $\tilde \Phi$. Elements of the form \eqref{rep1} commute with 
$\tilde\Phi$. Moreover, we have: 
$$ [ N, \tilde\Phi ]\,(x) = \frac{1}{(2\pi \sqrt{-1})} 
\left(\frac{(1+j-i)}{2} - \frac{(1+j-i+2)}{2} \right)\, x 
=-N(x) , $$ 
$$ [ \sigma_1 (w_1), \tilde\Phi ] (x) = \left(\frac{(1+j-i)}{2} - 
\frac{(1+j+i)}{2} \right)\, \sigma_1 (w_1)(x) = -i \sigma_1 
(w_1)(x). $$ 
$$ [ \ell, \tilde\Phi ]\, (x) = \left(\frac{(1+j-i)}{2} - 
\frac{(1+j+2-i)}{2} \right) \, (2\pi \sqrt{-1})^{-1} x \wedge \omega 
= - \ell(x) $$ 
and 
$$ [ \sigma_2 (w_2), \tilde\Phi ] (x) = \left(\frac{(1+j-i)}{2} - 
\frac{(1-j-i)}{2} \right)\, \sigma_2 (w)(x) = j \, \sigma_2 
(w)(x). $$ 
In particular, it follows that all the commutators that arise from the 
right representation are bounded operators (\cf Remark 
\ref{rem1}). 
 
\noindent $\diamond$ 

The Lefschetz representation $\sigma_2$ of $\SL(2,\R)$ on the odd cohomology 
descends to a representation of $\PSL(2,\R)$: 

\begin{cor} The element $\sigma_2(-id)\in{\rm A}$ acts trivially on
the odd cohomology $\H^{2q+1}(K^\cdot,d)$. 
\label{sign} 
\end{cor} 
 
\noindent {\bf Proof.} For $x\in K^{i,j}$ we have 
$$ \sigma \left\{ 1, \left( \begin{array}{rr} -1 & 0 \\ 0 & -1 
\end{array}\right) \right\} (x) = (-1)^j x. $$ 
Since $j+1 =q$, where $q$ is the degree of the differential forms, we 
obtain that the induced action is trivial on odd cohomology. 
 
\noindent $\diamond$ 

\begin{rem} {\em The operator $\tilde\Phi$ in the data $({\rm A},
\H^\cdot (K,d), \tilde\Phi)$ of Theorem \ref{repres} does not yet
satisfy all the properties of a Dirac operator. In fact, the
eigenspaces of $\tilde\Phi$, which coincide with the graded pieces
$gr_{2p}^\w  H^q(\X)$ of the cohomology, are not finite dimensional as the 
condition on the resolvent in Definition \ref{specDef} would 
imply. Therefore, it is necessary to restrict the structure $({\rm 
A}, \H^\cdot (K,d),\tilde\Phi)$ to a suitable 
subspace of $\H^\cdot(K,d)$, which still carries all the 
arithmetic information.  } 
\end{rem} 
 
\begin{defn}\label{PhiCone} The operator $\Phi$ on 
$\H^\cdot({\rm Cone}(N)^\cdot) = H^\cdot(\XX)$ is obtained by extending 
the action on its graded pieces 
\begin{equation}\label{Phi-cone} 
\Phi|_{gr_{2p}^\w H^q(\XX)} :=\left\{ \begin{array}{ll} p & q\geq 2p \\ 
p-1 & q\leq 2p-1. \end{array} \right. 
\end{equation} 
according to the decomposition $H^q(\XX) = 
\oplus_{p\in\Z}gr^\w_{2p}H^q(\XX)$. 
\end{defn} 
 
Notice that this definition is compatible with the operator 
$\tilde\Phi$ defined in \eqref{PhiK}, acting on the complex 
$K^{\cdot,\cdot}$ and with the induced operator $\Phi$ on 
$\H^\cdot(K,d)=H^\cdot (\X)$. In fact, from the Wang exact sequence 
\eqref{Wang} and Corollary \ref{behv-N} we know that, for $q\geq 2p-1$ 
$gr_{2p}^\w H^q(\XX)$ is identified with a subspace of $gr_{2p}^\w 
H^q (\X)$, hence the restriction of the operator $\Phi$ on $H^\cdot 
(\X)$ acts on $gr_{2p}^\w H^q(\XX)$ as multiplication by $p$. In the 
case when $q\leq 2(p-1)$, again using the exact sequence \eqref{Wang} 
(\cf \eqref{sequ2-ex} Corollary \ref{behv-N}), we can define $\Phi$ of 
an element in $gr_{2p}^\w H^q(\XX)$ as $\Phi$ of a preimage in 
$gr_{2(p-1)}^\w H^{q-1}(\X)$, hence as multiplication by $p-1$. This 
is obviously well defined, hence the definition \eqref{Phi-cone} is 
compatible with the exact sequences and the duality 
isomorphisms. Moreover, the operator $\Phi$ of \eqref{Phi-cone} 
agrees with the operator \eqref{Phi} on the subspace 
$H^\cdot(\X)^{N=0}$ of $H^\cdot(\XX)$. 
 
\begin{thm} \label{spectral3} Consider a Riemann surface $X_{/\C}$ and
the hyper-cohomology $H^\cdot (\XX)$ of ${\rm Cone}(N)$.
The inner product \eqref{innerprod} defined by the polarization $\psi$
induces an inner product on $H^\cdot (\XX)$. Moreover, the representation 
\eqref{rho-preS3} of $\SL(2,\R)$ induces an action of the real group ring
by bounded operators on the real Hilbert space completion 
of $H^q(\XX)$. For ${\rm A}$ the image under $\rho$ of the group ring, 
consider the data $({\rm A},H^\cdot (\XX),\Phi)$, with $\Phi$ as in
\eqref{Phi-cone}. The operator $\Phi$ satisfies the properties of a
1--summable Dirac operator, with bounded commutators with the elements
of ${\rm A}$.
\end{thm} 
 
\noindent {\bf Proof.} The Wang exact sequence \eqref{Wang} and 
Corollary \ref{behv-N} imply that the hyper-cohomology $H^\cdot 
(\XX)$ of ${\rm Cone}(N)$ injects or is mapped upon surjectively 
by the hyper-cohomology $\H^\cdot(K,d)$ of the complex, in a way 
which is compatible with the grading. Thus, we obtain an induced 
inner product and Hilbert space completion on $H^\cdot (\XX)$. 
Consider $Ker(N)\subset gr^\w_{2p} H^q(\X)$. By Corollary 
\ref{behv-N}, we know that this is non-trivial only if $q\ge 2p$, 
and in that case it is given by $gr^\w_{2p} H^q(\XX)$. Thus, we 
can show that there is an induced representation on $\oplus_{2p\le 
q} gr^\w_{2p} H^q(\XX)$ by showing that the representation of 
${\rm A}(\tilde\Gamma)$ on $H^q(\X)$ preserves $Ker(N)$. 
 
In the definition of the complex $K^{i,j,k}$ in \eqref{complex}, the 
indices $i,j,k$ and the integers $p,q$ are related by 
$2p=j+1-i$ and $q=j+1$. Thus, the condition $q\ge 2p$ corresponds to 
$i\geq 0$. The representation $\sigma_2$ of $\SL(2,\R)$ on $K^{i,j}$ 
preserves the subspace with $i\geq 0$. Similarly, by construction, the 
representation $\sigma_2$ preserves the subspaces $\oplus_j K^{i,j,k}$ 
of $\oplus_{j,k\ge i \ge 0} K^{i,j,k}$. This implies that the induced 
representation $\sigma_2$ on $\H^\cdot(K,d)$ preserves 
the summands of $gr^\w_{2p}H^q(\X)$ as in Remark \ref{rem3}, and in 
particular it preserves ${\rm Ker}(N)$. Thus we obtain a representation 
$\rho_{Ker(N)}$ mapping the real group ring to ${\rm A}$
in ${\mathcal B}({\rm Ker}(N))$.
 
The duality isomorphisms $N^{q-2p}$ of Proposition \ref{N-isos} 
determine duality isomorphisms between pieces of the 
hyper-cohomology $H^\cdot(\XX)$ of the cone: 
\begin{equation}\label{isoms-cone} \begin{array}{lr} 
\delta_0 :  gr^\w_{2p} H^0(\XX) \stackrel{\simeq}{\to} gr^\w_{2r} 
H^1(\XX), & p\le 0, r=-p+1 \ge 1 \\[3mm] 
\delta_1 : gr^\w_{2p} H^1(\XX) \stackrel{\simeq}{\to} gr^\w_{2r} 
H^2(\XX), & p\le 0, r=-p+2 \ge 2 \\[3mm] 
\delta_2 : gr^\w_{2p} H^2(\XX) \stackrel{\simeq}{\to} gr^\w_{2r} 
H^3(\XX), & p\le 1, r=-p+3 \ge 2.\end{array} \end{equation} We set 
$\delta= \oplus_{q=0}^2 \delta_q$ and we obtain an action of ${\rm 
A }(\tilde\Gamma)$ on $H^\cdot(\XX)$ by extending the 
representation $\rho_{Ker(N)}$ by $\delta\circ \rho_{Ker(N)} \circ 
\delta^{-1}$ on the part of $H^\cdot(\XX)$ dual to ${\rm Ker}(N)$. 
 
The operator $\tilde \Phi$ of \eqref{PhiK} induces the operator $\Phi$ 
of \eqref{Phi-cone} on $\H^\cdot({\rm Cone}(N))=H^\cdot(\XX)$. 
This has the properties of a Dirac operator: the eigenspaces are all finite 
dimensional by the result of Proposition \ref{descr}, and the 
commutators are bounded by Theorem \ref{repres}. 
The spectrum of $\Phi$ is given by $\Z$ with constant multiplicities,
so that $\Phi^{-1}$ on the complement of the zero modes is an
infinitesimal of order one. 
 
\noindent $\diamond$ 
 
\medskip 
 
\begin{rem} {\em We make a few important comments about the data $({\rm
A},H^\cdot (\XX),\Phi)$ of Theorem \ref{spectral3}. 
Though for the purpose of our paper we only consider arithmetic
surfaces, the results of Theorems \ref{repres} and \ref{spectral3}
admit a generalization to higher dimensional arithmetic varieties.
Moreover, notice that the data give a simplified cohomological version
of a spectral 
triple encoding the full geometric data at arithmetic infinity, which
should incorporate the spectral triple for the Hodge--Dirac operator
on $X_{/\C}$. In our setting, we restrict to forms harmonic with
respect to the harmonic theory defined by $\Box$ on the complex
$K^{\cdot,\cdot}$ (\cf Theorem  
\ref{main thm} and Corollary \ref{main cor}) that are ``square 
integrable'' with respect to the inner product \eqref{innerprod} 
given by the polarization. This, together with the action of the Lefschetz
$\SL(2,\R)$, is sufficient to recover the alternating product of the
local factor (see Theorem \ref{Gamma-zeta3}). In a more refined
construction of a spectral triple, which induces the data $({\rm
A},H^\cdot (\XX),\Phi)$ in cohomology, the Hilbert space will consist of
$L^2$-differential forms, possibly with additional geometric data,
where a $C^*$-algebra representing the algebra of functions on a
``geometric space at arithmetic infinity'' will act. } 
\end{rem}

\subsection{Simultaneous uniformization}\label{unifsec} 

We begin by recalling the following elementary fact of hyperbolic 
geometry. 
Let $\Gamma$ be a Kleinian group acting on $\P^1(\C)$. Let 
$\Omega \subset \P^1(\C)$ be a $\Gamma$-invariant domain. A subset 
$\Omega_0 \subset \Omega$ is $\Gamma$-stable if, for every $\gamma\in 
\Gamma$, either $\gamma(\Omega_0) =\Omega_0$ or $\gamma(\Omega_0)\cap 
\Omega_0=\emptyset$. The $\Gamma$-stabilizer of $\Omega_0$ is the 
subgroup $\Gamma_0$ of those $\gamma \in \Gamma$ such that 
$\gamma(\Omega_0) =\Omega_0$. Let $\pi_\Gamma$ denote the quotient map 
$\pi_\Gamma : \Omega \to \Gamma \backslash \Omega$. 
 
\smallskip 
 
\begin{claim} {\em (\cf Theorem 6.3.3 of \cite{Beardon}). Let 
$\Omega_0 \subset \Omega$ be an open $\Gamma$-stable subdomain and 
let $\Gamma_0$ be the $\Gamma$-stabilizer of $\Omega_0$. Then the 
quotient map $\pi_\Gamma$ induces a conformal equivalence 
$$ \Gamma_0 \backslash \Omega_0 \simeq \pi_\Gamma (\Omega_0). $$ } 
\label{claim1} 
\end{claim} 
 
\smallskip 
 
If $\Gamma$ is a Kleinian group, a {\em quasi circle} for $\Gamma$ is 
a Jordan curve $C$ in $\P^1(\C)$ which is invariant under the action 
of $\Gamma$. In particular, such curve contains the limit set 
$\Lambda_\Gamma$. 
 
\medskip 
 
In the case of Schottky groups, the following theorem shows that 
Bowen's construction of a quasi--circle for $\Gamma$ (\cf \cite{Bo}) 
determines a pair of Fuchsian Schottky groups $G_1, G_2 \subset 
\PSL(2,\R)$ associated to $\Gamma \subset \PSL(2,\C)$. 
The theorem describes the simultaneous 
uniformization by $\tilde\Gamma$ of the two Riemann surfaces with 
boundary $X_i =G_i \backslash \H^2$, where $\tilde\Gamma$ is the
$\Gamma$-stabilizer of the connected components of
$\P^1(\C)\smallsetminus C$. 

\begin{thm} Let $\Gamma\subset \PSL(2,\C)$ be a Schottky group of rank 
$g\geq 2$. Then 
the following properties are satisfied: 
\begin{enumerate} 
\item There exists a quasi--circle $C$ for $\Gamma$. 
\item There is a collection of curves $\hat C$ on the compact Riemann 
surface $X_{/\C} = \Gamma \backslash \Omega_\Gamma$ such that 
$$ X_{/\C} = X_1 \cup_{\partial X_1= \hat C= \partial X_2} 
X_2, $$ 
where $X_i = G_i \backslash \H^2$ are Riemann surfaces with 
boundary, and the $G_i\subset \PSL(2,\R)$ are Fuchsian Schottky 
groups. The $G_i$ are isomorphic to $\tilde \Gamma\subset 
\PSL(2,\C)$, the $\Gamma$--stabilizer of the 
two connected components $\Omega_i$ of $\P^1(\C) \backslash C$. 
\end{enumerate} 
\label{Schottky-Bers} 
\end{thm} 
 
\noindent {\bf Proof.} {\em 1.} For the construction of a quasi--circle we 
proceed as in \cite{Bo}. The choice of a set of 
generators $\{ g_i \}_{i=1}^g$ for $\Gamma$ determines $2g$ Jordan curves 
$\gamma_i$, $i=1\ldots 2g$ in $\P^1(\C)$ with pairwise disjoint interiors 
$D_i$ such that, if we write $g_{i+g} = g_i^{-1}$ for $i=1\ldots g$, 
the fractional linear transformation $g_i$ maps the interior of 
$\gamma_i$ to the exterior of $\gamma_{i+g \mod 2g}$. Now fix a choice 
of $2g$ pairs of points $\rho_i^{\pm}$ on the curves $\gamma_i$ in 
such a way that $g_i$ maps the two points $\rho_i^{\pm}$ to the two points 
$\rho_{i+g \mod 2g}^{\mp}$. Choose a collection $C_0$ of pairwise 
disjoint oriented 
arcs in $\P^1(\C)$ with the property that they do not intersect the 
interior of the disks $D_i$. Also assume that the oriented boundary of 
$C_0$ as a 1-chain is given by $\partial C_0 = \sum_i \rho_i^+ - 
\sum_i \rho_i^-$. Then the curve 
\begin{equation}\label{q-circle} 
C := \Lambda_\Gamma \cup \bigcup_{\gamma\in \Gamma} \gamma C_0 
\end{equation} 
is a quasi--circle for $\Gamma$. 
 
{\em 2.} The image of the curves $\gamma_i$ in the quotient $X_{/\C}= 
\Gamma \backslash \Omega_\Gamma$ consists of $g$ closed curves, whose 
homology classes $a_i$, $i=1\ldots g$, span the kernel $Ker(I_*)$ of 
the map 
$I_* : H_1(X_{/\C},\Z) \to H_1(\mX_\Gamma,\Z)$ induced by the inclusion 
of $X_{/\C}$ as the boundary at infinity in the compactification of 
$\mX_\Gamma$. The image under the quotient map of the collection of 
points $\{ \gamma \rho_i^{\pm} \}_{\gamma \in \Gamma, i=1\ldots 2g}$ 
consists of two points on each curve $a_i$, and the image of $C \cap 
\Omega_\Gamma$ consists of a collection $\hat C$ of pairwise disjoint 
arcs on $X_{/\C}$ connecting these $2g$ points. By cutting the surface 
$X_{/\C}$ along $\hat C$ we obtain two surfaces $X_i$, $i=1,2$, with 
boundary $\partial X_i = \hat C$. 
 
Since $C$ is $\Gamma$-invariant, the two connected components 
$\Omega_i$, $i=1,2$, of $\P^1(\C) \backslash C$ are 
$\Gamma$-stable. Let $\Gamma_i$ denote the $\Gamma$-stabilizer of 
$\Omega_i$. Notice that $\Gamma_1=\Gamma_2$. In fact, suppose there is 
$\gamma\in \Gamma$ such that $\gamma\in \Gamma_1$ and $\gamma \notin 
\Gamma_2$. Then $\gamma(\P^1(\C)) \subset \Omega_1 \cup C$, so that 
the image $\gamma(\P^1(\C))$ is contractible in $\P^1(\C)$. This would 
imply that $\gamma$ has topological degree zero, but an orientation 
preserving fractional linear transformation has topological degree 
one. 
 
We denote by $\tilde \Gamma$ the $\Gamma$-stabilizer $\tilde \Gamma = 
\Gamma_1=\Gamma_2$. Since the components $\Omega_i$ are open 
subdomains of the $\Gamma$-invariant domain $\Omega_\Gamma$, Claim 
\ref{claim1} implies that the quotients $\tilde \Gamma \backslash \Omega_i$ 
are conformally equivalent to the image $\pi_\Gamma(\Omega_i) \subset 
X_{/\C}$. By the explicit description of the surfaces with boundary 
$X_i$, it is easy to see that $\pi_\Gamma(\Omega_i) =X_i$. 
 
The quasi-circle $C$ is a Jordan curve in $\P^1(\C)$, 
hence by the Riemann mapping theorem there exist conformal maps 
$\alpha_i$ of the two connected components $\Omega_i$ to the two 
hemispheres $U_i$ of $\P^1(\C)\smallsetminus \P^1(\R)$, 
\begin{equation}\label{alphai} 
 \alpha_i : \Omega_i \stackrel{\simeq}{\longrightarrow} U_i  \ \ \ \ \ 
U_1\cup U_2 = \P^1(\C)\smallsetminus \P^1(\R). 
\end{equation} 
 
Consider the two groups 
$$ G_i :=\{ \alpha_i \gamma \alpha_i^{-1}: \, \gamma \in \tilde \Gamma 
\}. $$ 
These are isomorphic as groups to $\tilde \Gamma$, $G_i \simeq 
\tilde \Gamma$. Moreover, the $G_i$ preserve the upper/lower 
hemisphere $U_i$, hence they are Fuchsian groups, $G_i \subset 
\PSL(2,\R)$. 
 
The conformal equivalence $\tilde \Gamma \backslash \Omega_i \simeq X_i$ 
implies that the $G_i$ provide the Fuchsian uniformization of $X_i = 
G_i \backslash \H^2$, where $\H^2$ is identified with the upper/lower 
hemisphere $U_i$ in $\P^1(\C)\smallsetminus \P^1(\R)$. 
 
The group $\tilde\Gamma \subset \Gamma$ is itself a discrete 
purely loxodromic subgroup of $\PSL(2,\C)$ isomorphic to a free 
group, hence a Schottky group, so that the $G_i$ are Fuchsian 
Schottky groups. 
 
\noindent $\diamond$ 
 
Let $X_{/\R}$ be an orthosymmetric smooth real algebraic curve. 
In this case, we can apply the following refinement of the result 
of Theorem \ref{Schottky-Bers}. We refer to \cite{Alling}, \cite{Sepp} 
for a proof. 
 
\begin{prop}\label{real-Bers} 
Let $X_{/\R}$ be a smooth real orthosymmetric algebraic curve of genus 
$g\geq 2$. Then the following holds. 
\begin{enumerate} 
\item $X_{/\R}$ has a Schottky uniformization such that the domain of 
discontinuity $\Omega_\Gamma \subset \P^1(\C)$ is symmetric with 
respect to $\P^1(\R)\subset \P^1(\C)$. 
\item The reflection about $\P^1(\R)$ gives an involution on 
$\Omega_\Gamma$ that induces the involution $\iota: X_{/\R} \to 
X_{/\R} $ of the real structure. 
\item The circle $\P^1(\R)\subset \P^1(\C)$ is a quasi-circle for the 
Schottky group $\Gamma$, such that the image in $X_{/\R}$ of $\P^1(\R) 
\cap \Omega_\Gamma$ is the fixed point set $X_\iota$ of the involution. 
\item The Schottky group $\Gamma$ is a Fuchsian Schottky group. 
\end{enumerate} 
\end{prop} 
 
The choice of a lifting $\Gamma \subset \SL(2,\C)$ of the Schottky 
group determines corresponding lifts of $\tilde\Gamma \subset 
\SL(2,\C)$ and $G_i \subset \SL(2,\R)$. 

\medskip

\begin{rem} \label{delta<1} {\em In \cite{Man}, the condition 
$\dim_H(\Lambda_\Gamma) 
<1$ on the limit set was necessary in order to ensure convergence of 
the Poincar\'e series that gives the abelian differentials on $X_{/\C}$, 
hence in order to express the Green function on $X_{/\C}$ in 
terms of geodesics in the handlebody $\mX_\Gamma$. 
Notice that this condition is 
satisfied for an orthosymmetric smooth 
real algebraic curve $X_{/\R}$, with the choice of Schottky 
uniformization described above, where the limit set $\Lambda_\Gamma$ 
is contained in the rectifiable circle $\P^1(\R)$. } 
\end{rem} 
 
\bigskip

The above results on simultaneous uniformization provide a way of
implementing the datum of the Schottky uniformization into the
cohomological spectral data of \S \ref{cohsp3}, by letting the pair
$G_1\times G_2$ of Fuchsian Schottky groups in $\PSL(2,\R)$ act via
the $\SL(2,\R)\times \SL(2,\R)$ representation of the Lefschetz
module. 

\begin{lem} Let $\sigma: \SL(2,\R)\times \SL(2,\R) \to {\rm Aut}(K)$ 
be the representation associated to the bigraded Lefschetz module 
structure on the complex $K^{\cdot,\cdot}$. Let $\Gamma \subset 
\SL(2,\C)$ be a Schottky group that determines a Schottky 
uniformization of $X_{/\C}$. Let $\tilde\Gamma$ is the corresponding 
lift to $\SL(2,\C)$ of the $\Gamma$-stabilizer of the components 
$\Omega_i$ in the complement of a quasi-circle $C$. 
Then $(K^{\cdot,\cdot}, N, \ell, \psi)$ 
carries a left and a right action of $\tilde\Gamma$, 
\begin{equation} \label{left} \sigma_1 (\gamma) := \sigma\{ \alpha_1 
\gamma \alpha_1^{-1}, 1 \}  \end{equation} 
\begin{equation} \label{right} 
\sigma_2 (\gamma) := \sigma\{ 1, \alpha_2 \gamma \alpha_2^{-1} \}, 
\end{equation} 
where $\alpha_i$ are the conformal maps \eqref{alphai} of 
$\Omega_i$ to the two hemispheres in $\P^1(\C)\smallsetminus \P^1(\R)$. 
\label{repGamma} 
\end{lem} 
 
\noindent {\bf Proof.} By Theorem \ref{Schottky-Bers} we obtain 
Fuchsian Schottky groups $G_i =\{ \alpha_i \gamma \alpha_i^{-1}  , \gamma\in 
\tilde\Gamma \}$ in $\SL(2,\R)$. We consider the 
restriction of the representation $\sigma : \SL(2,\R)\times \SL(2,\R) 
\to {\rm Aut}(K)$ to $G_1 \times \{ 1 \}$ and $\{ 1 \} \times G_2$ 
as in \eqref{left} and \eqref{right}. 
 
\noindent $\diamond$ 
 
We can then adapt the result of Theorems \ref{repres} and
\ref{spectral3} to the restriction of the representation
\eqref{rho-preS3} to the group ring $\R[\tilde\Gamma]$. We denote by 
${\rm A}(\tilde\Gamma)\subset {\rm A}$ the image of the group ring
$\R[\tilde\Gamma]$ under the representation $\rho$.

\begin{thm} \label{spectral3Gamma}
Let $(K^{\cdot,\cdot}, d, N, \ell, \psi)$  be 
the polarized bigraded Lefschetz module associated to a Riemann 
surface $X_{/\C}$ of genus $g\geq 2$, 
and let $\Gamma\subset \SL(2,\C)$ be a choice of Schottky 
uniformization for $X_{/\C}$. Let $\tilde\Gamma$ be a lift to 
$\SL(2,\C)$ of the $\Gamma$--stabilizer of the two connected 
components of $\P^1(\C)\setminus C$ as in Theorem 
\ref{Schottky-Bers}. Consider the representation
\begin{equation}\label{rho-preS3Gamma} \rho: \R[\tilde\Gamma]  \to 
{\mathcal B}(\H^\cdot (K,d)) 
\end{equation} 
induced by \eqref{rho-preS3}, and the corresponding representation
$$ \rho : \R[\tilde\Gamma] \to {\mathcal B}({\rm 
Ker}(N)). $$ 
Then the results of Theorems \ref{repres} and
\ref{spectral3} hold for the data $(\R[\tilde\Gamma],H^\cdot
(\XX),\Phi)$, with ${\rm A}(\tilde\Gamma)=\rho(\R[\tilde\Gamma])$ and
$\Phi$ as in \eqref{Phi-cone}.
\end{thm}

Heuristically, the algebra ${\rm A}(\Gamma)$ represents a non-commutative 
version of the hyperbolic handlebody. In fact, 
if $\Gamma \subset \PSL(2,\C)$ is a Schottky group, the group ring of 
$\Gamma$, viewed as a non-commutative space, carries the complete 
topological information on the handlebody, which is the classifying space
of $\Gamma$. 

\bigskip 
 
If $X$ is an arithmetic surface over ${\rm Spec}(O_\K)$, where $O_\K$ 
is the ring of integers of a number field $\K$ with $n=[\K:\Q]$, the 
above result can be applied at each of the $n$ Archimedean  primes, by 
choosing at each prime $\alpha: \K\hookrightarrow \C$ a Schottky 
uniformization of the corresponding Riemann surface $X_{\alpha 
(\K)}$. At the primes that correspond to the $r$ real embeddings, 
$X_{\alpha (\K)}$ has a real structure. 
 
We have the following version of Theorem \ref{spectral3Gamma} for the
case of a real algebraic curve.  
 
\begin{prop}\label{real-rep} 
Let $X$ be an arithmetic surface over ${\rm Spec}(O_\K)$, with the 
property that, at all the real Archimedean  primes, the Riemann surface 
$X_{\alpha(\K)}$ is an orthosymmetric smooth real algebraic curve of 
genus $g\geq 2$. Let $(K^{\cdot,\cdot}, d, N, \ell, \psi)$  be 
the polarized bigraded Lefschetz module associated to 
$X_{/\R}=X_{\alpha(\K)}$. 
Then the representation $\sigma_2$ extends to representations 
$$ \rho: \R[\Gamma] \to {\mathcal B}(\H^\cdot (K,d)), $$ 
$$ \rho_N: \R[\Gamma] \to {\mathcal B}(\H^\cdot ({\rm 
Cone}(N))) $$ 
with the properties as in Theorems \ref{repres} and \ref{spectral3}, 
where $\Gamma$ is the Fuchsian Schottky uniformization for $X_{/\R}$
of Proposition \ref{real-Bers}. 
\end{prop} 
 
\medskip 
 
\begin{rem} \label{tunnel} {\em 
In this paper, the choice of dealing with the case of the 
Schottky group in Theorem \ref{spectral3Gamma} is 
motivated by the geometric setting proposed by Manin \cite{Man}. 
However, it is clear that the argument given in Theorem 
\ref{repres} holds in greater generality. This suggests that the 
picture of Arakelov geometry at the Archimedean  places may be 
further enriched by considering {\em tunnelling} phenomena between 
different Archimedean  places - something like higher order 
correlation functions - where, instead of filling each Riemann 
surface $X_{\alpha(\K)}$ by a handlebody, one can consider more 
general hyperbolic 3-manifolds with different boundary components 
at different Archimedean  primes. We leave the investigation of 
such phenomena to future work. } 
\end{rem}

\bigskip

\subsection{Some zeta functions and determinants} 
 
The duality isomorphisms $N^{q-2p}$ of Proposition \ref{N-isos} and 
the induced isomorphisms $\delta_q$ of \eqref{isoms-cone} give some 
further structure to the spectral triple. 
 
Define subspaces $H^\pm(\XX)$ of $H^\cdot(\XX)$ in the following way: 
\begin{equation}\label{Hpm} 
\begin{array}{l} 
H^-(\XX)= \oplus_{p\leq 0} gr_{2p}^\w H^0(\XX) \oplus \oplus_{p\le 
0} gr^\w_{2p} H^1(\XX) \oplus \oplus_{p\le 1} gr^\w_{2p} H^2(\XX), 
\\[3mm] 
 H^+(\XX):= \oplus_{p\ge 1} gr^\w_{2p} H^1(\XX) \oplus \oplus_{p\ge 
2} gr^\w_{2p} H^2(\XX) \oplus \oplus_{p\ge 2} gr^\w_{2p} H^3(\XX), 
\end{array} 
\end{equation} 
Let $\delta=\oplus_{q=0}^2 \delta_q$ be the duality isomorphism of 
\eqref{isoms-cone} and set 
$$ \omega = \left(\begin{array}{cc} 0 & \delta^{-1} \\ \delta & 0 
\end{array} \right). $$ 
The map $\omega$ interchanges the subspaces $H^\pm(\XX)$. 
 
\begin{lem}\label{even-q} 
The map $\omega$ has the following properties: 
\begin{itemize} 
\item $\omega^2= id$, $\omega^* =\omega$. 
\item $[\omega,a]=0$, for all $a\in {\rm A}$. 
\item $(\Phi \omega + \omega \Phi)|_{H^q(\XX)} = q\cdot id$. 
\end{itemize} 
\end{lem} 
 
\noindent {\bf Proof.} By construction (\cf Theorem 
\ref{spectral3}) the action of ${\rm A}$ commutes 
with $\omega$. By \eqref{Phi-cone}, for $q\geq 2p$ the operator 
$\Phi$ on $gr_{2p}^\w H^q(\XX)$ acts as multiplication by $p$. The 
duality isomorphism, mapping $gr_{2p}^\w H^q(\XX)$ to 
$gr_{2(q-p+1)}^\w H^{q+1}(\XX)$ (\cf Proposition \ref{N-isos}), and 
$\Phi$ acts on $gr_{2(q-p+1)}^\w H^{q+1}(\XX)$ as multiplication by 
$(q-p)$. Thus, we obtain 
$$ (\Phi \omega + \omega \Phi)|_{gr_{2p}^\w H^q(\XX)} (x) = (q-p)\cdot x + 
p \cdot x = q \cdot x. $$ 
 
\noindent $\diamond$ 
 
\begin{rem} {\em Recall that a spectral triple $({\mathcal 
A},{\mathcal H},D)$ is called {\em even} if there is an operator 
$\omega$ such that $\omega^2= id$ and $\omega^* =\omega$; the 
commutator $[\omega,a]=0$, for all $a\in {\mathcal A}$ and the Dirac 
operator satisfies $D \omega + \omega D =0$. The 
conditions of Lemma \ref{even-q} provide a weaker version of this 
notion, depending on the degree of the cohomology $H^q(\XX)$. } 
\end{rem} 
 
Due to the presence of this further structure on the spectral triple, 
determined by the duality isomorphisms, in addition to the family of 
zeta functions \eqref{aDzeta}, \eqref{aDzetas}, we can consider zeta 
functions of the form 
\begin{equation} \label{aDzetasP} 
\zeta_{a,P_\pm \Phi}(s,z) := \sum_{\lambda\in {\rm Spec}(P_\pm \Phi)} 
\Tr(a\, \Pi(\lambda,P_\pm \Phi)) (s-\lambda)^{-z}, 
\end{equation} 
where $P_{\pm}$ are the projections on $H^\pm(\XX)$. 
In this setting we can now recover the $\Gamma$-factors. 
 
\begin{thm}\label{Gamma-zeta3} 
Consider $a=\sigma_2(-id)$ as an element in $Aut(K^{\cdot,\cdot})$, 
acting on $H^\cdot (\XX)$ via the induced representation (\cf Theorem 
\ref{spectral3}). Then the zeta function \eqref{aDzetasP} 
$$ \zeta_{a, P_- \Phi}(s,z):= \sum_{\lambda \in {\rm Spec}(P_- \Phi)} 
\Tr(a \Pi(\lambda, P_- \Phi)) (s-\lambda)^{-z} $$ 
satisfies 
\begin{equation}\label{zeta-L} 
\exp\left( -\frac{d}{dz} \zeta_{a,P_- \Phi/(2\pi)}(s/(2\pi),z)|_{z=0} 
\right)^{-1} = \frac{ L_\C (H^1(X_{/\C},\C),s) }{L_\C (H^0(X_{/\C},\C),s) 
\cdot L_\C (H^2(X_{/\C},\C),s) }. 
\end{equation} 
\end{thm} 
 
\noindent {\bf Proof.} 
Notice that we have $P_- \Phi = \Phi|_{H^\cdot(\X)^{N=0}}$. Moreover, 
recall that the element $a=\sigma_2(-id)$ acts as $(-1)^{q-1}$ on 
differential forms of degree $q$. We have 
$$ \exp\left( -\frac{d}{dz} \zeta_{a,P_- \Phi/(2\pi)}(s/(2\pi),z)|_{z=0} 
\right) =\prod_{q=0}^2 \exp\left( -\frac{d}{dz} \zeta_{a,P_- 
\frac{\Phi}{2\pi}|_{H^q(\XX)}}(s/(2\pi),z)|_{z=0} \right) $$ 
$$ = \prod_{q=0}^2 \exp\left( (-1)^q \frac{d}{dz} 
\zeta_{\Phi_q}(s/(2\pi),z)|_{z=0} \right), $$ 
where $\Phi_q=\Phi|_{H^q(\X)^{N=0}}$. The result then follows by 
Proposition \ref{infdet}. 
 
\noindent $\diamond$ 
 
\medskip

\bigskip 
 
We give a few more examples of computations with zeta functions related 
to the arithmetic spectral triple. 
 
\begin{ex} For $Re(s) >>0$, the zeta function \eqref{zetaD} of the 
Dirac operator $\Phi$ is given by 
\begin{equation} \label{zetaPhiformula} 
\zeta_{\Phi}(s)=\Tr(|\Phi|^{-s}) = (4g+4) \zeta(s) + 1+ \frac{1}{2^s}, 
\end{equation} 
where $\zeta(s)$ is the Riemann zeta function 
$$ \zeta(s) =\sum_{n\geq 1} \frac{1}{n^s}. $$ 
\label{ex1zeta} 
\end{ex} 
 
\noindent {\bf Proof.} We compute explicitly $\Tr(|\Phi|^{-s})$. 
On the complement of the zero modes, the operator $|\Phi|$ has 
eigenvalues the positive integers, and the corresponding eigenspaces 
$E_n(|\Phi|)$ are described as follows: 
\begin{equation} \begin{array}{llll} 
E_1(|\Phi|) =& gr_{-2} H^0(X^*) \oplus & & \\ 
             & gr_{-2} H^1(X^*) \oplus & & gr_{4} H^1(X^*) \oplus \\ 
             & gr_{-2} H^2(X^*) \oplus & gr_{2} H^2(X^*)\oplus & 
                                                gr_{4} H^2(X^*) \oplus\\ 
             & & gr_{4} H^3(X^*), 
 \end{array} \label{E1Phi} 
\end{equation} 
\begin{equation} \begin{array}{llll} 
E_2(|\Phi|) =& gr_{-4} H^0(X^*) \oplus & & \\ 
             & gr_{-4} H^1(X^*) \oplus & & gr_{6} H^1(X^*) \oplus \\ 
             & gr_{-4} H^2(X^*) \oplus & gr_{4} H^2(X^*)\oplus & 
                                             gr_{6} H^2(X^*) \oplus \\ 
             & & & gr_{6} H^3(X^*), \end{array} \label{E2Phi} 
\end{equation} 
\begin{equation} \begin{array}{lll} 
E_n(|\Phi|) =& gr_{-2n} H^0(X^*) \oplus & \\ 
             & gr_{-2n} H^1(X^*) \oplus & gr_{2(n+1)} H^1(X^*) \oplus \\ 
             & gr_{-2n} H^2(X^*) \oplus & gr_{2(n+1)} H^2(X^*) \oplus \\ 
             &  & gr_{2(n+1)} H^3(X^*), \end{array} \label{EnPhi} 
\end{equation} 
for $n\geq 3$. 
 
Using then the result of Proposition \ref{descr} 
to compute the dimension of these eigenspaces, we have 
$$ \begin{array}{llll} 
\dim gr_{2p} H^0(X^*) = 1 &  (p \le 0) && \\[2mm] 
\dim gr_{2p} H^1(X^*) = 2g & (p \le 0) &  \dim gr_{2p} H^1(X^*) = 1 & (p 
\ge 1) \\[2mm] 
\dim gr_{2p} H^2(X^*) = 1 &  (p \le 1)  & \dim gr_{2p} H^2(X^*) = 2g & 
(p \ge 2)  \\ 
& & \dim gr_{2} H^2(X^*) = 1 & \\[2mm] 
&& \dim gr_{2p} H^3(X^*) = 1 &  (p \ge 2) 
\end{array} $$ 
We obtain 
\begin{equation}\label{eigendim} 
\dim E_n(|\Phi|) =\left\{ \begin{array}{ll} (4g+4) & n\geq 3 \\ 
(4g+5) & n=1,2 \end{array} \right. 
\end{equation} 
This completes the calculation. Thus, we obtain 
$$ \Tr(|\Phi|^{-s}) =\sum_{n\geq 1} \dim E_n(|\Phi|) n^{-s} = 
(4g+5) (1 + \frac{1}{2^s}) + (4g+4)\sum_{n\geq 3} n^{-s} $$ 
$$ = (4g+4) \zeta(s) + 1 + \frac{1}{2^s}. $$ 
 
\noindent $\diamond$ 
 
As an immediate consequence of this calculation we obtain the 
volume determined by the Dirac operator $\Phi$. 
 
\begin{ex} The volume in the metric determined by $\Phi$ is given by 
$V=(4g+4)$. \label{ex2zeta} 
\end{ex} 
 
\noindent {\bf Proof.} We compute the volume using the residue formula 
\eqref{VolRes}. Recall that the Riemann zeta function has residue $1$ at 
$s=1$. In fact, the well known formula 
$$ \lim_{s\to 1+} (s-1) \zeta_\K (s) = 2^{r_1} (2\pi)^{r_2} |d|^{-1/2} h 
R w^{-1}, $$ 
holds for an arbitrary number field $\K$, with class number $h$, 
$w$ roots of unity, discriminant $d$ and regulator $R$, with $r_1$ and 
$r_2$ counting the embeddings of $\K$ into $\R$ and $\C$. Applied to 
$\K=\Q$ this yields the result. This implies that, for the zeta function 
computed in \eqref{zetaPhiformula}, we obtain 
\begin{equation} \label{Vzeta} 
 V = \Tr_\omega (|\Phi|^{-1}) = {\rm Res}_{s=1} \Tr(|\Phi|^{-s}) = 
(4g+4) {\rm Res}_{s=1} \zeta(s)= (4g+4). 
\end{equation} 
 
\noindent $\diamond$ 
 
Notice how, while the handlebody $\mX_{\tilde\Gamma}$ in its 
natural hyperbolic metric has infinite volume, the Dirac operator 
$\Phi$ induces on ${\rm A}(\tilde\Gamma)$, which is our 
non-commutative version of the handlebody, a metric of finite 
volume. This is an effect of letting $\R[\tilde\Gamma]$ act via the
Lefschetz $\SL(2,\R)$ representation. 
 
\medskip 
 
It is evident from the calculation of the eigenspaces $E_n(|\Phi|)$ in 
Example \ref{ex1zeta} that the Dirac operator $\Phi$ has a spectral 
asymmetry (\cf \cite{APS}). This corresponds to an 
eta invariant, which can be computed easily from the dimensions of the 
eigenspaces in Example \ref{ex1zeta}, as follows. 
 
\begin{ex} The eta function of the Dirac operator $\Phi$ is given by 
\begin{equation}\label{eta} 
\eta_{\Phi} (s) := \sum_{0\neq \lambda \in \Sp(\Phi)} {\rm 
sign}(\lambda) \frac{1}{|\lambda|^s} =1 + \frac{1}{2^s}, 
\end{equation} 
The eta invariant $\eta_{\Phi}(0)=2$, measuring the spectral 
asymmetry, is independent of $g$. 
\end{ex} 
 
\smallskip

\subsection{Zeta function of the special fiber and Reidemeister 
torsion} \label{ReiZsec} 
 
In this paragraph we show that the expression \eqref{zeta-L} of 
Theorem \ref{Gamma-zeta3} can be interpreted as a {\em Reidemeister 
torsion}, and it is related to a zeta function for the fiber at 
arithmetic infinity. 
 
We begin by giving the definition of a zeta function 
of the special fiber of a semistable fibration, which motivates the 
analogous notion at arithmetic infinity. 
 
Let $X$ be a regular, proper and flat scheme over 
$\Sp(\Lambda)$, for 
$\Lambda$ a discrete valuation ring with quotient field $K$ and 
finite residue field $k$. Assume 
that $X$ has geometrically reduced, connected and 
one-dimensional fibers. 
Let us denote by $\eta$ and $v$ \resp the generic and the 
closed point of $\Sp(\Lambda)$ and by $\bar\eta$ and $\bar v$ the 
corresponding geometric points. Assume that the special fiber $X_v$ 
of $X$ is a connected, effective  Cartier divisor with reduced 
normal crossings defined over $k = k(v)$. This degeneration 
is sometime referred to as a {\it semistable} fibration over 
$\Sp(\Lambda)$. 
 
Let $N_v$ denote the cardinality of $k$. Then, define the 
zeta-function of the special fiber 
$X_v$ as follows ($u$ is an indeterminate) 
\begin{equation}\label{Zu} 
Z_{X_v}(u) = \frac{P_1(u)}{P_0(u)P_2(u)},\qquad P_i(u) = 
\det (1-f^\ast u~|~H^i(X_{\bar\eta},\Q_\ell)^{I_{\bar v}}), 
\end{equation} 
where $f^\ast$ is the geometric Frobenius \ie the map induced by the 
Frobenius morphism $f: X_{\bar v} \to X_{\bar v}$ on the cohomological 
inertia-invariants at $\bar v$.\vspace{.1in} 
 
The polynomials $P_i(u)$ are closely related to the characteristic 
polynomials of the Frobenius $$ F_i(u) = 
\det (u\cdot 1-f^\ast~|~H^i(X_{\bar\eta},\Q_\ell)^{I_{\bar 
v}})$$ through the formula 
\begin{equation}\label{PandF} 
P_i(u) = u^{b_i}F_i(u^{-1}),\qquad b_i = {\rm degree}(F_i). 
\end{equation} 
 
The zeta function $Z_{X_v}(u)$ generalizes on a semistable fiber the 
description of the Hasse-Weil zeta function of a smooth, projective 
curve over a finite field. 
 
\medskip 
 
Based on this construction we make the following definition for the 
fiber at an Archimedean  prime of an arithmetic surface: 
\begin{equation}\label{ZY} 
 Z_\Phi (u):=\frac{P_1(u)}{P_0(u) P_2(u)}, 
\end{equation} 
where we set 
\begin{equation}\label{PqPhi} 
 P_q(u) := \det_\infty \left( \frac{1}{2\pi}-u \frac{\Phi_q}{2\pi} 
\right), 
\end{equation} 
with $\Phi_q = \Phi|_{H^q(\X)^{N=0}}$. 
 
\medskip 
 
In order to see how this is related to the result of Theorem 
\ref{Gamma-zeta3}, we recall briefly a simple observation of Milnor (\cf \S 3 
\cite{Mil}). Suppose given a finite complex $L$ and an infinite cyclic 
covering $\tilde L$, with $H_*(\tilde L,\kappa)$ finitely generated over 
the coefficient field $\kappa$. Let $h: \pi_1 L \to \kappa(s)$ be the 
composition of the homomorphism $\pi_1 L \to \Pi$ associated to the 
cover with the inclusion $\Pi \subset {\rm Units}(\kappa (s))$. The 
Reidemeister torsion for this covering is given (up to 
multiplication by a unit of $\kappa \Pi$) by the alternating product 
of the characteristic polynomials $F_q(s)$ of the $\kappa$--linear map 
$$ s_*: H_q(\tilde L,\kappa) \to H_q(\tilde L,\kappa), $$ 
\begin{equation}\label{prodPitau} 
\tau(s) \simeq F_0(s) F_1(s)^{-1} F_2(s) \cdots F_n(s)^{\pm 1}. 
\end{equation} 
Moreover, for a map $T: L \to L$, let $\zeta_T (u)$ be the Weil zeta 
$$ \zeta_T (u) =P_0(u)^{-1} P_1(u) P_2(u)^{-1} \cdots P_n(u)^{\pm 
1}, $$ 
where the polynomials $P_q(u)$ of the map $T_*$ are related to the 
characteristic polynomials $F_q(s)$ by \eqref{PandF} and $b_q$ are the 
$q$-the Betti number of 
the complex $L$. By analogy with \eqref{prodPitau}, Milnor writes the 
Reidemeister torsion $\tau_T(s)$ (up to 
multiplication by a unit) as 
$$ \tau_T(s):=F_0(s) F_1(s)^{-1} F_2(s)  \cdots F_n(s)^{\mp 1}, $$ 
where $F_q(s)$ are the characteristic polynomials of the map $T_*$. 
Then the relation between zeta function and Reidemeister torsion is 
given by: 
\begin{equation}\label{zeta-torsion} 
\zeta_T (s^{-1})\tau_T (s) = s^{\chi(L)}, 
\end{equation} 
where $\chi(L)$ is the Euler characteristic of $L$. 
 
\bigskip 
 
Similarly, we can derive the relation between the zeta function of the 
fiber at infinity defined as in \eqref{ZY} and the alternating product 
of Gamma factors in \eqref{zeta-L}. 
 
First notice that the expression \eqref{zeta-L} is of the form 
\eqref{prodPitau}. Namely, we write 
\begin{equation} \label{Tors-Y} 
\frac{ L_\C (H^1(X_{/\C},\C),s) }{L_\C (H^0(X_{/\C},\C),s) 
\cdot L_\C (H^2(X_{/\C},\C),s) } = \frac{F_0(s) \cdot F_2(s)}{F_1(s)}, 
\end{equation} 
where we set 
\begin{equation}\label{FqPhi} 
F_q (s) := \det_\infty \left( \frac{s}{2\pi}-\frac{\Phi_q}{2\pi} 
\right), 
\end{equation} 
with $\Phi_q = \Phi|_{H^q(\X)^{N=0}}$. 
For this reason we may regard \eqref{Tors-Y} as the Reidemeister 
torsion of the fiber at arithmetic infinity: 
\begin{equation}\label{tauY} 
\tau_\Phi (s):= \frac{F_0(s) \cdot F_2(s)}{F_1(s)}. 
\end{equation} 
 
\smallskip 
 
The relation between zeta function and Reidemeister torsion is then 
given as follows. 
 
\begin{prop} 
The zeta function $Z_\Phi$ of \eqref{ZY} and the Reidemeister torsion 
$\tau_\Phi$ of \eqref{tauY} are related by 
$$ Z_\Phi (s^{-1}) \tau_\Phi (s) = s^{g-2} e^{\chi\,  s\log s}, $$ 
with $g$ is the genus of the Riemann surface $X_{/\C}$ and 
$\chi=2-2g$ its Euler characteristic. 
\end{prop} 
 
\noindent {\bf Proof.} The result follows by direct calculation of 
the regularized determinants as in Section \ref{2.5}. Namely, we 
compute (in the case $q=0,1$) 
$$ P_q(u)= \det_\infty\left( \frac{1}{2\pi}-u \frac{\Phi_q}{2\pi} 
\right) = \exp\left( -b_q \frac{d}{dz} ( (2\pi)^z \sum_{n\ge 0} (1+u 
n)^{-z})|_{z=0} \right) $$ 
$$ = \exp\left(b_q \left(\log \Gamma\left(\frac{1}{u}\right)+ 
\frac{\log 2\pi}{u} + \frac{\log u}{2} + \frac{\log u}{u} \right) 
\right)  = u^{b_q/2} e^{-b_q \frac{\log u}{u}} (2\pi)^{-1/u} 
\Gamma(1/u), $$ 
where $b_q$ are the Betti numbers of $X_{/\C}$. 
The case $q=2$ is analogous, but for the presence of the $+1$ 
eigenvalue in the spectrum of $\Phi_2$, hence we obtain 
$$ P_2(u) =\exp\left( - b_2 \frac{d}{dz} \left( (2\pi)^z u^{-z} \zeta(1/u, 
z) - \left(\frac{1}{u}-1 \right)^{-z} \right)_{z=0} \right) $$ 
$$ = \Gamma_\C \left(\frac{1}{u}-1 \right)^{-1} 
u^{-3/2} e^{\frac{\log u}{u}}.  $$ 
Thus, we obtain 
$$ Z_\Phi (s^{-1}) =  \frac{L_\C (H^0(X_{/\C},\C),s) 
\cdot L_\C (H^2(X_{/\C},\C),s) }{ L_\C (H^1(X_{/\C},\C),s) } 
s^{g-2} e^{\chi \, s\log s}. $$

\noindent $\diamond$

\section{Shift operator and dynamics.}\label{6} 
 
In this section we consider an arithmetic surface $X$ over 
$\Sp(O_\K)$, with $\K$ a number field, and a fixed Archimedean  prime 
which corresponds to a {\em real} embedding $\alpha: \K\hookrightarrow 
\R$. We also assume that the corresponding Riemann surface 
$X_{/\R}$ is an orthosymmetric smooth real algebraic 
curve of genus $g\geq 2$. 
 
\smallskip 
 
We consider a Schottky uniformization of the Riemann 
surface $X_{/\R}$ and the hyperbolic filling given by 
the handlebody $\mX_\Gamma$ that has $X_{/\R}$ as the 
conformal boundary at infinity, $\mX_\Gamma \cup X_{/\R} =\Gamma 
\backslash (\H^3 \cup \Omega_\Gamma)$. 
 
\medskip 
 
Geodesics in $\mX_\Gamma$ can be lifted to geodesics in $\H^3$ with 
ends on $\P^1(\C)$. Among these, geodesics with one or both ends on 
$\Omega_\Gamma \subset \P^1(\C)$ correspond to geodesics in 
$\mX_\Gamma$ that reach the boundary at infinity $X_{/\R}= \Gamma 
\backslash \Omega_\Gamma$ in infinite time. The geodesics in $\H^3$ 
with both ends on $\Lambda_\Gamma \subset \P^1(\C)$ project in the 
quotient to geodesics contained in the {\em convex core} ${\mathfrak 
C}_\Gamma= \Gamma \backslash \text{ 
Hull }(\Lambda_\Gamma)$ of $\mX_\Gamma$. Since the Schottky group 
$\Gamma$ is geometrically finite, ${\mathfrak C}_\Gamma$ is a 
bounded  region inside $\mX_\Gamma$ (\cf 
\cite{MaTa}). For this reason, geodesics in $\mX_\Gamma$ that lift to 
geodesics in $\H^3$ with both ends on $\Lambda_\Gamma$ are called 
{\em bounded geodesics}. 
 
\medskip 
 
We denote by $\Xi \subset \mX_\Gamma$ the 
image under the quotient map of all geodesics in $\H^3$ with 
endpoints on $\Lambda_\Gamma\subset \P^1(\C)$, endowed with the 
induced topology. Similarly, we denote by $\Xi_c \subset \Xi$ the 
image in $\mX_\Gamma$ of all geodesics in $\H^3$ with endpoints of the 
form $\{ z^-(h), z^+(h) \}$ for some primitive $h\in \Gamma$. Here 
$z^{\pm}(h)$ are the attractive and repelling fixed points of $h$, and 
the element $h$ is primitive in $\Gamma$ if it is not a power of some 
other element of $\Gamma$. 
 
\begin{defn} We denote by $\tilde\Xi$ the orientation 
double cover of $\Xi$ and we refer to it as the {\em infinite tangle 
of bounded geodesics}. Similarly, we define $\tilde\Xi_c$ to be the 
orientation double cover of $\Xi_c$ and we refer to it as the {\em 
tangle of primitive closed geodesics}. 
\label{tangle-defn} 
\end{defn} 
 
Since $\Xi$ is orientable, $\tilde\Xi \cong \Xi \times \Z/2$, where 
the second coordinate is the choice of an orientation on each 
geodesic. So, for instance, the geodesics in $\H^3$ with endpoints $\{ 
z^-(h), z^+(h) \}$ and $\{ z^+(h),z^-(h) \}=\{ z^-(h^{-1}), 
z^+(h^{-1}) \}$ correspond to the same closed geodesic in $\Xi_c$, 
while in the double cover $\tilde\Xi_c$ they give rise to the two 
different lifts of the geodesic in $\Xi_c$. 
 
More precisely, let $L_{\{ a, b \}}$ denote the geodesic in $\H^3\cong 
\C\times \R^+$ with endpoints $\{ a, b \}$ in the complement of the 
diagonal in $\P^1(\C)\times \P^1(\C)$. A parameterization for the
geodesic $L_{\{ a, b \}}$ is given by
\begin{equation}\label{geodesicH3} 
\tilde L_{ \{ a,b \} }(s) = \left( \frac{a e^s + b 
e^{-s}}{e^s + e^{-s}}, \frac{|a-b|}{e^s + e^{-s}}\right) \ \ \ \ 
s\in \R .
\end{equation} 
The parameter $s$ in 
\eqref{geodesicH3} determines a parameterization by arc length on the 
corresponding geodesic $\pi_\Gamma( L_{\{ a, b \}} )$ in $\tilde\Xi$, 
where $\pi_\Gamma : \H^3 \to \mX_\Gamma$ is the quotient map. 
Then we have 
$$ \tilde\Xi= \{ \pi_\Gamma(\tilde L_{\{ a, b \}}(s)): \, \,  s\in \R, \, \, 
(a,b)\in (\Lambda_\Gamma \times \Lambda_\Gamma)^0 \}, $$ 
where 
$$ (\Lambda_\Gamma\times\Lambda_\Gamma)^0 := 
(\Lambda_\Gamma\times\Lambda_\Gamma)\smallsetminus \Delta $$ 
denotes the complement of the diagonal in 
$\Lambda_\Gamma\times\Lambda_\Gamma$. 
 
\begin{rem}\label{Z2} {\em The $\Z/2$ involution on 
$\tilde \Xi$ that has $\Xi$ as quotient corresponds to the involution 
that exchanges the two factors in $\Lambda_\Gamma \times \Lambda_\Gamma$.} 
\end{rem} 
 
\begin{rem}\label{Real-str} {\em Notice that most constructions 
and results presented in this section are topological in nature, hence 
they do not require any assumption on the conformal structure of the 
Riemann surface $X_{/\R}$. In particular they are not sensitive to 
whether $X$ is a complex or real smooth algebraic curve. However, 
when we refer to the results of the previous sections, by 
comparing the dynamical and Archimedean  cohomology in Theorem 
\ref{map-ar-dyn}, we need to know that the $\Z/2$ cover $\tilde \Xi$ 
is compatible with the conformal structure on the Riemann 
surface. This requires the presence of the real structure, so that complex 
conjugation $z\mapsto \bar z$ on $\P^1(\C)$, which induces 
the change of orientation on the geodesics $L_{\{a,b\}}$ in $\H^3$, 
determines an involution on the $\Gamma$--quotient.} 
\end{rem}

\medskip 
 
According to \cite{Man}, the tangle of bounded 
geodesics $\tilde\Xi$ provides a geometric realization of the {\em 
dual graph} ${\mathcal G}_\infty$ of the maximally degenerate closed 
fiber at arithmetic infinity. 
 
\smallskip 
 
In this section, we consider a resolution of $\tilde\Xi$ by the mapping 
torus (suspension flow) ${\mathcal S}_T$ of a dynamical system $T$ 
(\cf Proposition \ref{tangle-torus}), with a surjection ${\mathcal 
S}_T \to \tilde\Xi$. 
 
In Theorem \ref{thm-homologydyn} we give an explicit description of the 
cohomology $H^1({\mathcal S}_T)$. Such cohomology 
is endowed with a filtration whose 
graded pieces depend uniquely on the coding of geodesics in 
$\tilde\Xi$ (\cf Proposition \ref{homology-dyn} and Proposition 
\ref{H1pairing}).

\subsection{The limit set and the shift operator}\label{subsecshift} 
 
Given a choice of a set of generators $\{ g_i \}_{i=1}^g$ for the 
Schottky group $\Gamma$, there is a bijection between the elements of 
$\Gamma$ and the set of all {\em reduced words} in the $\{ g_i 
\}_{i=1}^{2g}$, where we use the notation $g_{i+g} := g_{i}^{-1}$, 
for $i=1,\ldots, g$. Here by reduced words we mean all finite 
sequences $w = a_0 \ldots a_\ell$ in the $g_i$, for any $\ell \in \N$, 
satisfying $a_{i+1} \neq a_i^{-1}$ for all $i=0\ldots, \ell-1$. 
 
We also consider the set ${\mathcal S}^+$ of all right--infinite 
{\em reduced sequences} in the $\{ g_i \}_{i=1}^{2g}$, 
\begin{equation}\label{S+} {\mathcal S}^+= \{ a_0 a_1 \ldots 
a_\ell \ldots  \, \, | a_i 
\in \{ g_i \}_{i=1}^{2g}, \, \, a_{i+1} \neq a_i^{-1}, \forall i\in \N \}, 
\end{equation} 
and the set ${\mathcal S}$ of {\em doubly infinite reduced sequences} 
in the $\{ g_i \}_{i=1}^{2g}$, 
\begin{equation}\label{S} {\mathcal S}= \{\ldots a_{-m} \ldots a_{-1} 
a_0 a_1 \ldots a_\ell 
\ldots  \, \, | a_i \in \{ g_i \}_{i=1}^{2g}, \, \, a_{i+1} \neq a_i^{-1}, 
\forall i\in \Z \}. \end{equation} 
 
On the space ${\mathcal S}$ we consider the topology generated by the 
sets $W^s(x,\ell) =\{ y\in {\mathcal S} | x_k = y_k, k\geq 
\ell \}$, and the $W^u(x,\ell)=\{ y\in {\mathcal S} | x_k = y_k, k\leq 
\ell \}$ for $x\in {\mathcal S}$ and $\ell \in \Z$. This induces a 
topology with analogous properties on ${\mathcal S}^+$ by realizing it as a 
subset of ${\mathcal S}$, for instance, by extending each sequence to 
the left as a constant sequence. 
 
\medskip 
 
We define the {\em one-sided shift operator} $T$ on ${\mathcal S}^+$ 
as the map 
\begin{equation}\label{shift+} 
T( a_0 a_1 a_2 \ldots a_\ell \ldots) = a_1 a_2 \ldots 
a_\ell \ldots 
\end{equation} 
 
We also define a {\em two-sided shift operator} $T$ on ${\mathcal 
S}$ as the map 
\begin{equation}\label{shift} 
\begin{array}{rcccccccccl} T(& \ldots & a_{-m} & \ldots & a_{-1}& a_0 
& a_1 & 
\ldots & a_{\ell} & \ldots &) = \\ & \ldots & a_{-m+1} & \ldots & a_{0} & 
a_1 & a_2 & \ldots & a_{\ell+1} & \ldots & \end{array} \end{equation} 
 
\medskip 
 
Given a choice of a base point $x_0 \in \H^3\cup \Omega_\Gamma$ 
we can define a map $Z: {\mathcal S}^+ \to \Lambda_\Gamma$ in the 
following way. For an eventually periodic sequence 
$w\overline{a_0\ldots a_N} \in {\mathcal S}^+$, with an initial word 
$w$, we set 
\begin{equation}\label{Z-map1} 
 Z(w\overline{a_0\ldots a_N}) = w z^+(a_0\ldots a_N). 
\end{equation} 
Here we identify the finite reduced word $w$ with an element in 
$\Gamma$, hence $w z^+(a_0\ldots a_N)$ is the image under $w\in \Gamma$ of 
the attractive fixed point of the element $a_0\ldots a_N$ of $\Gamma$. 
For sequences $a_0 \ldots a_\ell\ldots$ that are not eventually 
periodic, we set 
\begin{equation}\label{Z-map2} 
 Z(a_0 \ldots a_\ell\ldots) = \lim_{\ell \to \infty} (a_0 \ldots 
a_\ell) x_0, 
\end{equation} 
where again we identify a finite reduced word $a_0 \ldots 
a_\ell$ with an element in $\Gamma$. 
 
\smallskip 
 
We also introduce the following notation. We denote by ${\mathcal 
S}^+(w) \subset {\mathcal S}^+$ the set 
of right infinite reduced sequences in the $\{ g_i \}_{i=1}^{2g}$ that 
begin with an assigned word $w=a_0 \ldots a_\ell$. If $g\in \Gamma$ is 
expressed as a reduced word $w$ in the $g_i$, we write $\Lambda_\Gamma 
(g):= Z({\mathcal S}^+(w))$. 
 
\smallskip 
 
In the following we denote by $\Lambda_\Gamma \times_\Gamma 
\Lambda_\Gamma$ the quotient by the diagonal action of 
$\Gamma$ of the {\em complement of the diagonal} 
$(\Lambda_\Gamma\times\Lambda_\Gamma)^0$. 
 
\smallskip 
 
The group $\Gamma$ acts on $\P^1(\C)$ by fractional linear 
transformations, hence on $\Lambda_\Gamma$, which is a 
$\Gamma$-invariant subset of $\P^1(\C)$. The group $\Gamma$ also acts 
on ${\mathcal 
S}^+$: an element $\gamma\in \Gamma$, identified with a reduced word 
$\gamma=c_0\ldots c_k$ in the $g_i$, maps a sequence $a_0\ldots 
a_\ell\ldots$ to the sequence obtained from $c_0\ldots c_k 
a_0\ldots a_\ell\ldots$ by making the necessary cancellations that 
yield a reduced sequence. 
 
\medskip 
 
\begin{lem} 
The following properties are satisfied. 
\begin{enumerate} 
 
\item The spaces ${\mathcal S}^+$ and ${\mathcal S}$ are topologically 
Cantor sets. The one-sided shift $T$ of \eqref{shift+} is a 
continuous surjective map on ${\mathcal S}^+$, while the two-sided 
shift $T$ of \eqref{shift} is a homeomorphism of ${\mathcal S}$. 
 
\item The limit set $\Lambda_\Gamma$ with the topology induced by the 
embedding in $\P^1(\C)$ is also a Cantor set, and the map $Z$ of 
\eqref{Z-map1} and \eqref{Z-map2} is a homeomorphism. 
 
The shift operator $T$ on ${\mathcal S}^+$ induces the map $ZTZ^{-1}: 
\Lambda_\Gamma \to \Lambda_\Gamma$ of the form 
$$ ZTZ^{-1} |_{\Lambda_\Gamma (g_i)}(z) = g_i^{-1} (z). $$ 
 
\item  The map $Z$ is $\Gamma$-equivariant. 
 
\end{enumerate} 
\label{S-Lamb} 
\end{lem} 
 
\noindent{\bf Proof.} 
{\em 1.} The first claim can be verified easily. 
 
{\em 2.} It is not hard to see that the correspondence 
$a_0 \ldots a_\ell \ldots \mapsto \lim_{\ell \to \infty} (a_0 \ldots 
a_\ell) x_0$ gives a bijection between 
the complement of eventually periodic sequences in ${\mathcal S}^+$ 
and the complement of the fixed points $\{ z^-(h), z^+(h) \}_{h\in 
\Gamma}$ in $\Lambda_\Gamma$ (\cf \cite{Dalbo} Prop.~1.2). 
To see that the correspondence $w \overline{a_0\ldots a_N}\mapsto w\, 
z^+(a_0\ldots a_N)=z^+(w a_0\ldots a_N w^{-1})$ is a bijection between 
the set of eventually periodic sequences and the set of fixed points 
$\{ z^-(h), z^+(h) \}_{h\in \Gamma}$ in $\Lambda_\Gamma$, we proceed 
as in \cite{Dalbo}. For any $h\in \Gamma$, $h$ can be written as a 
reduced word $h=a_0\ldots 
a_\ell$ in the $g_i$'s. If this word satisfies $a_\ell \neq a_0^{-1}$, then 
$Z^{-1}(z^+(h))=\overline{a_0\ldots a_\ell}$ and 
$Z^{-1}(z^-(h))=\overline{a_\ell^{-1} \ldots a_0^{-1}}$. If $a_{\ell} 
=a_0^{-1}$, then there is an element $\gamma \in \Gamma$, such that 
$h= \gamma a_{i_k}\ldots a_{i_{k+N}} \gamma^{-1}$, with $a_{i_{k+N}} 
\neq a_{i_k}^{-1}$. In this case $Z^{-1}(z^+(h))= \gamma 
\overline{a_{i_k}\ldots a_{i_{k+N}}}$ and $Z^{-1}(z^-(h))=\gamma 
\overline{a_{i_{k+N}}^{-1}\ldots a_{i_{k}}^{-1}}$. As a continuous 
bijection from a compact to a Hausdorff space, $Z$ is a homeomorphism. 
The expression for $ZTZ^{-1}$ is then immediate. 
 
{\em 3.} By continuity it is sufficient to check that the map $Z$ 
restricted to the dense subset $$\{ z^-(h),z^+(h) \}_{h\in \Gamma} \subset 
\Lambda_\Gamma$$ is $\Gamma$-equivariant. This was proved already in 
{\em 2.}, since $Z(\gamma \overline{a_0\ldots a_N}) = \gamma\, 
z^+(a_0\ldots a_N)$. 
 
\noindent $\diamond$

We recall two notions that will be useful in the following. 
 
\smallskip 
 
Let $A=(A_{ij})$ be an $N\times N$ elementary matrix. The {\em subshift of 
finite type} with {\em transition matrix} 
$A$ is the subset of the set of all doubly infinite sequences in the 
alphabet $\{ 1,\ldots, N \}$ of the form 
\begin{equation} \label{SA} {\mathcal S}_A : = \{ \ldots i_{-m} \ldots 
i_{-1} i_0 i_1 \ldots i_\ell\ldots \, | 1\leq i_k \leq N, \, A_{i_k 
i_{k+1}}=1, \forall k\in \Z \}. \end{equation} 
A double sided shift operator of the form \eqref{shift} can be defined 
on any subshift of finite type. 
 
In the following we will consider the case where the elementary matrix 
$A$ is the symmetric $2g\times 2g$ matrix 
with $A_{ij}=0$ for $|i-j|=g$ and $A_{ij}=1$ otherwise. 
 
\begin{defn}\label{Smale-def} 
The pair $({\mathcal S},T)$ of a space and a homeomorphism is a {\em Smale 
space} if locally ${\mathcal S}$ can be decomposed as the product of 
expanding and contracting directions for $T$. Namely, the following 
properties are satisfied. 
\begin{enumerate} 
\item For every point $x\in {\mathcal S}$ there exist subsets $W^s(x)$ 
and $W^u(x)$ of ${\mathcal S}$, such that $W^s(x) \times  W^u(x)$ is 
homeomorphic to a neighborhood of $x$. 
\item The map $T$ is contracting on $W^s(x)$ and expanding on 
$W^u(x)$, and $W^s(Tx)$ and $T(W^s(x))$ agree in some neighborhood of 
$x$, and so do $W^u(Tx)$ and $T(W^u(x))$. 
\end{enumerate} 
\end{defn} 
 
\begin{lem} \label{S-Lamb2} 
The following properties are satisfied. 
\begin{enumerate} 
\item The map $Q: {\mathcal S} \to \Lambda_\Gamma \times \Lambda_\Gamma$ 
$$ Q(\ldots a_{-m}\ldots a_{-1} a_0 a_1 \ldots a_\ell \ldots ) = \left( 
Z(a_{-1}^{-1} a_{-2}^{-1} \ldots a_{-m}^{-1} \ldots) \, , \, Z(a_{0} a_1 
a_2 \ldots a_\ell \ldots) \right) $$ 
is an embedding of the space ${\mathcal S}$ in the Cartesian product 
$\Lambda_\Gamma \times \Lambda_\Gamma$. 
The image of the embedding $Q$ is given by 
${\rm Im} (Q) = \bigcup_{i\neq j} \Lambda_\Gamma(g_i) \times 
\Lambda_\Gamma(g_j)$. 
On ${\rm Im} (Q)$ the two-sided shift operator $T$ of \eqref{shift} 
induces the map $Q T Q^{-1}$ 
$$ \begin{array}{rllllcllll} Q T Q^{-1}|_{\Lambda_\Gamma (g_i)\times 
\Lambda_\Gamma(g_j)}( & Z(g_i b_1 b_2 \ldots b_m \ldots), &  & 
Z(g_j a_1 \ldots  a_\ell \ldots) &) = \\[2mm] 
( & Z(g_j^{-1} g_i b_1 \ldots b_m \ldots), & g_j^{-1} g_j & Z(a_1 a_2 
\ldots a_\ell \ldots) & ). \end{array} $$ 
 
\item The map $Q: {\mathcal S} \to \Lambda_\Gamma \times 
\Lambda_\Gamma$  descends to a homeomorphism of the quotients 
$$ \bar Q: {\mathcal S}/T \stackrel{\simeq}{\to} \Lambda_\Gamma 
\times_\Gamma \Lambda_\Gamma. $$ 
 
\item The space ${\mathcal S}$ can be identified with the subshift of 
finite type ${\mathcal S}_A$ with the symmetric $2g\times 2g$ matrix 
$A=(A_{ij})$ with $A_{ij}=0$ for $|i-j|=g$ and $A_{ij}=1$ otherwise. 
 
The two-sided shift operator $T$ on ${\mathcal S}$ of 
\eqref{shift} decomposes ${\mathcal S}$ in a product of expanding and 
contracting directions, so that $( {\mathcal S}, T)$ is a Smale 
space. 
\end{enumerate} 
\end{lem} 
 
\noindent{\bf Proof.} 
{\em 1.} The first claim follows directly from the definitions and 
Lemma \ref{S-Lamb}. 
 
{\em 2.} Notice that ${\rm Im} (Q)$ intersects each $\Gamma$ orbit in 
the complement of the diagonal $(\Lambda_\Gamma \times 
\Lambda_\Gamma)^0$. Moreover, if $(a,b) =Q(x)$, then for $\gamma \in 
\Gamma$ the element $(\gamma a, \gamma b)$ is in ${\rm Im} (Q)$ iff 
$b\in \Lambda_\Gamma(\gamma^{-1})$, and in that case $(\gamma a, 
\gamma b) = Q(T^n x)$, for $n=\text{ length }(\gamma)$ as a reduced 
word in the $g_i$. The statement on $\bar Q$ then follows easily. 
 
{\em 3.} We write ${\mathcal S}$ as the 
shift of finite type ${\mathcal S}_A$. Namely, we identify a reduced sequence 
$$\ldots a_{-m} \ldots a_{-1}\, a_0\, a_1 \ldots a_\ell \ldots$$ where each 
$a_k = g_{i_k} \in \{ g_i \}_{i=1}^{2g}$ with the sequence 
$\ldots i_{-m} \ldots i_{-1} i_0 i_1 \ldots i_\ell \ldots$ satisfying 
$A_{i_k i_{k+1}}=1$, for all $k\in \Z$. 
 
Then, for subshifts of finite type with the double sided shift 
\eqref{shift}, the sets 
$$ W^u(x) = \cup_{\ell \in \Z} W^u(x,\ell), $$ 
with 
$$ W^u(x,\ell) :=\{ y\in {\mathcal S} | x_k = y_k, k\leq \ell \} $$ 
and 
$$ W^s(x) = \cup_{\ell \in \Z} W^s(x,\ell). $$ 
with 
$$ W^s(x,\ell) := \{ y\in {\mathcal S} | x_k = y_k, k \geq \ell \} $$ 
give the expanding and contracting directions, so that 
$({\mathcal S},T)$  satisfies the properties of a Smale space (\cf 
\cite{Putn}). 
 
\noindent $\diamond$ 

The homeomorphism of Lemma \ref{S-Lamb2}.{\em 2} at the level of the
quotient spaces is sufficient for our purposes, but the identification
could be strengthened at the level of the groupoids of the equivalence
relations rather than on the quotients themselves, hence giving an
actual identification of noncommutative spaces. 
 
\begin{rem} \label{SmaleCalg} {\em It is well known that one can 
associate different ${\rm C}^*$--algebras to Smale spaces (\cf \cite{Rue2}, 
\cite{Putn}, \cite{PutSpi}). For the Smale space $({\mathcal S},T)$ we 
consider, there are four possibilities: the crossed product algebra 
$C({\mathcal S})\rtimes_T \Z$ and the ${\rm C}^*$--algebras ${\rm 
C}^*({\mathcal G}^s)\rtimes_T \Z$, ${\rm C}^*({\mathcal G}^u)\rtimes_T 
\Z$, ${\rm C}^*({\mathcal G}^a)\rtimes_T \Z$ obtained by considering 
the action of the shift $T$ on the groupoid ${\rm C}^*$--algebra (\cf 
\cite{Ren} for the definition of such algebra) associated to the 
groupoids ${\mathcal G}^s$, ${\mathcal G}^u$, ${\mathcal G}^a$ of the 
stable, unstable, and asymptotic equivalence relations on $({\mathcal 
S},T)$. } 
\end{rem} 
 
In the following we consider the algebras $C({\mathcal 
S})\rtimes_T \Z$ and ${\rm C}^*({\mathcal G}^u)\rtimes_T \Z$ and show 
that the first is a non-commutative space describing the 
quotient $\Lambda_\Gamma\times_\Gamma \Lambda_\Gamma$ and 
the second is a non-commutative space describing the quotient 
$\Lambda_\Gamma /\Gamma$.

\subsection{Coding of geodesics} \label{sectcoding}
 
Since the Schottky group $\Gamma$ is a free group consisting of only 
loxodromic elements, the coding of geodesics in $\mX_\Gamma$ in terms 
of the dynamical system $({\mathcal S},T)$ is particularly simple. The 
following facts are well known. We recall them briefly for 
convenience. 
 
We denote by ${\mathcal S}^p \subset {\mathcal S}$ the set of periodic 
reduced sequences in the $g_i$, \ie the set of periodic points of the 
shift $T$. We define 
$$ \hat \Xi:= \left\{ \pi_\Gamma( L_{\{ a, b \}} ) : (a,b)\in 
(\Lambda_\Gamma \times \Lambda_\Gamma)^0 
\right\}, $$ 
and 
$$ \hat \Xi_c := \left\{ \pi_\Gamma( L_{\{ z^+(h), z^-(h) \}} ) : h\in 
\Gamma\setminus id \right\}, $$ 
where $L_{\{a,b\}}$ denotes the geodesic in $\H^3\cong 
\C\times \R^+$ with endpoints $\{ a, b \}$.
 
The following Lemma gives a coding of the primitive closed geodesics 
in $\tilde\Xi_c$ by the quotient ${\mathcal S}^p /T$. 
 
\begin{lem} 
The correspondence 
${\mathcal L}_c: w\overline{a_0\ldots a_N} \mapsto \pi_\Gamma( L_{ 
\{ w\, z^+(a_0\ldots a_N), w\, z^-(a_0\ldots a_N) \} })$ 
induces a bijection between ${\mathcal S}^p /T$ and $\hat \Xi_c$ 
\label{iso-coding1} 
\end{lem} 
 
\noindent{\bf Proof.} Arguing as in Lemma \ref{S-Lamb}, we see that every 
closed geodesic in $\tilde\Xi_c$ is of the form 
$$\pi_\Gamma( L_{ \{ z^+(a_0\ldots a_N), z^-(a_0\ldots a_N) \} })$$ for 
some reduced sequence $a_0\ldots a_N$ with $a_N \neq 
a_0^{-1}$. If two elements 
$\overline{a_0,a_2,\ldots,a_N}$ and $\overline{b_0,b_1,\ldots,b_M}$ 
of ${\mathcal S}^p$ represent the same primitive closed geodesic, 
then the elements $h_a =a_0 a_2 \cdots a_N$ and $h_b = b_0 
b_1\cdots b_M$ are conjugate in $\Gamma$, $h_a= g h_b g^{-1}$, by an 
element $g=c_1 c_2 \cdots c_k$. It is easy to see that 
this implies $c_k= b_0^{-1}$, $c_{k-1}=b_1^{-1}$, etc.~ so 
that, for some $1\leq N_0 \leq N$, we have 
\[ 
\overline{b_0,b_1,\ldots,b_M}= T^{N_0} ( 
\overline{a_0,a_1,\ldots,a_N} ), 
\] 
that is, the two sequences are in the same equivalence class modulo 
the action of $T$. We refer to \cite{Dalbo} for details. 
 
\noindent $\diamond$ 
 
Via the map ${\mathcal L}_c$ of Lemma \ref{iso-coding1} we can 
define on $\hat \Xi_c$ a topology which makes it homeomorphic to the 
quotient ${\mathcal S}^p /T$. 
 
Similarly, we obtain a coding of geodesics in $\tilde\Xi$ by 
the quotient space ${\mathcal S}/T$. 
 
\begin{lem} 
The map ${\mathcal L}: {\mathcal S} \to \hat \Xi$, 
${\mathcal L}: x \mapsto \pi_\Gamma( L_{ \{ a,b \} })$, 
for $x\in {\mathcal S}$ and $(a,b)=Q(x)$ in $\Lambda_\Gamma \times 
\Lambda_\Gamma$, induces a bijection between ${\mathcal S}/T$ 
and $\hat \Xi$. 
\label{geod-coding} 
\end{lem} 
 
\noindent{\bf Proof.} 
Any geodesic in $\tilde\Xi$ lifts to a geodesic in 
$\H^3$ with ends $(a,b)$ in the complement of the diagonal 
$(\Lambda_\Gamma \times \Lambda_\Gamma)^0$. 
Notice that in $\H^3$ we have $\gamma L_{\{ a,b \}} = L_{\{ \gamma a, 
\gamma b \}}$. The claim then follows easily. 
 
\noindent $\diamond$ 
 
Via the map ${\mathcal L}$ of Lemma \ref{geod-coding} we can 
define on $\hat \Xi$ a topology which makes it homeomorphic to the 
quotient ${\mathcal S}/T$. 
 
\bigskip 
 
We introduce a topological space defined in terms of the Smale space 
$({\mathcal S},T)$, which we consider as a graph associated to the 
fiber at arithmetic infinity. This maps onto the dual graph 
$\tilde\Xi$ considered in \cite{Man}. 
 
\begin{defn} The mapping torus (suspension flow) of the dynamical 
system $( {\mathcal S}, T)$ is defined as 
\begin{equation} \label{suspensionT}  {\mathcal S}_T := {\mathcal S} 
\times [0,1] / (x,0)\sim (Tx,1) \end{equation} 
\end{defn} 
 
Consider the map $\tilde Q : {\mathcal S}_T \to \tilde\Xi$, defined by 
\begin{equation} \label{tildeQ} 
\tilde Q([x,t]) = \pi_\Gamma  \tilde L_{\{ a, b \}} (s(x,t)), 
\end{equation} 
where $(a,b)=Q(x)$ in $(\Lambda_\Gamma \times
\Lambda_\Gamma)^0$. Notice that the map $\tilde Q$ yields a geodesic
in the handlebody, but the parameterization induced by the time
coordinate $t$ on the mapping torus, in general, will not agree with
the natural parameterization of the geodesic by the arc length $s$. In
fact, the induced parameterization, here denoted by $s(x,t)$, has the
property that  the geodesic line $\tilde L_{\{ a, b \}} 
(s(x,t))$ in $\H^3$ crosses a fundamental domain for the action of 
$\Gamma$ in time $t\in [0,1]$. 
 
\begin{prop} 
The map $\tilde Q : {\mathcal S}_T \to \tilde\Xi$ of \eqref{tildeQ} is 
a continuous surjection. It is a bijection away from the intersection 
points of different geodesics in $\tilde\Xi$. 
\label{tangle-torus} 
\end{prop} 
 
\noindent{\bf Proof.} 
The parameterization $s(t)$ is chosen in such a way that the map 
\eqref{tildeQ} is well defined on equivalence classes. By {\em 2.} of 
Lemma \ref{S-Lamb2} and Lemma \ref{geod-coding} the map is a 
continuous surjection. 
 
\noindent $\diamond$

\subsection{Cohomology and homology of ${\mathcal S}_T$} \label{combinat} 
 
We give an explicit description of the cohomology 
$H^1({\mathcal S}_T,\Z)$. 
 
\begin{thm} 
The cohomology $H^1({\mathcal S}_T,\Z)$ satisfies the following 
properties. 
\begin{enumerate} 
\item There is an identification of $H^1({\mathcal S}_T,\Z)$ with the 
$K_0$-group of the crossed product ${\rm C}^*$-algebra for the 
action of $T$ on ${\mathcal S}$, 
\begin{equation}  H^1({\mathcal S}_T,\Z) \cong K_0({\rm C}({\mathcal 
S}) \rtimes_T \Z).  \label{H1K0} \end{equation} 
\item The identification \eqref{H1K0} endows $H^1({\mathcal S}_T,\Z)$ 
with a filtration by free abelian groups $F_0 
\hookrightarrow F_1 \hookrightarrow \cdots F_n\hookrightarrow \cdots$, 
with ${\rm rank} F_0=2g$ and 
${\rm rank} F_n = 2g(2g-1)^{n-1}(2g-2) +1$, for $n\geq 1$, so that 
$$ H^1({\mathcal S}_T,\Z)= \varinjlim_{n} F_n. $$ 
\end{enumerate} 
\label{thm-homologydyn} 
\end{thm} 
 
\noindent{\bf Proof.} {\em 1.} The shift $T$ acting on ${\mathcal S}$ 
induces an automorphism of the 
$C^*$--algebra of continuous functions $C({\mathcal S})$. With an 
abuse of notation we still denote it by $T$. Consider the crossed 
product $C^*$--algebra $C({\mathcal S}) \rtimes_T \Z$. This is 
a suitable norm completion of $C({\mathcal S}) [T, T^{-1}]$ with 
product $(V*W)_k = \sum_{r\in\Z} V_k \cdot (T^r W_{r+k})$, for 
$V=\sum_k V_k T^k$, $W=\sum_k W_k T^k$, and $V*W=\sum_k (V*W)_k T^k$. 
 
The $K$--theory group $K_0({\rm C}({\mathcal 
S}) \rtimes_T \Z)$ is described by the co--invariants of the 
action of $T$ (\cf \cite{BoHa},\cite{PaTu}). Namely, 
let ${\rm C}({\mathcal S},\Z)$ be the set of continuous functions 
from ${\mathcal S}$ to the integers. This is an abelian group 
generated by characteristic functions of clopen sets 
of ${\mathcal S}$. The {\it invariants} and {\it 
co--invariants}  are given respectively by ${\rm C}({\mathcal 
S},\Z)^T := \{ f \in {\rm C}({\mathcal S},\Z) \, | f- f\circ T =0 
\}$ and ${\rm C}({\mathcal S},\Z)_T := {\rm C}({\mathcal S},\Z) / 
{\rm B}({\mathcal S},\Z)$, with ${\rm B}({\mathcal S},\Z) := \{ f- 
f\circ T \, | f \in {\rm C}({\mathcal S},\Z) \}$. 
We have the following result (\cf \cite{BoHa}): 
\begin{itemize} 
\item The $C^*$--algebra ${\rm C}({\mathcal S})$ is a commutative 
AF--algebra (approximately finite dimensional), obtained as the direct 
limit of the finite dimensional commutative $C^*$--algebras generated 
by characteristic functions of a covering of ${\mathcal S}$. Thus, 
$K_0({\rm C}({\mathcal S})) \cong {\rm C}({\mathcal 
S},\Z)$, being the direct limit of the $K_0$-groups of the finite 
dimensional commutative $C^*$--algebras, 
and $K_1({\rm C}({\mathcal S}))=0$ for the same reason. 
\item The Pimsner--Voiculescu exact sequence (\cf \cite{PV}) then 
becomes of the form 
\begin{equation}\label{PV} 0 \to K_1({\rm C}({\mathcal S})\rtimes_T 
\Z) \to {\rm C}({\mathcal S},\Z) \stackrel{I-T_*}{\to} {\rm 
C}({\mathcal S},\Z) \to K_0({\rm C}({\mathcal S})\rtimes_T \Z) \to 
0, \end{equation} 
with $K_0({\rm C}({\mathcal S})\rtimes_T \Z) \cong {\rm C}({\mathcal 
S},\Z)_T$. Since the shift $T$ is {\it topologically transitive}, 
i.e. it has a dense orbit, we also have $K_1({\rm C}({\mathcal 
S})\rtimes_T \Z) \cong {\rm C}({\mathcal S},\Z)^T \cong \Z$. 
\end{itemize} 
 
Now consider the cohomology group $H^1({\mathcal S}_T,\Z)$. Via the 
identification with \v{C}ech cohomology, we can identify 
$H^1({\mathcal S}_T,\Z)$ with the group of homotopy classes 
of continuous maps of ${\mathcal S}_T$ to the circle. The isomorphism 
\begin{equation}\label{isomCTH1} 
{\rm C}({\mathcal S},\Z)_T \cong H^1({\mathcal S}_T,\Z) 
\end{equation} 
is then given explicitly by 
\begin{equation}\label{isomCTH1-map} 
f \mapsto [ \exp(2\pi i t f(x)) ], 
\end{equation} 
for $f\in {\rm C}({\mathcal S},\Z)$ and with 
$[\cdot ]$ the homotopy class. The map is well defined on the 
equivalence class of $f$ mod $B({\mathcal S},\Z)$ since, for an element 
$f -h+h\circ T$ the function 
$$ \exp(2\pi i t (f-h+h\circ T)(x)) =\exp(2\pi i tf(x)) 
\exp(2\pi i ((1-t)h(x) + t h(T(x)))), $$ 
since $h$ is integer valued, but $\exp(2\pi i ((1-t)h(x) + t 
h(T(x))))$ is homotopic to the constant function equal to $1$. 
It is not hard to see that \eqref{isomCTH1-map} gives the desired 
isomorphism \eqref{isomCTH1} (\cf \S 4-5 \cite{BoHa}). 
This proves the first statement. 
 
\medskip 
 
{\em 2.} There is a filtration on the set of coinvariants ${\rm 
C}({\mathcal S},\Z)_T$ (\cf Theorem 19 \S 4 of \cite{PaTu}). It is 
obtained in the following way. First, it is possible to identify 
\begin{equation}\label{CT-P} 
 {\rm C}({\mathcal S},\Z)_T ={\rm C}({\mathcal S},\Z)/B({\mathcal 
S},\Z) \cong {\mathcal P}/\delta {\mathcal P}, \end{equation} 
where ${\mathcal P} \subset {\rm C}({\mathcal S},\Z)$ is the set of 
functions that depend only on future coordinates, and $\delta$ is the 
operator $\delta(f) = f- f\circ T$. In fact, since 
characteristic functions of clopen sets in ${\mathcal S}$ depend only 
on finitely many coordinates, any function in 
${\rm C}({\mathcal S},\Z)$, when composed with a sufficiently high 
power of $T$ becomes a function only of the future coordinates, \ie of 
${\mathcal S}^+$, the set of (right) infinite reduced 
sequences in the generators of $\Gamma$ and their inverses, $\{ g_i 
\}_{i=1}^{2g}$. Since all the $f\circ T^k$, $k\geq 0$, 
define the same equivalence class in ${\rm C}({\mathcal S},\Z)_T$, we 
have the identification of \eqref{CT-P}. 
 
Then ${\mathcal P}$ can be identified with ${\rm 
C}({\mathcal S}^+, \Z)$ viewed as the submodule of the $\Z$-module 
${\rm C}({\mathcal S},\Z)$ of functions that only depend on future 
coordinates. As such, it is generated by characteristic 
functions of clopen subsets of ${\mathcal S}^+$. A basis of clopen 
sets for the topology of ${\mathcal S}^+$ is given by the sets 
${\mathcal S}^+(w)\subset {\mathcal S}^+$, where $w=a_0\ldots a_N$ is 
a reduced word in the $g_i$ and ${\mathcal S}^+(w)$ is the set of 
reduced right infinite sequences $b_0 b_{1} b_{2}\ldots b_{n}\ldots $ 
such that $b_{k} = a_k$ for $k=0\ldots N$. Thus, ${\mathcal P}$ has a 
filtration ${\mathcal P}=\cup_{n=0}^\infty {\mathcal P}_n$, where 
${\mathcal P}_n$ is generated by the characteristic functions of 
${\mathcal S}^+(w)$ with $w$ of length at most $n+1$. Taking into 
account the relations between these, we obtain that ${\mathcal P}_n$ 
is a free abelian group generated by the characteristic functions of 
${\mathcal S}^+(w)$ with $w$ of length exactly $n+1$. The number of such 
words is $2g(2g-1)^n$, hence ${\rm rank} {\mathcal P}_n = 2g(2g-1)^n$. 
 
The map $\delta$ satisfies $\delta : {\mathcal P}_n \to {\mathcal 
P}_{n+1}$, with a 1-dimensional kernel given by the constant 
functions. More precisely, if we write $f(a_0...a_n)$ for a 
function in ${\mathcal P}_n$, then 
$$ (\delta f)(a_0\ldots a_n a_{n+1}) = f(a_0\ldots a_n) - f(a_1\ldots 
a_{n+1}). $$ 
The resulting quotients 
$$ F_n = {\mathcal P}_{n} / \delta {\mathcal P}_{n-1} $$ 
are torsion free (\cf Theorem 19 \S 4 of \cite{PaTu}) and have ranks 
$$ {\rm rank} F_n = 2g(2g-1)^{n-1}(2g-2) +1 $$ 
for $n\geq 1$, while $F_0 \cong {\mathcal P}_0$ is of rank $2g$. 
There is an injection $F_n \hookrightarrow F_{n+1}$ induced by the 
inclusion ${\mathcal P}_{n}\subset {\mathcal P}_{n+1}$, and 
${\mathcal P}/\delta {\mathcal P}$ is the direct limit of the $F_n$ 
under these inclusions. Thus we obtain the filtration on 
$H^1({\mathcal S}_T,\Z)$: 
$$ H^1({\mathcal S}_T,\Z)= \varinjlim_{n} F_n. $$ 
This proves the second statement. 
 
\noindent $\diamond$ 
 
\begin{rem} {\em There is an interesting degree shift between 
$K$--group $K_0({\rm C}({\mathcal S})\rtimes_T \Z)$ of the crossed 
product algebras associated to the 
Smale space $({\mathcal S},T)$ and the cohomology $H^1({\mathcal 
S}_T,\Z)$. This degree shift is a general phenomenon related to the
Thom isomorphism \eqref{Thomiso} as we shall discuss in \S
\ref{modinftysect} (\cf \cite{Co-Thom}, \cite{Connes-tr}). }
\end{rem}

The following result computes the first homology of ${\mathcal S}_T$.

\begin{prop} 
The homology group $H_1({\mathcal S}_T,\Z)$ has 
a filtration by free abelian groups ${\mathcal K}_N$, 
\begin{equation}\label{hom-lim-dir} 
 H_1({\mathcal S}_T,\Z) = \varinjlim_N {\mathcal K}_N, 
\end{equation} 
with 
$$ K_N = {\rm rank}( {\mathcal K}_N ) = \left\{ \begin{array}{lr} 
(2g-1)^N +1 & N \text{ even } \\ (2g-1)^N + (2g-1) & N \text{ odd } 
\end{array}\right. $$ 
The group $H_1({\mathcal S}_T,\Z)$ can also be written as 
$$ H_1({\mathcal S}_T,\Z) = \oplus_{N=0}^\infty {\mathcal R}_N $$ 
where ${\mathcal R}_n$ is a free abelian group of ranks $R_1=2g$ and 
$$ 
 R_N = {\rm rank}({\mathcal R}_N) =\frac{1}{N} \sum_{d|N} \mu(d)\, 
(2g-1)^{N/d}, 
$$ 
for $N>1$, with $\mu$ the M\"obius function. This is isomorphic to the 
free abelian group on countably many generators, $\Z\langle {\mathcal 
S}^p /T \rangle$. 
\label{homology-dyn} 
\end{prop} 
 
\noindent{\bf Proof.} The lines $\{ [x,t]: t\in \R \}$ are pairwise 
disjoint in the mapping torus ${\mathcal S}_T$, hence $H_1({\mathcal 
S}_T,\Z)$ is generated by {\em closed curves} $\{ 
[x,t]: t\in \R \}$. By construction, a closed curve in ${\mathcal 
S}_T$ corresponds to a point $x\in {\mathcal S}$ such that $T^N x =x$ 
for some $N\geq 1$. More precisely, if $x= \overline{ a_0 \ldots a_N 
}$ is a doubly infinite periodic sequence in ${\mathcal S}$, obtained 
by repeating the word $a_0 \ldots a_N$, so that 
$T^Nx=x$, then the map $c_x: S^1 \to {\mathcal S}_T$ of the 
form 
\begin{equation}\label{c-x} 
c_x : e^{2\pi i t} \mapsto [(x,Nt)],  \  \  \   T^N x= x. 
\end{equation} 
defines a closed curve in ${\mathcal S}_T$, which is a non-trivial 
homology class in $H_1({\mathcal S}_T,\Z)$. Since all non-trivial 
homology classes can be obtained this way, the homology $H_1({\mathcal 
S}_T,\Z)$ is generated by all $c_x$ as in \eqref{c-x}, for $x$ a 
reduced word $a_0\ldots a_N$ such that $a_N \neq a_0^{-1}$. 
 
Let ${\mathcal K}_N$ be the free 
abelian groups generated by all the reduced words $a_0 \ldots a_N$ of 
length $N+1$ satisfying $a_N \neq a_0^{-1}$. When we identify the 
elements of the ${\mathcal K}_N$ with homology classes via 
\eqref{c-x}, we introduce relations between the ${\mathcal K}_N$ for 
different $N$, namely we have embeddings $k {\mathcal 
K}_N \hookrightarrow {\mathcal K}_{kN}$, 
\begin{equation}\label{N-kN} 
 k \langle a_0 \ldots a_N \rangle \mapsto  \, \, \langle 
\underbrace{a_0 \ldots a_N a_0 \ldots a_N\ldots a_0 \ldots 
a_N}_{k-times} \rangle. 
\end{equation} 
Thus, the homology $H_1({\mathcal S}_T,\Z)$ is computed as the limit 
$H_1({\mathcal S}_T,\Z)=\lim_N {\mathcal K}_N$ with respect to the 
maps $J_k: {\mathcal K}_N \longrightarrow {\mathcal K}_{kN}$ that send 
the $c_x$ of \eqref{c-x} to the composite 
$$ S^1\stackrel{z\mapsto z^k}{\longrightarrow} S^1 \stackrel{c_x}{\to} 
{\mathcal S}_T. $$ 
 
Thus, the homology $H_1({\mathcal S}_T,\Z)$ can be identified 
with the $\Z$-module 
generated by the elements of ${\mathcal S}^p /T$. This quotient can be 
written as a 
disjoint union ${\mathcal S}^p /T =\cup_{n=0}^\infty {\mathcal S}^p_n 
/T$, where 
${\mathcal S}^p_n \subset {\mathcal S}^p$ is the subset of (primitive) 
periodic sequences with period of length $n+1$. This gives the 
description of $H_1({\mathcal S}_T,\Z)$ as direct sum of the 
${\mathcal R}_n = 
\Z \langle {\mathcal S}^p_n /T \rangle$. 
 
The computation of the ranks of the ${\mathcal K}_N$ and of the 
${\mathcal R}_n$ is obtained as follows. 
 
For a fixed $a_0$, denote by $p(N,k)$ the number of reduced sequences 
$a_0\ldots a_N$, such that the last $k$ terms are all equal to 
$a_0^{-1}$. Then $2g \cdot p(N,0)= {\rm rank} {\mathcal 
K}_N$. Moreover, since the total number of all reduced sequences of 
length $N+1$ is $2g(2g-1)^N$, we have $\sum_{k=0}^N p(N,k)= (2g-1)^N$. 
It is not hard to see from the definition that the $p(N,k)$ 
satisfy $p(N,k)=p(N+1,k+1)$, $p(N,N)=0$, $p(1,0)=2g-1$, 
$p(N,0)= p(N+1,1)+p(N-1,0)$. 
 
The calculation of $p(N,0)$ then follows inductively, using 
$p(N,1)=p(N-1,0) - p(N-2,0)$, and the sum $\sum_{k=0}^N p(N,k)= 
(2g-1)^N$, where $p(N,k)=p(N-k+1,1)$. 
 
The rank of the ${\mathcal R}_N$ can be computed by first considering 
that $K_N = \sum_{d|N} d\, R_d$, since the total number of reduced 
sequences $a_0\ldots a_N$ with $a_N\neq a_0^{-1}$ is the sum of the 
cardinalities of the ${\mathcal S}^p_d$ over all $d$ dividing 
$N$. These satisfy $\#{\mathcal S}^p_d =d\, R_d$. Then we obtain 
$$ R_N = \frac{1}{N}\sum_{d|N} \mu(d) K_{N/d}=\frac{1}{N}\sum_{d|N} 
\mu(d) (2g-1)^{N/d}, $$ 
using the M\"obius inversion formula in the first equality and the 
fact that $\sum_{d|N}\mu(d) =\delta_{N,1}$ in the second. 
 
\noindent $\diamond$ 
 
Combining Theorem \ref{thm-homologydyn} with Proposition 
\ref{homology-dyn}, we can compute explicitly the pairing of homology 
and cohomology for ${\mathcal S}_T$. This relates the 
filtration by the $F_n$ on the cohomology of ${\mathcal S}_T$ to the 
coding of closed geodesics in $\tilde\Xi_c$. 
 
\begin{prop} 
Let $F_n$ and ${\mathcal K}_N$ be the filtrations defined, 
respectively, in Theorem \ref{thm-homologydyn} and Proposition 
\ref{homology-dyn}. There is a pairing 
\begin{equation}\label{pairing1} 
\langle \cdot, \cdot \rangle : F_n \times 
{\mathcal K}_N \to \Z   \ \ \ \ \  \langle [f], x \rangle = N \cdot 
f(\bar x), \end{equation} 
with $x=a_0 \ldots a_N$. Here the 
representative $f\in [f]$ is a function that depends on the first 
$n+1$ terms $a_0\ldots a_n$ of sequences in ${\mathcal S}$, and $\bar 
x$ is the truncation of the periodic sequence $\overline{a_0\ldots 
a_N}$ after the first $n$ terms. 
This pairing descends to the direct limits of the filtrations, where 
it agrees with the classical cohomology/homology pairing 
\begin{equation}\label{pairing2} 
\langle \cdot, \cdot \rangle : H^1({\mathcal S}_T,\Z) \times 
H_1({\mathcal S}_T,\Z) \to \Z. 
\end{equation} 
\label{H1pairing} 
\end{prop} 
 
\noindent{\bf Proof.} 
First, it is not hard to check that the pairing \eqref{pairing1} is 
compatible with the maps $F_n \hookrightarrow 
F_{n+1}$ and $J_k: {\mathcal K}_N \rightarrow {\mathcal 
K}_{kN}$. In fact, \eqref{pairing1} is invariant under the maps 
$F_n \hookrightarrow F_{n+1}$, while under the map 
$J_k: {\mathcal K}_N \hookrightarrow {\mathcal K}_{kN}$ we have 
$$ \langle [f], J_k(x) \rangle = kN f(\bar x) = k\,  \langle [f], x 
\rangle. $$ 
Thus, \eqref{pairing1} induces a pairing of the direct 
limits 
\begin{equation}\label{pairing3} 
\langle \cdot, \cdot \rangle : \varinjlim_n F_n \times \varinjlim_N 
{\mathcal K}_N \to \Z. 
\end{equation} 
In order to check that \eqref{pairing3} agrees with the 
cohomology/homology pairing \eqref{pairing2}, notice that 
a class $c$ in the homology $H_1({\mathcal S}_T,\Z)$ is realized as a 
finite linear combination of oriented circles in ${\mathcal 
S}_T$, where each such circle is described by a map $c_x : S^1 \to 
{\mathcal S}_T$ of the form \eqref{c-x}. 
 
The pairing $\langle u, c \rangle$ 
of an element $u$ of the cohomology $H^1({\mathcal S}_T,\Z)$ with a 
generator $c$ of the homology $H_1({\mathcal S}_T,\Z)$ is given by the 
homotopy class $[u\circ c]$ of 
$$ u\circ c : S^1 \longrightarrow S^1. $$ 
 
We write $u(x,t)=[2\pi i t f(x)]$, for a generator of 
$H^1({\mathcal S}_T,\Z)$ in $F_n$, where $f$ is an element in 
${\mathcal P}_n$, which depends only on the first $n+1$ terms in the 
sequences $a_0\ldots a_n\ldots$ in ${\mathcal S}_T$ (\cf Theorem 
\ref{thm-homologydyn}). If $y=a_0 \ldots a_d$ is an element in 
${\mathcal R}_d$ of period $d$, and $c$ 
is the corresponding generator of $H_1({\mathcal S}_T,\Z)$ of the form 
$c_x(t) = [(x,d\cdot t) : x=\overline{a_0\ldots a_d} ]$, 
then the homotopy class $[ u\circ c_x ] \in \pi_1(S^1)=\Z$ is 
equal to $d\,\cdot  f(x)$ and this proves the claim. 
 
\noindent $\diamond$ 
 
\section{Dynamical (co)homology of the fiber at 
infinity}\label{Sect-dyn-cohom} 

In this paragraph we consider the filtered vector space
${\mathcal P}_\kappa = {\mathcal P}\otimes_\Z \kappa$, for $\kappa
=\R$ or $\C$, where ${\mathcal P}$ is the filtered $\Z$ module 
${\mathcal P}\subset {\rm C}({\mathcal S},\Z)$ of functions depending
on future coordinates, as in Theorem \ref{thm-homologydyn}.
Thus, ${\mathcal P}_\kappa$ can be identified with the subspace of
${\rm C}(\Lambda_\Gamma)$ of locally constant $\kappa$-valued
functions. With a slight abuse of notation, we drop the explicit
mention of $\kappa$ and use the same term ${\mathcal P}$ to denote the
vector space, and the notation ${\mathcal P}_n$ for its  
finite dimensional linear subspaces of $\kappa$-valued functions that 
are constant on $\Lambda_\Gamma (\gamma)\subset \Lambda_\Gamma$, for 
all $\gamma\in \Gamma$ of word length $|\gamma|> n+1$. 

We provide a choice of a linear subspace ${\mathcal V}$ 
of the filtered vector space ${\mathcal P}$, which is isomorphic to
the Archimedean  cohomology of Section \ref{3}, compatibly with
the graded structure of the Archimedean  cohomology and the grading
associated to the orthogonal projections of $L^2(\Lambda_\Gamma,\mu)$
onto the subspaces ${\mathcal P}_n$. 

The space ${\mathcal P}$ carries an action of an 
involution, which we denote $\bar F_\infty$ by analogy with the 
real Frobenius acting on the cohomologies of \S \ref{3}.  In 
Theorem \ref{map-ar-dyn} we also show that the embedding 
of $H^1(\X)^{N=0}$ in ${\mathcal P}$ is equivariant with respect 
to the action of the real Frobenius. 

The image of the subspace ${\mathcal V}$ under the quotient map by
the image of the coboundary $\delta = 1-T$ determines a subspace 
$\bar{\mathcal V}$ of the dynamical cohomology isomorphic to the 
Archimedean  cohomology. 
The {\em dynamical cohomology} $H^1_{dyn}$ is defined as the graded 
vector space given by the sum of the 
graded pieces of the filtration of $H^1({\mathcal S}_T)$, introduced 
in Theorem \ref{thm-homologydyn}. These graded pieces ${\rm Gr}_p$ are 
considered with coefficients in the $p$-th Hodge--Tate twist 
$\R(p)$. 
 
Similarly, we define a {\em dynamical homology} $H_1^{dyn}$ as the 
graded vector space given by the 
sum of the terms in the filtration of $H_1({\mathcal S}_T)$, 
introduced in Proposition \ref{homology-dyn}. These vector spaces are 
again considered with twisted $\R(p)$-coefficients. The resulting 
graded vector space also has an action of a real 
Frobenius $\bar F_\infty$. Theorem \ref{map-ar-dyn-2} shows that 
there is an identification of the dual of 
$H^1(\X)^{N=0}$, under the duality isomorphisms in $\H^\cdot (\XX)$ of 
\eqref{Nisos}, with a subspace of $H_1^{dyn}$. The identification is 
compatible with the $\bar F_\infty$ action induced by the change of 
orientation $z\mapsto \bar z$ on $X_{/\R}$. 
 
The pair 
\begin{equation}\label{hom-and-cohom} 
 H^1_{dyn} \oplus H_1^{dyn} 
\end{equation} 
provides a geometric setting, defined in terms of the dynamics of the 
shift operator $T$, which contains a copy of the Archimedean  
cohomology $H^1(\XX)^{N=0}$ and of its dual under the duality 
isomorphisms acting on $\H^\cdot({\rm 
Cone}(N))$. In Theorem \ref{map-ar-dyn-2} we also prove that, 
under these identifications, the duality isomorphism 
corresponds to the homology/cohomology pairing between $H^1({\mathcal 
S}_T)$ and $H_1({\mathcal S}_T)$. 
 
This construction also shows that the map $1-T$ plays, in this 
dynamical setting, a role dual to the monodromy map $N$ of 
the arithmetic construction of Section \ref{3} (\cf Remark 
\ref{Ker-Coker-(1-T)}). 
 
\medskip 

\subsection{Dynamical (co)homology} \label{sectdefcohomdyn}

The terms $F_n$ in the filtration of theorem
\ref{thm-homologydyn} define real (or complex)   
vector spaces, which we still denote $F_n$, in the filtration
of the cohomology $H^1({\mathcal S}_T,\C)$. Since, as $\Z$-modules, 
the $F_n$ are torsion free, the vector spaces obtained by tensoring
with $\R$ or $\C$ are of dimension 
$\dim F_n = 2g(2g-1)^{n-1} (2g-2) +1$ for $n\geq 
1$ and $\dim F_0 =2g$.  
We make the following definition. 
 
\begin{defn} 
Let $H^1({\mathcal S}_T,\R) =\varinjlim_n F_n$, for a filtration $F_n$ as in 
Theorem \ref{thm-homologydyn}, with real coefficients. Let ${\rm 
Gr}_n =F_n/F_{n-1}$ be the corresponding graded pieces, with ${\rm 
Gr}_0 = F_0$. We define the {\em dynamical cohomology} as 
\begin{equation} \label{H-dyn} 
H^1_{dyn} := \oplus_{p\leq 0} gr_{2p}^\Gamma H^1_{dyn}, 
\end{equation} 
where we set 
\begin{equation} \label{H-dyn-p} 
gr_{2p}^\Gamma H^1_{dyn} := {\rm Gr}_{-p} \otimes_{\R} \R(p) 
\end{equation} 
with $\R(p)= (2\pi \sqrt{-1})^p \R$. 
\label{def-H-dyn} 
\end{defn} 
 
Let $\iota_{\Xi} : \tilde\Xi \to \tilde\Xi$ be the 
involution on the orientation double cover $\tilde\Xi$ of $\Xi$ given 
by the $\Z/2$-action. This determines an involution 
$\iota_{{\mathcal S}} :{\mathcal S}\to {\mathcal S}$. The induced 
map $\iota_{{\mathcal S}}^*: {\mathcal P} \to {\mathcal P}$ preserves 
the subspaces ${\mathcal P}_n$ and 
commutes with the coboundary $\delta$, hence it descends to an 
induced involution $\iota_{{\mathcal S}}^*: H^1({\mathcal S}_T,\R) 
\to H^1({\mathcal S}_T,\R)$ which preserves the $F_n$ and induces
a map $\bar \iota^*: H^1_{dyn} \to H^1_{dyn}$.  
 
\begin{defn} \label{FinftyH1}
We define the action of the real Frobenius $\bar F_\infty$ on
$H^1_{dyn}$ as the composition of the involution $\bar \iota^*$ 
induced by the $\Z/2$-action on $\tilde\Xi$ and the action by $(-1)^p$ 
on $\R(p)$. 
\end{defn} 

\medskip 

Consider then the $\R$-vector space ${\mathcal K}_N$ generated by all 
reduced sequences $a_0\ldots a_N$ in the $\{ g_i \}_{i=1}^{2g}$ with 
the condition $a_0 \neq a_N^{-1}$ (\cf Proposition 
\ref{homology-dyn}). 
 
\begin{defn} \label{defdynhom}
We define the {\em dynamical homology} $H_1^{dyn}$ as 
\begin{equation}\label{H1dyn-hom} 
H_1^{dyn}:= \oplus_{p\geq 1} gr_{2p}^\Gamma H_1^{dyn}, 
\end{equation} 
where we set 
\begin{equation}\label{H-hom-dyn-gr} 
gr_{2p}^\Gamma H_1^{dyn} := {\mathcal K}_{p-1}\otimes_\R \R(p) 
\end{equation} 
for $\R(p)= (2\pi \sqrt{-1})^p \R$. The action of $\bar F_\infty$ on 
$H_1^{dyn}$ is given by 
\begin{equation} \label{H-hom-Frob} 
\bar F_\infty ( (2\pi \sqrt{-1})^p \, \, a_0\ldots a_{p-1} ) = 
(-1)^p  (2\pi \sqrt{-1})^p \, \, a_{p-1}^{-1} \ldots a_0^{-1}. 
\end{equation}  
\end{defn}

\subsection{Hilbert completions} 

It is convenient to introduce Hilbert space completions of the vector 
spaces ${\mathcal P}$  and
$H_1^{dyn}$. This will allow us to treat the graded structures and
filtrations in terms of orthogonal projections. It will also play an
important role later, when we consider actions of operator algebras.
The Hilbert spaces can be chosen real or complex.
When we want to preserve the information given by the twisted 
coefficients $\R(p)$, we can choose to work with real coefficients. 

\smallskip

By \eqref{PV}, we can
describe the cohomology $H^1({\mathcal S}_T)$ through the exact sequence
\begin{equation}\label{PVC}
0 \to \C \to {\rm C}({\mathcal S},\Z)\otimes \C
\stackrel{\delta=I-T}{\longrightarrow} {\rm  
C}({\mathcal S},\Z)\otimes \C \to H^1({\mathcal S}_T,\C) \to 0, 
\end{equation}
where ${\rm C}({\mathcal S},\Z)\otimes \C$ are the locally constant,
complex valued functions on ${\mathcal S}$. The same holds with $\R$
instead of $\C$ coefficients. By the argument of Theorem
\ref{thm-homologydyn}, we can replace in this sequence the space of 
locally constant functions on ${\mathcal S}$ by the space ${\mathcal
P}$ of locally constant functions of
future coordinates. These can be identified with locally constant
functions on $\Lambda_\Gamma$,
\begin{equation}\label{PVC2}
0 \to \C \to {\mathcal P} \stackrel{\delta}{\longrightarrow} {\mathcal
P} \to H^1({\mathcal S}_T,\C) \to 0.
\end{equation}

We can consider the (real) Hilbert space
$L^2(\Lambda_\Gamma,\mu)$, of functions of $\Lambda_\Gamma$ that are
square integrable with respect to the Patterson--Sullivan measure
$\mu$ (\cf \cite{Sull}) satisfying
\begin{equation}\label{PSmeas}
(\gamma^* d\mu)(x)= |\gamma^\prime (x)|^{\delta_H} \, d\mu(x), \ \ \
\forall \gamma\in \Gamma. 
\end{equation}
The subspace ${\mathcal P}\subset L^2(\Lambda_\Gamma,\mu)$ is norm
dense, hence we will use ${\mathcal L}=L^2(\Lambda_\Gamma,\mu)$ as
Hilbert space completion of ${\mathcal P}$.

\medskip

On the homology $H_1^{dyn}$ a (real) Hilbert space structure is 
obtained in the  
following way. Each summand ${\mathcal K}_N$ can be regarded as a 
finite dimensional linear subspace of the (real) Hilbert space 
$\ell^2(\Gamma)$, by identifying the generators $a_0\ldots a_N$ with 
$a_0\neq a_N^{-1}$ with a subset of the set of all finite reduces 
words, which is a complete basis of $\ell^2(\Gamma)$. This determines 
the inner product on each ${\mathcal K}_N$. We denote the 
corresponding norm by $\| \cdot \|_{{\mathcal K}_N}$. 
We obtain a real Hilbert space structure on $H_1^{dyn}$ with norm 
\begin{equation}\label{norm-hom} 
 \| \sum_p x_p \| := \left( \sum_p \| x_p \|_{{\mathcal K}_N}^2 
\right)^{1/2}. 
\end{equation} 
We denote the Hilbert space completion of $H_1^{dyn}$ in this norm by 
${\mathcal H}_1^{dyn}$.

\subsection{Arithmetic cohomology and dynamics}\label{arithmdynsect}

We construct a linear subspace ${\mathcal V}$ of ${\mathcal P}$ 
isomorphic to the Archimedean  cohomology of Section \ref{3},
compatibly with the action of the real Frobenius $\bar F_\infty$. 

Let $\Pi_n$ denote the orthogonal projections onto the subspaces
${\mathcal P}_n$, with respect to the inner product on the Hilbert
space ${\mathcal L}=L^2(\Lambda_\Gamma, \mu)$. We denote by
$\hat\Pi_n$ the projections $\hat\Pi_n = \Pi_n - \Pi_{n-1}$ onto the
subspaces ${\mathcal P}_n \cap {\mathcal P}_{n-1}^\perp$. These
determine an associated grading operator, the unbounded self--adjoint
operator $D: {\mathcal L}\to {\mathcal L}$ 
\begin{equation}\label{D-L}
D = \sum_n n\, \hat\Pi_n.
\end{equation} 

Consider again the involution $\iota_{\Xi} : \tilde\Xi \to \tilde\Xi$ 
and the induced $\iota_{{\mathcal S}} :{\mathcal S}\to {\mathcal S}$. 

\begin{defn}\label{FinftydefP}
We define the action of the real Frobenius $\bar F_\infty$ on the 
space ${\mathcal P}$ by
\begin{equation}\label{FinftyP}
 \bar F_\infty = (-1)^D \iota_{{\mathcal S}}.
\end{equation}
This extends to a bounded operator on ${\mathcal L}$.
\end{defn}

Notice that $\iota_{{\mathcal S}} \circ \hat\Pi_n = \hat\Pi_n \circ
\iota_{{\mathcal S}}$, for all $n\geq 1$, hence the involution $\bar
F_\infty$ commutes with all the projections $\hat\Pi_n$. 

Similarly, we have an analog of the Tate twist by $\R(p)$ on the
coefficients of the Archimedean  cohomology of Section \ref{3}, given
by the action on ${\mathcal P}$ of the linear operator $(2\pi
\sqrt{-1})^D$. 

\medskip

Let $\chi_{{\mathcal S}^+(w_{n,k})}$ denote the characteristic functions 
of the sets ${\mathcal S}^+(w_{n,k}) \subset {\mathcal S}^+$, where 
$w_{n,k}$ is a word in the $\{ g_j \}_{j=1}^{2g}$ of the form 
$$ w_{n,k} = \underbrace{g_k g_k\cdots g_k}_{n-times}. $$ 
 
\begin{lem} 
The functions $\chi_{{\mathcal S}^+(w_{n,k})}\in {\mathcal P}_n$ have 
the following properties: 
\begin{enumerate}
\item The elements $\hat\Pi_n \chi_{{\mathcal S}^+(w_{n,k})}$ are
linearly independent in ${\mathcal P}_n \cap {\mathcal
P}_{n-1}^\perp$. 
\item The images under the quotient map
\begin{equation}\label{classes-Gn} 
\chi_{n,k}:=[ \chi_{{\mathcal S}^+(w_{n,k})} ] \in {\rm Gr}_{n-1} 
\end{equation} 
are all linearly independent, hence they space a $2g$ 
dimensional subspace in each ${\rm Gr}_{n-1}\subset H^1_{dyn}$.
\end{enumerate}
\end{lem} 
 
\noindent{\bf Proof.} {\em 1.} The characteristic 
functions $\chi_{{\mathcal S}^+(w_{n,k})}$ for $n\geq 1$ and 
$k=1,\ldots, 2g$ are all linearly independent in ${\mathcal P}_n$. 
Moreover, no linear combination of the $\chi_{{\mathcal
S}^+(w_{n,k})}$ lies in ${\mathcal P}_{n-1}$.  
{\em 2.} The pairing of Proposition \ref{H1pairing} 
with $T$-invariant elements $g_ig_i\ldots g_i\ldots$ in the
dynamical homology shows that no linear combination of the
$\chi_{{\mathcal S}^+(w_{n,k})}$ lies in the 
image of $\delta=1-T$. Thus, passing to  
equivalence classes modulo the image of the map $\delta = 1-T$, we obtain 
linearly independent elements 
\begin{equation}\label{classes-Fn} 
\left( \chi_{{\mathcal S}^+(w_{n,k})} \mod (1-T) \right) \in F_{n-1} 
\subset H^1({\mathcal S}_T). 
\end{equation} 
Again, as in {\em 1.}, for any fixed $n$, no linear 
combination of the classes \eqref{classes-Fn} lies in $F_{n-2}$. 
Thus, by further taking the equivalence classes of the 
\eqref{classes-Fn} modulo $F_{n-2}$, we obtain $2g$ linearly 
independent elements \eqref{classes-Gn} in each ${\rm Gr}_{n-1}$. 
 
\noindent $\diamond$ 
 
\medskip 
 
We obtain elements in $H^1_{dyn}$ by considering 
\begin{equation}\label{classes-H1dyn} 
 (2\pi \sqrt{-1})^p \chi_{-p+1,k} \in gr_{2p}^\Gamma H^1_{dyn}. 
\end{equation} 
 
\begin{defn}\label{Vdef} 
We denote by ${\mathcal V}\subset {\mathcal P}$ the linear vector
space spanned by the elements $\hat\Pi_n \chi_{{\mathcal
S}^+(w_{n,k})}$. This is a graded vector space 
\begin{equation}\label{VdefP}
 {\mathcal V}:= \oplus_{p\leq 0} gr_{2p}^\Gamma {\mathcal V}, 
\end{equation}
with $gr_{2p}^\Gamma {\mathcal V} = \hat\Pi_{|p|} {\mathcal V}$. 
Similarly, we denote by $\bar{\mathcal V}\subset H^1_{dyn}$ the graded
subspace  
$$\bar{\mathcal V}:= \oplus_{p\leq 0} gr_{2p}^\Gamma \bar{\mathcal V}$$ 
where $gr_{2p}^\Gamma \bar{\mathcal V}$ is the subspace of $gr_{2p}^\Gamma 
H^1_{dyn}$ spanned by the elements of the form \eqref{classes-H1dyn}, 
for $k=1,\ldots, 2g$. 
\end{defn} 

Notice that in \eqref{VdefP} we have added a sign to the grading
($p\leq 0$), in order to match the sign of the grading of the
Archimedean  cohomology. This means that we have to introduce a sign in
the grading operator \eqref{D-L}. We discuss this more precisely when
we introduce the dynamical spectral triple. 
 
\medskip 
 
Now we show that there is a natural definition of a map from the
Archimedean  cohomology to the space ${\mathcal V}$ and to the
dynamical cohomology. This involves a basis of holomorphic
differentials determined by the Schottky uniformization. 

\smallskip
 
For $k=1,\ldots, g$, let $\eta_k$ denote a basis of holomorphic 
differentials on the Riemann surface $X_{/\R}$ satisfying the 
normalization 
\begin{equation}\label{normaliz} 
\int_{a_j} \eta_k =\delta_{jk}, 
\end{equation} where the $a_k$ 
are a basis of $Ker(I_*)$ for $I_* : H_1(X_{/\R},\Z) \to 
H_1(\mX_\Gamma,\Z)$ induced by the inclusion of $X_{/\R}$ as the 
boundary at infinity of the handlebody $\mX_\Gamma$. We refer to this 
basis as the {\em canonical basis} of holomorphic differentials. 
 
Recall that, since we are considering an orthosymmetric smooth real 
algebraic curve, $X_{/\R}$ has a Schottky uniformization by a Fuchsian 
Schottky group as in Proposition \ref{real-Bers}. 
It is known then (\cf \cite{Man} \cite{Sepp}) that a 
holomorphic differential $\eta$ on $X_{/\R}$ can be obtained as 
Poincar\'e series with exponent $1$, $\eta = \Theta^1(f)$, where $f$ 
is a meromorphic function on $\P^1(\C)$ with divisor $D(f)\subset 
\Omega_\Gamma$, and 
\begin{equation}\label{Poin-ser} 
\Theta^m (f)(z) :=\sum_{\gamma \in \Gamma} f(\gamma(z)) 
\left(\frac{\partial \gamma(z)}{\partial z}\right)^m 
\end{equation} 
is the Poincar\'e series with exponent $m$. The fact that the 
Hausdorff dimension of the limit set satisfies 
$\dim_H(\Lambda_\Gamma)<1$ ensures absolute convergence on 
compact sets in $\Omega_\Gamma$ (\cf Remark \ref{delta<1}). 
In particular, consider the automorphic series 
\begin{equation}\label{omega-k} 
\omega_k =\sum_{h \in C(\cdot | g_k)} d_z \log \langle h z^+(g_k), h 
z^-(g_k), z, z_0 \rangle, 
\end{equation} 
where $\langle a,b,c,d \rangle$ is the cross ratio of points in 
$\P^1(\C)$, and $C(\cdot | g_k)$ denotes a set of representatives of 
the coset classes $\Gamma / \Z\langle g_k \rangle$, for $\{ g_k 
\}_{k=1}^g $ the generators of $\Gamma$, and $z_0$ a base point in 
$\Omega_\Gamma$. By Lemma 8.2 of \cite{Man} (\cf Proposition 1.5.2 of 
\cite{ManMar2}), the expression \eqref{omega-k} gives the canonical 
basis of holomorphic differentials satisfying the normalization 
condition \eqref{normaliz}, in the form 
\begin{equation}\label{omega-eta} 
 \eta_k = \frac{1}{2\pi \sqrt{-1}} \omega_k,  \ \ \ \text{ so that } 
\ \ \int_{a_j} \omega_k =(2\pi \sqrt{-1}) \delta_{jk}. 
\end{equation} 
 
The formula \eqref{omega-k} gives the explicit correspondence between 
the set of generators of $\Gamma$ and the canonical basis of holomorphic 
differentials $g_k \mapsto \omega_k$, for $k=1,\ldots, g$, which we 
use in order to produce the following identification. 
 
\begin{thm}\label{map-ar-dyn} 
Consider the map
\begin{equation}\label{map-U} 
U: gr_{2p}^\w H^1(\X)^{N=0} \longrightarrow 
gr_{2p}^\Gamma {\mathcal V}, 
\end{equation} 
given by
\begin{equation}\label{map-U2} 
U((2\pi \sqrt{-1})^{p-1} \varphi_k):= (2\pi\sqrt{-1})^p \,\,
\hat\Pi_{|p|} 
\frac{\chi_{{\mathcal S}^+(w_{n,k})} - \chi_{{\mathcal
S}^+(w_{n,k+g})}}{2}  
\end{equation} 
\begin{equation}\label{map-U3} 
U((2\pi \sqrt{-1})^{p-1} \varphi_{k+g}):= (2\pi\sqrt{-1})^p \,\,
\hat\Pi_{|p|} 
\frac{\chi_{{\mathcal S}^+(w_{n,k})} + \chi_{{\mathcal
S}^+(w_{n,k+g})}}{2},
\end{equation} 
for $k=1,\ldots,g$ and $p\leq 0$, where we set 
\begin{equation}\label{varphi} 
\varphi_k=(\omega_k + \bar\omega_k)/2 \ \ \ \text{ and }  \ \ \ 
\varphi_{g+k}=-i(\omega_k -\bar\omega_k)/2. 
\end{equation}
Consider also the map
\begin{equation}\label{map-Ubar} 
\bar U: gr_{2p}^\w H^1(\X)^{N=0} \longrightarrow 
gr_{2p}^\Gamma \bar{\mathcal V}, 
\end{equation} 
given by 
\begin{equation}\label{map-U2bar} 
\bar  U((2\pi \sqrt{-1})^{p-1} \varphi_k):= (2\pi\sqrt{-1})^p \,\, 
\frac{\chi_{-p+1,k} - \chi_{-p+1,k+g}}{2} 
\end{equation} 
\begin{equation}\label{map-U3bar} 
\bar  U((2\pi \sqrt{-1})^{p-1} \varphi_{k+g}):= (2\pi\sqrt{-1})^p \,\, 
\frac{\chi_{-p+1,k} + \chi_{-p+1,k+g}}{2}. 
\end{equation} 
The map $U$ is an isomorphism of 
$H^1(\X)^{N=0}$ and ${\mathcal V}\subset {\mathcal P}$.
It is equivariant with respect 
to the action of the real Frobenius $\bar F_\infty$. The Tate twist
that gives the grading of $H^1(\X)^{N=0}$ corresponds to the action of
$(2\pi \sqrt{-1})^{-D}$, for $D$ in \eqref{D-L}, on ${\mathcal V}$. 
The map $\bar U$ is an isomorphism of 
$H^1(\X)^{N=0}$ and $\bar{\mathcal V}$ as graded vector spaces, which
is equivariant with respect to the action of the real Frobenius $\bar
F_\infty$. 
\end{thm} 
 
\noindent{\bf Proof.} 
The elements \eqref{varphi} give a basis of 
$H^1_{DR}(X_{/\R},\R(1))$, where the twist is due to the choice of 
normalization \eqref{normaliz} and the relation \eqref{omega-eta}. A 
basis for $H^1(X_{/\R},\R(p))$ is 
then given by the $(2\pi \sqrt{-1})^{p-1} \varphi_k$, $k=1,\ldots, 
2g$. 
 
By \eqref{HY}, we have $H^1(\X)^{N=0} = \oplus_{p\leq 0} gr_{2p}^\w 
H^1(\X)^{N=0}$, with $gr_{2p}^\w H^1(\X)^{N=0} 
=H^1(X_{/\R},\R(p))$. Thus, we obtain a basis for $H^1(\X)^{N=0}$, of 
the form 
$$ \{ (2\pi \sqrt{-1})^{p-1} \varphi_k :\,  k=1,\ldots, 2g, \, p\leq 0 
\}. $$ 
 
By construction, the maps \eqref{map-U} and \eqref{map-Ubar} 
define isomorphisms of graded vector spaces. 
We need to check that they are equivariant with respect to the 
action of the real Frobenius. 
 
In the case of a real $X_{/\R}$, the action of complex 
conjugation $z\mapsto \bar z$ corresponds geometrically to a change of 
orientation on $X_{/\R}$, which induces a change 
of orientation on the handlebody $\mX_\Gamma$. If $L_{\{a,b\}}$ is 
the geodesic in $\tilde\Xi$ such that the orientation of the geodesic 
at the endpoint $b\in \P^1(\C)$ agrees with the outward pointing 
normal vector, then under the change of orientation induced by 
$z\mapsto \bar z$ the geodesic $L_{\{a,b\}}$ is exchanged with 
$L_{\{b,a\}}$, which is exactly the effect of the involution on 
$\tilde\Xi$. The induced involution on ${\mathcal P}$ exchanges 
$\chi_{{\mathcal S}^+(w_{n,k})}$ and $\chi_{{\mathcal 
S}^+(w_{n,g+k})}$. Since we have $\bar F_\infty \hat\Pi_n =\hat\Pi_n
\bar F_\infty$, we obtain
$$ \bar F_\infty ((2\pi \sqrt{-1})^p \hat\Pi_{|p|} \chi_{{\mathcal
S}^+(w_{n,k})} = (-1)^p (2\pi \sqrt{-1})^p  \hat\Pi_{|p|} \chi_{{\mathcal
S}^+(w_{n,k+g})}. $$ 
Similarly, the involution on $\bar{\mathcal 
V}$ is given by 
$$ \bar F_\infty ((2\pi \sqrt{-1})^p \chi_{-p+1,k}) = (-1)^p (2\pi 
\sqrt{-1})^p  \chi_{-p+1,g+k}. $$ 
 
On the other hand, under the action of the real Frobenius $\bar 
F_\infty$ we have $H^1(X_{/\R},\R)= E_1 \oplus E_{-1}$ with $\dim 
E_{\pm 1}=g$ (\cf Remark \ref{rem4}), generated respectively by the 
$(2\pi \sqrt{-1})^{-1}\varphi_k$ and $(2\pi 
\sqrt{-1})^{-1}\varphi_{g+k}$, for $k=1,\ldots, g$. This gives the 
corresponding splitting into eigenspaces $H^1(\X)^{N=0}=E^+ \oplus 
E^-$ as in \eqref{deco}. Thus, we see that \eqref{map-U2} and 
\eqref{map-U3} and \eqref{map-U2bar} and 
\eqref{map-U3bar} are $\bar F_\infty$-equivariant, since we have 
\begin{align*} 
\bar F_\infty \left( (2\pi \sqrt{-1})^p \, \frac{1}{2} (
\chi_{{\mathcal S}^+(w_{n,k})} \pm \chi_{{\mathcal S}^+(w_{n,k+g})} )
\right) & = \pm 1 \cdot (-1)^p \cdot  \\ 
& (2\pi \sqrt{-1})^p \, \frac{1}{2} ( 
\chi_{{\mathcal S}^+(w_{n,k})} \pm \chi_{{\mathcal S}^+(w_{n,k+g})}
 ) \\[3mm] 
\bar F_\infty \left( (2\pi \sqrt{-1})^p \, \frac{1}{2} ( \chi_{-p+1,k} 
\pm \chi_{-p+1,k+g} ) \right) & = \pm 1 \cdot (-1)^p \cdot  \\ 
& (2\pi \sqrt{-1})^p \, \frac{1}{2} ( \chi_{-p+1,k} \pm 
\chi_{-p+1,k+g} ) \\[3mm] 
\bar F_\infty\left( (2\pi\sqrt{-1})^{p-1} \varphi_k \right) 
& = +1 \cdot (-1)^{p-1}  (2\pi\sqrt{-1})^{p-1} \varphi_k \\[3mm] 
\bar F_\infty\left( (2\pi \sqrt{-1})^{p-1} \varphi_{k+g} \right) 
& = -1 \cdot (-1)^{p-1} (2\pi\sqrt{-1})^{p-1} \varphi_{k+g}. 
\end{align*} 
 
\noindent $\diamond$ 
 
\medskip 
 
\begin{rem} {\em The characteristic function $\chi_{{\mathcal 
S}^+(w_{n,k})}$ can be regarded as the ``best 
approximation'' within ${\mathcal P}_n$ (\cf Theorem 
\ref{thm-homologydyn}.{\em 2.}) to a function (non-continuous) 
supported on the periodic sequence of period $g_k$, 
$$ f(g_k g_k g_k g_k g_k \ldots) =1 \ \ \text{ and } \ \ f(a_0 a_1 a_2 a_3 
\ldots)=0 \, \, \text{ otherwise.} $$ 
In fact, the grading by $n$ should be regarded as the choice of a {\em
cutoff} on the cohomology ${\mathcal H}^1$, corresponding to the
choice of a ``mesh'' on $\Lambda_\Gamma$.
The periodic sequence $g_k g_k g_k g_k g_k \ldots$ 
represents under the correspondence of Lemma \ref{iso-coding1} the 
closed geodesic in $\tilde\Xi$ that is an oriented core handle of the 
handlebody $\mX_\Gamma$. Thus, the elements \eqref{classes-H1dyn} that 
span the subspace $gr_{2p}^\Gamma \bar{\mathcal V} \subset 
gr_{2p}^\Gamma H^1_{dyn}$ can be regarded as the ``best 
approximations'' within $gr_{2p}^\Gamma H^1_{dyn}$ to cohomology 
classes supported on the core handles of the handlebody. 
In other words, we may regard the index $p\leq 0$ in the graded 
structure ${\mathcal V}=\oplus_p gr_{2p}^\Gamma {\mathcal V}$ as 
measuring a way of ``zooming in'', with increasing precision for 
larger $|p|$, on the core handles of the handlebody $\mX_\Gamma$. } 
\end{rem}

Notice that, while the Archimedean  cohomology $H^1(\X)^{N=0}$ 
is identified with the {\em kernel} of the monodromy 
map, the dynamical cohomology is constructed by 
considering the {\em cokernel} of the map $\delta =1-T$. This suggests a 
duality between the monodromy $N$ and the map $1-T$. 
This will be made more precise in the next paragraph. 
 
\subsection{Duality isomorphisms}\label{arithmdynsect2}  

We identify a copy of the dual of the Archimedean  cohomology inside
the dynamical homology $H_1^{dyn}$.
 
\begin{defn} 
We define the linear subspace ${\mathcal W}\subset H_1^{dyn}$ to be
the graded vector ${\mathcal W}=\oplus_{p\geq 1} gr_{2p}^\Gamma
{\mathcal W}$, where $gr_{2p}^\Gamma {\mathcal W}$ is generated by the
$2g$ elements  
$$ (2\pi \sqrt{-1})^p \, \, \underbrace{g_k g_k \ldots 
g_k}_{p-times} . $$ 
\end{defn} 
 
\begin{rem} {\em Notice that the generators of $H_1^{dyn}$ are 
periodic sequences $\overline{a_0\ldots a_N}$, hence elements in 
$Ker(1-T^d)$ for $d$ the period length, $d|N$. Notice in particular 
that the subspace ${\mathcal W}$ can be identified with the part of 
$H_1^{dyn}$ that is generated by elements in $Ker(1-T)$, \ie periodic 
sequences with period length $d=1$. } \label{Ker(1-T)} 
\end{rem} 
 
The subspace ${\mathcal W}$ of the dynamical 
homology is related both to the subspace $\bar{\mathcal V}$ of the 
dynamical cohomology, and to the space 
$$ \oplus_{r\geq 2} gr_{2r}^\w H^3_Y(X) \cong \delta_1( \oplus_{p\leq 
0} gr_{2p}^\w H^1(\X)^{N=0} ), $$ 
for $r=-p+2$, and with $\delta_1$ the duality isomorphism of 
\eqref{isoms-cone} (\cf Proposition \ref{N-isos}). 
 
\begin{lem}\label{map-VW-lem} 
The homology/cohomology pairing \eqref{pairing1} induces an 
identification 
\begin{equation}\label{tildedelta1}
 \tilde \delta_1: gr_{2p}^\Gamma \bar{\mathcal V} 
\stackrel{\simeq}{\longrightarrow} gr_{2(-p+1)}^\Gamma {\mathcal W} 
\end{equation} 
of the vector spaces $gr_{2p}^\Gamma {\mathcal W}$ and 
$gr_{2p}^\Gamma \bar{\mathcal V}$, for all $p\leq 0$.  The map 
$\tilde \delta_1$ satisfies $\tilde \delta_1\circ \bar F_\infty + \bar 
F_\infty \circ \tilde \delta_1 =0$. 
\end{lem} 
 
\noindent{\bf Proof.} The pairing \eqref{pairing1} of the class of the 
characteristic function $\chi_{{\mathcal S}^+(w_{-p+1,k})}$ with the 
element $\underbrace{g_j g_j \ldots 
g_j}_{(-p+1)-times}$ is 
$$ \langle \chi_{{\mathcal S}^+(w_{-p+1,k})}, \underbrace{g_j 
g_j \ldots g_j}_{(-p+1)-times} \rangle = (-p+1) \, 
\delta_{jk}, $$ 
for $j,k  =  1\ldots 2g$. This induces a pairing 
\begin{equation}\label{pairR(1)} 
\langle \cdot, \cdot \rangle : gr_{2p}^\Gamma \bar{\mathcal V} \times 
gr_{2(-p+1)}^\Gamma {\mathcal W} \to \R(1) 
\end{equation} 
$$ \langle (2\pi \sqrt{-1})^p \chi_{-p+1,k}, (2\pi \sqrt{-1})^{(-p+1)} 
\, \underbrace{g_j g_j \ldots g_j}_{(-p+1)-times} \rangle = 
(2\pi\sqrt{-1}) (-p+1) \delta_{jk}. $$ 
Via this pairing, we obtain an identification 
$$ \tilde \delta_1: gr_{2p}^\Gamma \bar{\mathcal V} 
\stackrel{\simeq}{\longrightarrow} gr_{2(-p+1)}^\Gamma {\mathcal W} $$ 
of the form 
\begin{equation}\label{VW-map} 
\begin{array}{ll} 
\tilde \delta_1 : & (2\pi \sqrt{-1})^p \frac{1}{2}(\chi_{ -p+1,k }\pm 
\chi_{-p+1,k+g})  \\[3mm] \mapsto & 
\frac{(2\pi \sqrt{-1})^{(-p+1)}}{(-p+1)} \frac{1}{2}\left( 
\underbrace{g_k g_k \ldots g_k}_{(-p+1)-times} \pm 
\underbrace{g_k^{-1} g_k^{-1} \ldots g_k^{-1}}_{(-p+1)-times} 
\right). 
\end{array} 
\end{equation} 
The relation $\tilde \delta_1\circ \bar F_\infty + \bar 
F_\infty \circ \tilde \delta_1 =0$ follows by construction. 
 
\noindent $\diamond$ 
 
We then obtain the following result. 
 
\begin{thm}\label{map-ar-dyn-2} 
For $p\leq 0$, consider the map 
\begin{equation}\label{map-tildeU} 
\tilde U: gr_{2(-p+2)}^\w H^2(\XX) \longrightarrow 
gr_{2(-p+1)}^\Gamma {\mathcal W}, 
\end{equation} 
given by 
\begin{equation}\label{map-tildeU2} 
 \tilde U ((2\pi \sqrt{-1})^{-p} \varphi_k):= (2\pi\sqrt{-1})^{-p+1} \,\, 
\frac{1}{2(-p+1)}\left( \underbrace{g_k \ldots g_k}_{(-p+1)-times} - 
\underbrace{g_k^{-1} \ldots g_k^{-1}}_{(-p+1)-times} \right) 
\end{equation} 
\begin{equation}\label{map-tildeU3} 
 \tilde U((2\pi \sqrt{-1})^{-p} \varphi_{k+g}):= (2\pi\sqrt{-1})^{-p+1} \,\, 
\frac{1}{2(-p+1)} \left( \underbrace{g_k \ldots g_k}_{(-p+1)-times} + 
\underbrace{g_k^{-1} \ldots g_k^{-1}}_{(-p+1)-times} \right), 
\end{equation} 
for $k=1,\ldots,g$, with $\varphi$ as in \eqref{varphi}. 
The map $\tilde U$ is a $\bar 
F_\infty$-equivariant isomorphism of the graded vector 
spaces $\oplus_{r=-p+2\geq 
2}  gr_{2r}^\w H^2(\XX)$ and ${\mathcal W}$. Moreover, the following 
diagram commutes and is compatible with the action of $\bar F_\infty$: 
\begin{eqnarray} 
\diagram 
gr_{2p}^\w H^1(\X)^{N=0} \dto^{\bar U}\rto^{\delta_1} &  gr_{2(-p+2)}^\w 
H^2(\XX)\dto^{\tilde U} \\ 
gr_{2p}^\Gamma \bar{\mathcal V}\rto^{\tilde \delta_1} & gr_{2(-p+1)}^\Gamma 
{\mathcal W} 
\enddiagram 
\label{diagr1} 
\end{eqnarray} 
\end{thm} 
 
\noindent{\bf Proof.} We have $gr_{2(-p+2)}^\w H^2(\XX) \cong 
H^1(X_{/\R},\R(-p+1))$ by Propositions \ref{descr} and \ref{grad}. 
The duality isomorphism \eqref{isoms-cone} (\cf Proposition 
\ref{N-isos}) 
$$ \delta_1 : gr_{2p}^\w H^1(\X)^{N=0} \stackrel{\simeq}{\to} 
gr_{2(-p+2)}^\w H^2(\XX), $$ 
is given by 
$$ N^{2p-1}: H^1(X_{/\R},\R(p)) \stackrel{\simeq}{\to} 
H^1(X_{/\R},\R(-p+1)), $$ 
$$ (2\pi \sqrt{-1})^{p-1} \varphi_k \mapsto (2\pi \sqrt{-1})^{-p} 
\varphi_k. $$ 
Notice that the duality isomorphism $\delta_1$ also satisfies 
$\delta_1 \circ \bar F_\infty + \bar F_\infty \circ \delta_1=0$. The result 
then follows immediately.

\noindent $\diamond$ 
 
\begin{rem} {\em There is an intrinsic duality in the identifications 
of diagram \eqref{diagr1}. In fact, as discussed in the Remark 
\ref{Ker(1-T)}, the space $H^1(\X)^{N=0}$ 
corresponds to ${\rm Ker}(N)$ is identified with ${\mathcal V}$, which 
is obtained by considering the ${\rm Coker}(1-T)$, while the image of 
$H^1(\X)^{N=0}$ under the  duality isomorphism $\delta_1$ is obtained 
by taking the ${\rm Coker}(N)$ in $\H^\cdot(\XX)$ and is identified 
with ${\mathcal W}$, which can be identified with ${\rm Ker}(1-T)$. Modulo 
this duality, we have a correspondence between the monodromy map $N$ 
and the dynamical map $1-T$. The presence of this duality is not 
surprising, considering that the cohomological construction of Section 
\ref{3} is a theory of the special fiber, while the dynamical 
construction of Section \ref{6} and \cite{Man} is a theory of the {\em 
dual} graph.} \label{Ker-Coker-(1-T)} 
\end{rem}

\section{Dynamical spectral triple.}\label{6bis} 

In Proposition \ref{bound-commOA} and
Theorem \ref{dyn-SP3OA} we construct a spectral triple $({\mathcal A},
{\mathcal H}, D)$ associated to the dual graph ${\mathcal S}_T$, where
the Hilbert space is given by the cochains of the dynamical cohomology
and the algebra is the crossed product ${\rm C}(\Lambda_\Gamma)\rtimes
\Gamma$, describing the action of the Schottky group on its
limit set. 

We recall the construction and basic properties of the Cuntz-Krieger 
algebra ${\mathcal O}_A$ associated to the shift of finite type 
$({\mathcal S},T)$. In Theorem \ref{OAprod} we show that this algebra 
describes, as a non-commutative space, the quotient of the limit set 
$\Lambda_\Gamma$ by the action of the Schottky group 
$\Gamma$. 

The Dirac operator $D$, defined by the grading operator \eqref{D-L} and
a sign, restricts to the subspace ${\mathcal V}$ 
isomorphic to the Archimedean  cohomology $H^1(\X)^{N=0}$ to the 
Frobenius-type operator $\Phi$ of Section \ref{3}. This ensures that
we can recover the local factor at arithmetic infinity from the
spectral geometry.
We refer to these data as {\em dynamical spectral 
triple}. 

We then show, in \S \ref{modinftysect}, that the homotopy quotient
$\Lambda_\Gamma \times_\Gamma \H^3$ provides an analog, in the 
$\infty$-adic case, of the $p$-adic reduction map obtained 
by considering the reductions mod $p^k$ (\cf \cite{Man} \cite{Mum}).

\subsection{Cuntz--Krieger algebra}\label{CKsect}

A partial isometry is a linear operator $S$ satisfying the relation 
$S= S S^* S$. The Cuntz--Krieger algebra ${\mathcal O}_A$ (\cf 
\cite{Cu} \cite{CuKrie}) is defined 
as the universal ${\rm C}^*$--algebra generated by partial isometries 
$S_1, \ldots, S_{2g}$, satisfying the relations 
\begin{equation} \label{CK1rel} \sum_j S_j S_j^* =I \end{equation} 
\begin{equation} \label{CK2rel} S_i^* S_i =\sum_j A_{ij} \, S_j 
S_j^*, \end{equation} 
where $A=(A_{ij})$ is the $2g\times 2g$ transition matrix of the 
subshift of finite type $({\mathcal S},T)$, namely the matrix whose 
entries are $A_{ij}=1$ whenever $|i-j|\neq g$, and $A_{ij}=0$ 
otherwise. 
 
We give a more explicit description of the generators of the 
Cuntz--Krieger algebra ${\mathcal O}_A$. 
 
Consider the following operators, acting on the Hilbert space 
$L^2(\Lambda_\Gamma, \mu)$, where $\mu$ is the Patterson--Sullivan 
measure on the limit set (\cf \cite{Sull}): 
\begin{equation} (T_{\gamma^{-1}} f)(x)
:=|\gamma^\prime(x)|^{\delta_H/2}\, f(\gamma x),  \ \ \ \text{ and } \ \ \
(P_\gamma f)(x) := \chi_\gamma (x) f(x), \label{TandP}  
\end{equation} 
where $\delta_H$ is the Hausdorff dimension of $\Lambda_\Gamma$ and the
element $\gamma\in \Gamma$ is identified with a reduced word in the 
generators $\{ g_j \}_{j=1}^g$ and their inverses, and 
$\chi_\gamma$ is the characteristic function of the cylinder 
$\Lambda_\Gamma(\gamma)$ of all (right) infinite reduced words 
that begin with the word $\gamma$. Then, for all $\gamma\in \Gamma$, 
$T_\gamma$ is a unitary operator with $T^*_\gamma = T_{\gamma^{-1}}$ 
and $P_{\gamma}$ is a projector. 
In particular, with the usual notation 
$g_{j+g} = g_j^{-1}$, for $j=1,\ldots,g$, we write 
$$ T_j := T_{g_j} \ \ \ \text{ and } \ \ \ 
P_j := P_{g_j} \ \ \ \text{ for } j=1,\ldots, 2g. $$ 
 
\begin{prop} \label{actionSi} 
The operators 
\begin{equation} \label{Si} S_i := \sum_{j} A_{ij} T_i^* P_j 
\end{equation} 
are bounded operators on $L^2(\Lambda_\Gamma, \mu)$ satisfying the 
relations \eqref{CK1rel} and \eqref{CK2rel}. Thus, the Cuntz-Krieger 
algebra ${\mathcal O}_A$ can be identified with the 
subalgebra of bounded operators on the Hilbert space 
$L^2(\Lambda_\Gamma, \mu)$ generated by the $S_i$ as in \eqref{Si}. 
\end{prop} 
 
\noindent{\bf Proof.} 
The operators $P_i$ are orthogonal projectors, \ie $P_i P_j 
=\delta_{ij} P_j$.  The composite $T_i^* P_j T_i$ satisfies 
$$ \sum_j A_{ij} (T_i^* P_j T_i \, f)(x) = \left\{ \begin{array}{lr} 
f(x) & \text{ if } P_i(x)=x \\[2mm] 
0 & \text{ otherwise. } \end{array} \right. $$ 
In fact, we have 
$$ \sum_j A_{ij} (T_i^* P_j T_i f)(x)= \sum_j A_{ij} 
T_i^* P_j \,|g_i^\prime (g_i^{-1} x)|^{-\delta_H/2}\, f (g_i^{-1} 
x) $$
$$ = \left\{\begin{array}{ll} \sum_j A_{ij}T_i^* P_j \, |g_i^\prime
(T x)|^{-\delta_H/2}\, f (T x) \,\, = f(x) & a_0 = g_i \\[2mm]
0 &  a_0\neq g_i. \end{array} \right. $$
This implies that the $S_i$ and $S_i^*$ satisfy 
$$ S_i S_i^* = \sum_j A_{ij} T_i^* P_j T_i = P_i. $$ 
Since the projectors $P_i$ satisfy 
$\sum_i P_i = I$, we obtain the relation \eqref{CK1rel}. 
Moreover, since $T_i T_i^* = 1$, and the entries $A_{ij}$ are all 
zeroes and ones, we also obtain 
$$ S_i^* S_i = \sum_{j,k} A_{ij}A_{ik} P_k T_i T_i^* P_j = \sum_j 
(A_{ij})^2 \,P_j = \sum_j A_{ij} P_j. $$ 
Replacing $P_j = S_j S_j^*$ from \eqref{CK1rel} we then obtain 
\eqref{CK2rel}. 
 
\noindent $\diamond$ 
 
The Cuntz--Krieger algebra ${\mathcal O}_A$ can be described in terms 
of the action of the free group $\Gamma$ on its limit set 
$\Lambda_\Gamma$ (\cf \cite{Rob}, \cite{Spi}), so that we can regard 
${\mathcal O}_A$ as a noncommutative space replacing the classical 
quotient $\Lambda_\Gamma / \Gamma$. 
 
In fact (\cf Lemma \ref{S-Lamb}), the action of $\Gamma$ on 
$\Lambda_\Gamma \subset \P^1(\C)$ determines a unitary representation 
\begin{equation}\label{AutGamma} \Gamma \to {\rm Aut}({\rm 
C}(\Lambda_\Gamma)) \ \ \ \  (T_{\gamma^{-1}} f)(x) =
|\gamma^\prime(x)|^{\delta_H/2}\, f(\gamma x),  
\end{equation} 
where ${\rm C}(\Lambda_\Gamma)$ is the ${\rm C}^*$--algebra of 
continuous functions on $\Lambda_\Gamma$. Thus, 
we can form the (reduced) crossed product ${\rm C}^*$--algebra 
${\rm C}(\Lambda_\Gamma)\rtimes \Gamma$. 
 
In the following Theorem we construct an explicit identification between 
the algebra ${\rm C}(\Lambda_\Gamma)\rtimes \Gamma$ and the subalgebra 
of bounded operators on $L^2(\Lambda_\Gamma,\mu)$ generated by the 
$S_i$, which is isomorphic to ${\mathcal O}_A$. 
 
\begin{thm} 
The Cuntz--Krieger algebra ${\mathcal O}_A$ satisfies the following 
properties. 
\begin{enumerate} 
\item There is an injection ${\rm C}(\Lambda_\Gamma) \to {\mathcal 
O}_A$ which identifies ${\rm C}(\Lambda_\Gamma)$ with the maximal 
commutative subalgebra of ${\mathcal O}_A$ generated by the 
$P_\gamma$ as in \eqref{TandP}. 
\item The generators $S_i$ of ${\mathcal O}_A$ given in \eqref{Si} 
realize ${\mathcal O}_A$ as a subalgebra of ${\rm 
C}(\Lambda_\Gamma)\rtimes \Gamma$. 
\item The operators $T_\gamma$ of \eqref{TandP} are elements in 
${\mathcal O}_A$, hence the injection of ${\mathcal O}_A$ inside ${\rm 
C}(\Lambda_\Gamma)\rtimes \Gamma$ is an isomorphism, 
\begin{equation} \label{treefree} {\mathcal O}_A \cong 
C(\Lambda_\Gamma) \rtimes \Gamma. \end{equation} 
\end{enumerate} 
\label{OAprod} 
\end{thm} 
 
\noindent{\bf Proof.} 
{\em 1.} The algebra ${\rm C}(\Lambda_\Gamma)$ acts on 
$L^2(\Lambda_\Gamma,\mu)$ as multiplication operators. We identify the 
characteristic function $\chi_{\Lambda_\Gamma(\gamma)}$ of the subset 
$\Lambda_\Gamma(\gamma)\subset \Lambda_\Gamma$ with the projector 
$P_\gamma$ defined in \eqref{TandP}. A direct calculation shows that, 
for any $\gamma\in \Gamma$, the projector $P_\gamma$ satisfies 
$P_\gamma = S_{i_1} \cdots S_{i_k} S_{i_k}^* \cdots S_{i_1}^*$, for 
$\gamma = g_{i_1}\cdots g_{i_k}$, hence it is in the algebra 
${\mathcal O}_A$. 
 
For a multi-index $\mu=\{ i_1,\ldots, i_k \}$, the {\em range 
projection} $P_\mu$ is the element $P_\mu =S_\mu S_\mu^*$ in 
${\mathcal O}_A$. We have identified $\chi_{\Lambda_\Gamma(\gamma)}$ 
in ${\rm C}(\Lambda_\Gamma)$ with the range projection $P_\mu$, for 
$\mu=\{ i_1,\ldots, i_k \}$ the multi-index of $\gamma = g_{i_1}\cdots 
g_{i_k}$. Thus, we have identified ${\rm C}(\Lambda_\Gamma)$ with the 
maximal commutative subalgebra of ${\mathcal O}_A$ generated by the 
range projections (\cf \cite{Cu} \cite{CuKrie}). 
 
{\em 2.} The operators $S_i$ defined in \eqref{Si} determine 
elements in ${\rm C}(\Lambda_\Gamma)\rtimes \Gamma$, by identifying 
the projectors $P_j$ with $\chi_{\Lambda_\Gamma(g_j)}$ as 
in (1), and the operators $T_i$ with the corresponding elements in the 
unitary representation \eqref{AutGamma}. These elements still satisfy 
the relations \eqref{CK1rel} \eqref{CK2rel}, hence the algebra 
generated by the $S_i$ can be regarded as a subalgebra of ${\rm 
C}(\Lambda_\Gamma)\rtimes \Gamma$. 
 
{\em 3.} Given {\em 1.} and {\em 2.}, in order to prove the isomorphism 
\eqref{treefree}, it is enough to show that, for any $\gamma \in 
\Gamma$, the operators $T_\gamma$ in the unitary representation 
\eqref{AutGamma} are in the subalgebra of ${\rm 
C}(\Lambda_\Gamma)\rtimes \Gamma$ generated by the $S_i$. In fact, 
since this subalgebra contains ${\rm C}(\Lambda_\Gamma)$, if we know 
it also contains the $T_\gamma$, it has to be the whole of ${\rm 
C}(\Lambda_\Gamma)\rtimes \Gamma$. Again this follows by a direct 
calculation: $T_\gamma = T_{i_1} \cdots T_{i_k}$ and $T_i = S_{g+i}+ 
S_i^*$, with $g_{g+i}=g_i^{-1}$, since
$$ T_i f (x) = | g_i^\prime (g_i^{-1}x)|^{-\delta_H/2} f(g_i^{-1}x), $$
$$ S_i f(x)=(1- \chi_{g_i^{-1}}(x)) \, | g_i^\prime(x) |^{\delta_H/2}
f(g_i x), $$ 
and
$$ S_i^* f(x)=\chi_{g_i}(x) \, |g_i^\prime (Tx)|^{-\delta_H/2}
f(Tx). $$ 
 
\noindent $\diamond$

\subsection{Spectral triple}

We construct a spectral triple for the noncommutative space ${\mathcal
O}_A$. As a Hilbert space, we want to consider a space that contains
naturally a copy of the Archimedean  cohomology $H^1(\X)^{N=0}$ 
and of the ${\rm Coker}(N)$ in $\H^2(\XX)$. Instead of realizing these
cohomology spaces inside dynamical cohomology and homology as in
\eqref{hom-and-cohom}, we will work with two copies of the chain
complex for the dynamical cohomology, identifying $H^1(\X)^{N=0}$ with
a copy of the subspace ${\mathcal V}$ and the ${\rm Coker}(N)$ in
$\H^2(\XX)$ with the other copy of ${\mathcal V}$, via the  
duality isomorphism $\delta_1$ of Proposition \ref{N-isos}. 

On the Hilbert space ${\mathcal H}={\mathcal L}\oplus {\mathcal L}$,
we consider the unbounded linear operator 
\begin{equation}\label{D-oper3L}
D|_{{\mathcal L}\oplus 0} = \sum_n (n+1) \, (\hat\Pi_n \oplus
0)  \ \ \ \  
D|_{0\oplus {\mathcal L}} = -\sum_n n \,\, (0\oplus\hat\Pi_n).
\end{equation} 

\begin{rem}\label{Kclass}  {\em
Notice that the shift by $+1$ in the
spectrum of $D$ on the positive part ${\mathcal L}\oplus 0$ takes into
account the shift by one in the grading between
dynamical homology and cohomology, as in \eqref{tildedelta1}. This
reflects, in turn, the shift by one in grading introduced by the duality
isomorphism \eqref{isoms-cone} between the cohomology of kernel and
cokernel of the monodromy map in the cohomology of the cone, \cf
Proposition \ref{N-isos}. This shift introduces a spectral asymmetry
in the Dirac operator \eqref{D-oper3L}, hence a nontrivial eta 
invariant associated to the spectral triple. There is another possible
natural choice of the sign for the Dirac operator $D$ on ${\mathcal
H}={\mathcal L}\oplus {\mathcal L}$, instead of the one in
\eqref{D-oper3L}. Namely, one can define the sign of $D$ using the 
duality isomorphism $\delta_1$ of \eqref{isoms-cone}, instead of using
the sign of the operator $\Phi$ on the archimedean
cohomology and its dual. With this other choice, the sign would be
given by the operator that permutes the two copies of ${\mathcal L}$
in ${\mathcal H}$. The choice of sign determined by $\Phi$, as in 
\eqref{D-oper3L}, gives rise to a spectral triple that is degenerate as
$K$-homology class, hence, in this respect, using the sign induced by the
duality $\delta_1$ may be preferable. }
\end{rem} 

\smallskip

The operator $D$ of \eqref{D-oper3L} has the following properties:

\begin{prop}\label{comp-resol} 
Consider the data $({\mathcal H},D)$ as in 
\eqref{D-oper3L}. The operator $D$ is self adjoint. Moreover,
for all $\lambda\notin \R$, the resolvent 
$R_\lambda(D):=(D-\lambda)^{-1}$ is a compact operator on ${\mathcal 
H}$. The restriction of the operator $D$ to ${\mathcal V}\oplus 
{\mathcal V}\subset {\mathcal H}$ agrees with the Frobenius-type 
operator $\Phi$ of \eqref{Phi-cone}: 
$$ \begin{array}{rll} U\,\, \Phi|_{gr_{2p}^\w H^1(\X)^{N=0}}\,  U^{-1} & = 
D|_{0\oplus gr_{2p}^\Gamma {\mathcal V}} & p\leq 0 \\[3mm] 
 U \delta_1^{-1} \,\, \Phi|_{gr_{2p}^\w H^2(\XX)}\,\,  
\delta_1  U^{-1}
& =D|_{gr_{2(-p+2)}^\Gamma {\mathcal V} \oplus 0} & p\geq 2,
\end{array} $$  
with the map $U$ of \eqref{map-U} of Theorem
\ref{map-ar-dyn-2}, and $\delta_1$ the duality isomorphism of
\eqref{isoms-cone}.  
\end{prop} 

\noindent{\bf Proof.} The operator $D$ of
\eqref{D-oper3L} is already given in diagonal form and is clearly
symmetric with  
respect to the inner product of ${\mathcal H}$, hence it is 
self-adjoint on the domain ${\rm dom}(D)=\{ X\in {\mathcal H}: DX\in 
{\mathcal H} \}$. For $\lambda \notin \Sp(D)$, the operators 
$R_\lambda(D)$ are bounded, and related by the resolvent equation 
$$ R_\lambda(D)- R_{\lambda'}(D)= 
(\lambda'-\lambda)R_\lambda(D)R_{\lambda'}(D). $$ 
This implies that, if $R_\lambda(D)$ is compact for one $\lambda \notin 
\Sp(D)$, then all the other $R_{\lambda'}(D)$, $\lambda' \notin \Sp(D)$ 
are also compact. In our case we have $\Sp(D)=\Z$, with finite
multiplicities, hence, for instance, $R_{1/2}(D)$ is a compact operator with 
spectrum $\{ (n+1/2)^{-1}: n\in \Z \}\cup \{ 0 \}$.  
The multiplicities grow exponentially, as shown in \S
\ref{combinat}, namely the eigenspaces of $D$ have dimensions
$\dim E_{n+1} = 2g(2g-1)^{n-1} (2g-2)$ for $n\geq 1$, $\dim E_n =2g
(2g-1)^{-n-1} (2g-2)$, for $n\leq -1$ and $\dim E_0=2g$. Thus,
the Dirac operator will not be finitely
summable on ${\mathcal H}$, while the restriction of $D$ to ${\mathcal
V}\oplus {\mathcal V}$ has constant multiplicities.
Furthermore, it is clear that the restriction of $D$ to ${\mathcal V}$
agrees with the restriction of $\Phi$ to $H^1(\X)^{N=0}$. On the other hand, 
recall that the operator $\Phi$ defined as in \eqref{Phi-cone} acts on 
$gr_{2p}^\w H^2(\XX)$, for $p\geq 2$, as multiplication by 
$p-1$. The map $\tilde U$ of \eqref{map-tildeU} identifies $gr_{2p}^\w 
H^2(\XX)$ with $gr_{2(p-1)}^\Gamma {\mathcal W}$. Theorem
\ref{map-ar-dyn-2} shows that we get 
$$  U \delta_1^{-1}\, \Phi|_{gr_{2p}^\w H^2(\XX)} \,\delta_1  U^{-1} 
=  D |_{gr_{2(-p+2)}^\Gamma {\mathcal V}\oplus 0}. $$
 
\noindent $\diamond$

We also consider the diagonal action of the algebra ${\mathcal O}_A$
on ${\mathcal H}$, via the representation \eqref{TandP},
\begin{equation}\label{rhoH}
\rho: {\mathcal O}_A \to {\mathcal B}({\mathcal H}).
\end{equation}
The operator $D$ and the algebra ${\mathcal O}_A$ satisfy
the following compatibility condition.

\begin{prop}\label{bound-commOA}
Assume that the Hausdorff dimension of $\Lambda_\Gamma$ is $\delta_H
<1$. Then the set of elements $a\in {\mathcal O}_A$ for which the
commutator $[D, \rho(a)]$ is a bounded operator on ${\mathcal H}$ is 
norm dense in ${\mathcal O}_A$.
\end{prop}

\noindent{\bf Proof.} 
In addition to the operators $S_i$,
we consider, for $k\geq 0$, operators of the form
\begin{equation}\label{SiepsK}
(S_{i,k} f)(x) = h_{i,k}^{+}(x) \,
(1-\chi_{g_i^{-1}})(x) \, f(g_i x), 
\end{equation}
and
\begin{equation}\label{hatSiepsK}
(\hat S_{i,k} f)(x) = h_{i,k}^{-}(Tx) \,
\chi_{g_i}(x) \, f(Tx), 
\end{equation}
where $h_{i,k}^{\pm} \in {\mathcal P}_k$ is the function
\begin{equation}\label{hikeps}
h_{i,k}^{\pm}(x) = \Pi_k \left( | g_i^\prime (x) |^{\pm \delta_H /2}
\right).  
\end{equation}
Notice that these operators satisfy
\begin{equation}\label{SiNfiltr}
S_{i,k} : {\mathcal P}_{n+1} \to {\mathcal P}_{n}  \ \ \ \
\forall n \geq k.  
\end{equation}
\begin{equation}\label{hatSiNfiltr}
\hat S_{i,k} : {\mathcal P}_n \to {\mathcal P}_{n+1}  \ \ \ \
\forall n \geq k.  
\end{equation}
Moreover, the operators $S_{i,k}^*$ and $\hat
S_{i,k}^*$ satisfy
\begin{equation}\label{SiNfiltr2}
S_{i,k}^* : {\mathcal P}_n^\perp \to {\mathcal
P}_{n+1}^\perp  \ \ \ \ \forall n \geq k,  
\end{equation}
\begin{equation}\label{hatSiNfiltr2}
\hat S_{i,k}^* : {\mathcal P}_{n+1}^\perp \to {\mathcal
P}_{n}^\perp  \ \ \ \  \forall n \geq k.  
\end{equation}

We want to estimate the commutator $[D,S_i]$. We have
$$ [D,S_i]=\sum_k k [\Pi_k,S_i]- \sum_k k
[\Pi_{k-1},S_i] = -S_i (1- \Pi_0) + \sum_{k\geq 0}
(S_i \Pi_{k+1} - \Pi_k S_i). $$
We can write this in the form
$$  -S_i (1- \Pi_0) + \sum_{k\geq 0}
((S_i -S_{i,k}) \Pi_{k+1} - \Pi_k
(S_i -S_{i,k})) - \sum_{k\geq 0} \Pi_k
S_{i,k} (1-\Pi_{k+1}), $$
where, in the last term, we have used \eqref{SiNfiltr}. We further
write it as
$$  -S_i (1- \Pi_0) + \sum_{k\geq 0}
((S_i -S_{i,k}) \Pi_{k+1} - \Pi_k
(S_i -S_{i,k})) - \sum_{k\geq 0} \Pi_k
(S_{i,k} - \hat S_{i,k}^*) (1-\Pi_{k+1}), $$ 
using \eqref{hatSiNfiltr}. This means that we can estimate the first
sum in terms of the series
\begin{equation}\label{series1}
\sum_{k\geq 0} \| S_i  -S_{i,k}  \|,
\end{equation}
and the second sum in terms of the series
\begin{equation}\label{series2}
\sum_{k\geq 0} \| S_{i,k}  - \hat S_{i,k}^* \|, 
\end{equation}
That is, if \eqref{series1} and \eqref{series2} converge, then the
commutator $[D,S_i]$ is bounded in the operator norm, with a
bound given in terms of the sum of the series \eqref{series1} and
\eqref{series2}. 
The series \eqref{series1} can be estimated in terms of
the series 
\begin{equation}\label{seriesgprime}
\sum_{k\geq 0} \left\| |g_i^\prime|^{\delta_H/2} -
h_{i,k}^{+}  \right\|_\infty. 
\end{equation}

For $\gamma=a_0\ldots a_k$ and $f\in {\mathcal L}$, we have 
$$ (\Pi_k f)(\gamma) = \frac{1}{Vol(\Lambda_\Gamma(\gamma))}
\int_{\Lambda_\Gamma(\gamma)} f(x) d\mu(x), $$
where $\Lambda_\Gamma(\gamma)$ is the set ${\mathcal S}^+(a_0\ldots
a_k)$ of irreducible (half) infinite words that start with $\gamma$.
Thus, we can estimate \eqref{seriesgprime} by the series
\begin{equation}\label{seriesg}
\sum_k \,\, \sup_{\gamma\in \Gamma, |\gamma|=k}\,\, \sup_{x,y\in
\Lambda_\Gamma(\gamma)} \left| 
|g_i^\prime|^{\delta_H /2}(x) -
|g_i^\prime|^{\delta_H /2}(y) \right|. 
\end{equation}
If the generator $g_i$ of $\Gamma$ is represented by a matrix
$$ g_i = \left( \begin{array}{cc} a_i & b_i \\ c_i & d_i
\end{array}\right) \in \SL(2,\C), $$
then
$$ | g_i^\prime(x) |^{\delta_H /2} = |c_i x +
d_i|^{-\delta_H}. $$

Since for a Schottky group $\Gamma$ the limit set
$\Lambda_\Gamma$ is always strictly contained in $\P^1(\C)$, we can
assume, without loss of generality, that in fact $\Lambda_\Gamma \subset
\C$. Since $\infty \notin \Lambda_\Gamma$, the $|g_i'|$ have no poles on
$\Lambda_\Gamma$. Then the functions $|g_i'|$ are Lipschitz functions on 
$\Lambda_\Gamma$, satisfying an estimate
\begin{equation}\label{holder} 
 \left| |g_i^\prime|^{\delta_H /2}(x) -
|g_i^\prime|^{\delta_H/2}(y) \right| \leq C_i
|x-y|, 
\end{equation}
for some constant $C_i>0$.
Thus, we can estimate \eqref{seriesg} in terms of the series 
\begin{equation}\label{seriesdiam1}
\sum_k \sup_{|\gamma|=k} {\rm
diam}(\Lambda_\Gamma(\gamma)). 
\end{equation}
Since $\Lambda_\Gamma(\gamma) =\gamma
(\Lambda_\Gamma\smallsetminus \Lambda_\Gamma(\gamma^{-1}))$, and
$$ \frac{a_\gamma x + b_\gamma}{c_\gamma x + d_\gamma} -
\frac{a_\gamma y + b_\gamma}{c_\gamma y + d_\gamma} =
\frac{x-y}{(c_\gamma x + d_\gamma)(c_\gamma y + d_\gamma)}, $$
for
$$  \gamma  = \left( \begin{array}{cc} a_\gamma & b_\gamma \\ c_\gamma
& d_\gamma \end{array}\right) \in \SL(2,\C), $$
we can estimate, for $x,y\in \Lambda_\Gamma(\gamma)$, 
\begin{equation}\label{diamest}
 |x-y| \leq {\rm diam}(\Lambda_\Gamma) \,\, |\gamma'(x)|^{1/2} \,\,
|\gamma'(y)|^{1/2}. 
\end{equation}

We can then give a very crude estimate of the series
\eqref{seriesdiam1} in terms of the series
\begin{equation}\label{seriesdiam}
\sum_{\gamma\in \Gamma} {\rm
diam}(\Lambda_\Gamma(\gamma)). 
\end{equation}
Notice that \eqref{seriesdiam1} in fact has a much better convergence
than \eqref{seriesdiam}, by an exponential factor.

By \eqref{diamest} we can then estimate the series \eqref{seriesdiam}
in terms of the Poincar\'e series of the Schottky group
\begin{equation}\label{Poincare}
\sum_{\gamma\in\Gamma} |\gamma'|^s,   \ \ \ \ \
s=1 > \delta_H.
\end{equation}
It is known (\cf \cite{Bo}) that the Poincar\'e series
\eqref{Poincare} converges absolutely for all $s>\delta_H$. Thus,
we obtain that the series \eqref{series1} and \eqref{series2} 
converge and the commutator $[D,S_i]$ is bounded. 

The series \eqref{series2} can be estimated in a similar way, in terms
of the series 
$$ \sum_k \| h_{i,k}^{+} - |g_i'|^{\delta_H}
h_{i,k}^{-} \|_{\infty}, $$
which can be deal with by the same argument.

A completely analogous argument also shows that the commutators
$[D,S_i^*]$ are bounded, using the operator 
$\hat S_{i,k}$ instead of $S_{i,k}$ and
$S_{i,k}^*$ instead of $\hat S_{i,k}^*$ in the
calculation.

Thus, we have shown that the dense subalgebra ${\mathcal O}_A^{alg}$ 
generated algebraically by the operators $S_i$ and $S_i^*$ has bounded 
commutator with $D$. 

\noindent $\diamond$  

In order to treat also the case $\delta_H>1$, one would need a more
delicate estimate of \eqref{seriesdiam1}. 
To our purposes, the condition on the Hausdorff dimension $\delta_H<1$
is sufficiently general, since it is satisfied in the cases 
we are interested in, namely for an archimedean prime that is a real
embedding, with $X_{/\R}$ an orthosymmetric smooth real algebraic
curve, \cf  Proposition \ref{real-Bers} and Remark \ref{delta<1}.

\medskip

We can then construct the ``dynamical spectral triple'' at arithmetic
infinity as follows:

\begin{thm}\label{dyn-SP3OA}
For $\delta_H=\dim_H(\Lambda_\Gamma)<1$, 
the data $({\mathcal O}_A,{\mathcal H},D)$ form a spectral
triple. 
\end{thm}

\noindent{\bf Proof.}
By Proposition \ref{comp-resol} we know that $D$ is self adjoint
and $(1+D^2)^{-1/2}$ is compact. 
Proposition \ref{bound-commOA} then
provides the required compatibility between the Dirac operator 
$D$ and the algebra. 

\noindent $\diamond$  

\begin{rem} \label{FA-OA} {\em The dynamical spectral triple we
constructed in Theorem \ref{dyn-SP3OA} is not finitely summable. 
This is necessarily the case, if the Dirac operator has bounded
commutator with group elements in $\Gamma$, since a result of Connes
\cite{Connes2} shows that non amenable discrete groups (as is
the case for the Schottky group $\Gamma$) do not admit finitely
summable spectral triples. However, it is shown in \cite{Cu}, 
\cite{CuKrie}, and \cite{Putn} that the Cuntz--Krieger 
algebra ${\mathcal O}_A$ also admits a second description as a crossed
product algebra. Namely, up to stabilization (\ie tensoring with compact 
operators) we have 
\begin{equation} \label{AF-T} {\mathcal O}_A \simeq 
{\mathcal F}_A \rtimes_T \Z, \end{equation} 
where ${\mathcal F}_A$ is an approximately finite 
dimensional (AF) algebra stably isomorphic to the groupoid 
${\rm C}^*$--algebra ${\rm C}^*({\mathcal G}^u)$ of Remark 
\ref{SmaleCalg}. It can be shown that the shift operator $T$ induces 
an action by automorphisms on ${\mathcal F}_A$, so that the crossed 
product algebra \eqref{AF-T} corresponds to ${\rm C}^*({\mathcal 
G}^u)\rtimes_T \Z$ associated to the Smale space $({\mathcal 
S},T)$. (\cf Remark \ref{SmaleCalg}). This means that, again by
Connes' result on hyperfiniteness (\cf \cite{Connes2}), it may be
possible to construct a finitely summable spectral triple using
the description \eqref{AF-T} of the algebra. It
is an interesting question whether the construction of a finitely
summable triple can be carried out in a way that is of arithmetic
significance. } 
\end{rem}

\subsection{Archimedean  factors from dynamics} 

The dynamical spectral triple we constructed in Theorem
\ref{dyn-SP3OA} is not finitely summable.
However, it is still possible to recover from these data the local
factor at arithmetic infinity.

\medskip

As in the previous sections, we consider a fixed Archimedean  prime 
given by a real embedding $\alpha: \K\hookrightarrow \R$, such 
that the corresponding Riemann surface $X_{/\R}$ is an 
orthosymmetric smooth real algebraic curve of genus $g\geq 2$. The 
dynamical spectral triple provides another interpretation of the 
Archimedean  factor $L_\R(H^1(X_{/\R},\R),s)=\Gamma_\C (s)^g$. 

\begin{prop}\label{L-factor2} 
Consider the zeta functions
\begin{equation}\label{zeta-proj1} 
\zeta_{\pi({\mathcal V}), D} (s,z):= 
\sum_{\lambda \in \Sp(D)} \Tr \left( \pi({\mathcal V}) 
\Pi(\lambda, D) \right) (s-\lambda)^{-z}, 
\end{equation} 
for $\pi({\mathcal V})$ the orthogonal projection on
the norm closure of $0\oplus {\mathcal V}$ in ${\mathcal H}$, and 
\begin{equation}\label{zeta-proj2} 
\zeta_{\pi({\mathcal V}, \bar F_\infty=id), D} (s,z):= 
\sum_{\lambda \in \Sp(D)} \Tr \left( \pi({\mathcal V}, \bar 
F_\infty=id) \Pi(\lambda, D) \right) (s-\lambda)^{-z}, 
\end{equation} 
for $\pi({\mathcal V}, \bar F_\infty=id)$ the orthogonal 
projection on the norm closure of $0\oplus {\mathcal V}^{\bar 
F_\infty=id}$. The corresponding regularized 
determinants satisfy 
\begin{equation}\label{L-det21} 
\exp\left( - \frac{d}{dz} \zeta_{\pi({\mathcal V}), D/2\pi} 
(s/2\pi,z)|_{z=0} \right)^{-1} =  
L_\C(H^1(X),s),
\end{equation}
\begin{equation}\label{L-det22} 
\exp\left( - \frac{d}{dz} \zeta_{\pi({\mathcal V}, \bar 
F_\infty=id), D/2\pi} (s/2\pi,z)|_{z=0} \right)^{-1} = 
L_\R(H^1(X),s). 
\end{equation} 
Moreover, the operator $\pi({\mathcal V})$ acts as 
projections in the AF algebra ${\mathcal F}_A$
compressed by the spectral projections
$\Pi(\lambda, D)$.   
\end{prop} 

\noindent{\bf Proof.} Let ${\mathcal V}_{even}=\oplus_{p=2k} 
gr_{2p}^\Gamma {\mathcal V}$ and ${\mathcal 
V}_{odd}=\oplus_{p=2k+1} gr_{2p}^\Gamma {\mathcal V}$. Further, we 
denote by ${\mathcal V}_{even}^\pm$ and ${\mathcal V}_{odd}^\pm$ 
the $\pm 1$ eigenspaces of the change of orientation involution. 
The $+1$ eigenspace of the Frobenius $\bar F_\infty$ is then given 
by ${\mathcal V}_{even}^+ \oplus {\mathcal V}_{odd}^-$. Let 
$\pi({\mathcal V}_{even}^+)$ and $\pi({\mathcal V}_{odd}^-)$ 
denote the corresponding orthogonal projections. 
 
We then compute explicitly 
$$ \zeta_{\pi({\mathcal V}), D/2\pi} (s/2\pi,z) =
\sum_{\lambda \in \Sp(D)} \Tr \left( \pi({\mathcal V}) 
\Pi(\lambda, D) \right) (s-\lambda)^{-z} = 2g (2\pi)^z \zeta(s,z) $$
and  
$$ \begin{array}{c} 
\zeta_{\pi({\mathcal V}, \bar F_\infty=id), D/2\pi} (s/2\pi,z) = 
\\[2mm] 
\sum_{\lambda \in \Sp(D)} \Tr \left( \pi({\mathcal V}_{even}^+) 
\Pi(\lambda, D) \right) (s-\lambda)^{-z} + \sum_{\lambda \in 
\Sp(D)} \Tr \left( \pi({\mathcal V}_{odd}^-) \Pi(\lambda, D) 
\right) 
(s-\lambda)^{-z}\\[2mm] = 
g \pi^z \zeta(s/2,z) + g \pi^z \zeta ((s-1)/2,z), \end{array} $$ 
so that we obtain the identities \eqref{L-det21} and \eqref{L-det22} 
as in Proposition  \ref{detR}. 

Consider the operators $Q_{i,n}=S_i^n {S_i^*}^n$ in the 
Cuntz--Krieger algebra ${\mathcal O}_A$. These are projections 
in ${\mathcal F}_A$ that act as multiplication by the 
characteristic function $\chi_{{\mathcal S}^+(w_{n,i})}$. Thus, the
operator $Q_{|p|}:=\sum_i Q_{i,|p|}$ has the property that the
compression $\hat\Pi_{|p|}Q_{|p|} \hat\Pi_{|p|}$ by the spectral
projections of $D$, acts as the orthogonal projection onto
$0\oplus gr^\Gamma_{2p} {\mathcal V}$. 

\noindent $\diamond$

After recovering the Archimedean  factor, one can ask 
a more refined question, namely whether it is possible
to recover the algebraic curve $X_{/\R}$ from the non-commutative
data described in this section. One can perhaps relate the crossed
product algebra ${\mathcal F}_A\rtimes_T \Z$ of 
\eqref{AF-T} to some geodesic lamination on $X_{/\R}$ that lifts, under
the Schottky uniformization map of Proposition \ref{real-Bers}, to a
collection of hyperbolic geodesics in $\P^1(\C)\backslash \P^1(\R)$
with ends on $\Lambda_\Gamma \subset \P^1(\R)$. This may provide an
approach to determining the periods of the curve 
and also bridge between the two constructions presented in the first
and second half of this paper. We hope to return to these
ideas in the future.

\section{Reduction mod $\infty$ and homotopy quotients} \label{modinftysect} 

In the previous sections we have described the (noncommutative) 
geometry of the fiber at arithmetic infinity of an arithmetic surface
in terms of its dual graph, which we obtained from two quotient
spaces: the spaces
\begin{equation}\label{quotients}
\Lambda_\Gamma /\Gamma \ \ \ \text{ and } \ \ \ \Lambda_\Gamma
\times_\Gamma \Lambda_\Gamma \simeq {\mathcal S}/\Z,  
\end{equation}
with $\Z$ acting via the invertible shift $T$, \cf \S
\ref{sectcoding}, which we can think of as the sets of vertices and 
edges of the dual graph. We analyzed their noncommutative geometry in
terms of Connes' theory of spectral triples. 

\smallskip

Another fundamental construction in noncommutative geometry (\cf
\cite{Connes-tr}) is that of {\em homotopy
quotients}. These are commutative spaces, which provide, up to
homotopy, geometric models for the corresponding noncommutative
spaces. The noncommutative spaces themselves, as we are going to show
in our case, appear as quotient spaces of foliations on the homotopy
quotients with contractible leaves.

\smallskip

The crucial point in our setting is that the homotopy quotient for the
noncommutative space ${\mathcal S}/\Z$ is precisely the mapping torus
\eqref{suspensionT} which gives the geometric model of the dual graph,
\begin{equation}\label{htpyquot2}
{\mathcal S}_T = {\mathcal S}\times_\Z \R,
\end{equation}
where the noncommutative space ${\mathcal S}/\Z$ can be
identified with the quotient space of the natural foliation on
\eqref{htpyquot2} whose generic leaf is contractible (a copy of $\R$).
On the other hand, the case of the noncommutative space 
$\Lambda_\Gamma /\Gamma$ is also extremely interesting. In fact, in
this case the homotopy quotient appears very naturally
and it describes what Manin refers to in \cite{Man} as the
``reduction mod $\infty$''.  

\medskip 

We recall briefly how the reduction map works in the non-Archimedean  
setting of Mumford curves (\cf \cite{Man} \cite{Mum}). Let ${\mathfrak 
K}$ be a finite extension of $\Q_p$ and let ${\mathfrak a}$ be its 
ring of integers. The correct analog for the 
Archimedean  case is obtained by ``passing to a limit', replacing 
${\mathfrak K}$ with its Tate closure in $\C_p$ (\cf \cite{Man} \S 
3.1), however, for our purposes here it is sufficient to illustrate 
the case of a finite extension. 
 
The role of the hyperbolic space $\H^3$ in the non-Archimedean  case is 
played by the Bruhat-Tits tree ${\mathcal T}_{BT}$  with vertices 
$$ {\mathcal T}^0_{BT} =\{ {\mathfrak a}-\text{lattices of rank 2 in a 
2-dim ${\mathfrak K}$-space} \}/ {\mathfrak K}^*. $$ 
Vertices 
in ${\mathcal T}_{BT}$ have valence $| \P^1({\mathfrak a}/{\mathfrak 
m})|$, where ${\mathfrak m}$ is the maximal ideal. Each edge in 
${\mathcal T}_{BT}$ has length $\log |{\mathfrak a}/{\mathfrak 
m}|$. The set of ends of ${\mathcal T}_{BT}$ is identified with 
$X({\mathfrak K})=\P^1({\mathfrak K})$. This is the analog of the 
conformal boundary $\P^1(\C)$ of $\H^3$. Geodesics correspond to 
doubly infinite paths in ${\mathcal T}_{BT}$ without backtracking. 
 
Fix a vertex $v_0$ on ${\mathcal T}_{BT}$. This corresponds to the 
closed fiber $X_{\mathfrak a} \otimes ({\mathfrak a}/{\mathfrak 
m})$ for the chosen ${\mathfrak a}$-structure $X_{\mathfrak a}$. Each 
$x\in \P^1({\mathfrak K})$ determines a unique choice of a subgraph 
$e(v_0,x)$ in ${\mathcal T}_{BT}$ with vertices $(v_0, v_1, v_2, 
\ldots )$ along the half infinite path without backtracking which has 
end $x$. The subgraphs $e(v_0,x)_k$ with vertices $(v_0, v_1,\ldots 
v_k)$ correspond to the reduction mod ${\mathfrak m}^k$, namely 
$$ \{ e(v_0,x)_k: x\in  X({\mathfrak K}) \} \leftrightsquigarrow 
X_{\mathfrak a}( {\mathfrak a}/{\mathfrak m}^k). $$ 
Thus the finite graphs $e(v_0,x)_k$ represent ${\mathfrak 
a}/{\mathfrak m}^k$ points, and the infinite graph $e(v_0,x)$ 
represents the reduction of $x$. 
 
\smallskip 
 
A Schottky group $\Gamma$, in this non-Archimedean  setting, is a purely 
loxodromic free discrete subgroup of $\PSL(2,{\mathfrak K})$ in $g$ 
generators. The doubly infinite 
paths in ${\mathcal T}_{BT}$ with ends at the pairs of fixed points 
$x^\pm(\gamma)$ of the elements $\gamma\in \Gamma$ produce a copy of the 
combinatorial tree ${\mathcal T}$ of the group $\Gamma$ in ${\mathcal 
T}_{BT}$. This is the analog of regarding $\H^3$ as the union of the 
translates of a fundamental domain for the action of the Schottky 
group, which can be thought of as a `tubular neighborhood' of a copy 
of the Cayley graph ${\mathcal T}$ of $\Gamma$ embedded in $\H^3$. 
The ends of the tree ${\mathcal T}\subset {\mathcal T}_{BT}$ 
constitute the limit set $\Lambda_\Gamma \subset \P^1({\mathfrak K})$. 
The complement $\Omega_\Gamma = \P^1({\mathfrak 
K})\smallsetminus \Lambda_\Gamma$ 
gives the uniformization of the Mumford curve $X({\mathfrak K})\simeq 
\Omega_\Gamma /\Gamma$. In turn, $X({\mathfrak K})$ can be identified 
with the ends of the quotient graph ${\mathcal T}_{BT}/\Gamma$, just 
as in the Archimedean  case the Riemann surface is the 
conformal boundary at infinity of the handlebody $\mX_\Gamma$. 
 
The reduction map is then obtained by considering the half infinite paths 
$e(v,x)$ in ${\mathcal T}_{BT}/\Gamma$ that start at a vertex $v$ of 
the finite graph ${\mathcal T}/\Gamma$ and whose end $x$ is 
a point of $X({\mathfrak K})$, while the finite graphs $e(v,x)_k$ 
provide the ${\mathfrak a}/{\mathfrak m}^k$ points. 
 
\medskip 
 
This suggests that the correct analog of the reduction map in the 
Archimedean  case is obtained by considering geodesics in $\H^3$ with 
an end on $\Omega_\Gamma$ and the other on $\Lambda_\Gamma$, as described in 
\cite{Man}. Arguing as in Lemma \ref{geod-coding}, we see that the set 
of such geodesics can be identified with the quotient $\Omega_\Gamma 
\times_\Gamma \Lambda_\Gamma$. The analog of the finite graphs 
$e(v,x)_k$ that define the reductions modulo ${\mathfrak m}^k$ is then
given by the quotient $\H^3 \times_\Gamma \Lambda_\Gamma$. 
 
\smallskip 

Notice then that the quotient space
\begin{equation}\label{htpyquot1}
\Lambda_\Gamma \times_\Gamma \H^3 = \Lambda_\Gamma \times_\Gamma
\underline{E}\Gamma, 
\end{equation}
is precisely the homotopy quotient of $\Lambda_\Gamma$ with respect to
the action of $\Gamma$, with $\underline{E}\Gamma =\H^3$ and the
classifying space $\underline{B}\Gamma = \H^3/\Gamma =\mX_\Gamma$,
(\cf \cite{Connes-tr}). In this case also we find
that the noncommutative space $\Lambda_\Gamma /\Gamma$ is the quotient
space of a foliation on the homotopy quotient \eqref{htpyquot1} with
contractible leaves $\H^3$.

\smallskip

The relation between the noncommutative spaces \eqref{quotients} and
the homotopy quotients \eqref{htpyquot2} \eqref{htpyquot1} is an
instance of a very general and powerful construction, namely the
$\mu$-map (\cf \cite{BaumConnes} \cite{Connes-tr}). In particular, in
the case of the noncommutative space ${\rm C}({\mathcal S})\rtimes_T
\Z$, the $\mu$-map 
\begin{equation}\label{Thomiso}
\mu: K^{*+1} ({\mathcal S}_T)\cong H^{*+1}({\mathcal S}_T,\Z) \to
K_*({\rm C}({\mathcal S})\rtimes_T \Z)
\end{equation}
is the Thom isomorphism that gives the identification of
\eqref{isomCTH1}, \eqref{isomCTH1-map} and recovers the
Pimsner--Voiculescu exact sequence \eqref{PV} as in \cite{Co-Thom}.  
The map $\mu$ of \eqref{Thomiso} assigns to a $K$-theory class 
${\mathcal E}\in K^{*+1} ({\mathcal S}\times_\Z \R)$ the index of the
longitudinal Dirac operator $\dirac_{\mathcal E}$ with coefficients
${\mathcal E}$. This index is an element of the $K$-theory of the
crossed product algebra ${\rm C}({\mathcal S})\rtimes_T
\Z$ and the $\mu$-map is an isomorphism.
Similarly, in the case of the noncommutative space ${\rm
C}(\Lambda_\Gamma) \rtimes \Gamma$, where we have a foliation
on the total space with leaves $\H^3$, the $\mu$-map 
\begin{equation}\label{mu-map}
\mu: K^{*+1} (\Lambda_\Gamma \times_\Gamma \H^3) \to
K_*({\rm C}(\Lambda_\Gamma)\rtimes \Gamma)
\end{equation}
is again given by the index of the longitudinal Dirac operator
$\dirac_{\mathcal E}$ with coefficients ${\mathcal E} \in K^{*+1}
(\Lambda_\Gamma \times_\Gamma \H^3)$. In this case the map is an
isomorphism because the Baum--Connes conjecture with coefficients
holds for the case of $G=SO_0(3,1)$, with $\H^3=G/K$ and
$\Gamma\subset G$ the Schottky group, \cf \cite{Kas}. 

\smallskip

In particular, analyzing the noncommutative space ${\rm
C}(\Lambda_\Gamma)\rtimes \Gamma$ from the point of view of the theory
of spectral triples provides cycles to pair with $K$-theory classes
constructed geometrically via the $\mu$-map.

\smallskip

To complete the analogy with the reduction map in the case of
Mumford curves, one should also consider the half infinite paths
$e(v,x)$ corresponding to the geodesics in $\mX_\Gamma$ parameterized
by $\Lambda_\Gamma \times_\Gamma \Omega_\Gamma$, in addition to the
finite graphs $e(v,x)_k$ that correspond to the homotopy quotient
\eqref{htpyquot2}. This means that the space that completely describes
the ``reduction modulo infinity'' is a compactification of the
homotopy quotient 
\begin{equation}\label{htpyquotC}
\Lambda_\Gamma \times_\Gamma (\H^3\cup \Omega_\Gamma),
\end{equation}
where  $\overline{\underline{E}\Gamma}=\H^3\cup
\Omega_\Gamma$ corresponds to the compactification of the classifying space
$\underline{B}\Gamma = \H^3/\Gamma =\mX_\Gamma$ to
$\overline{\underline{B}\Gamma }= (\H^3\cup \Omega_\Gamma)/\Gamma =
\mX_\Gamma\cup X_{/\C}$, obtained by adding the conformal boundary at
infinity of the hyperbolic handlebody. This is not the only instance
where it is natural to consider compactifications for 
$\underline{E}\Gamma$ and the homotopy quotients, \cf \eg
\cite{Weinb}.  

\medskip

\section{Further structure at arithmetic infinity}\label{7} 
 
We point out some other interesting 
aspects of the relation we have developed between the two approaches 
of Manin and Deninger to the theory of the closed fiber at arithmetic 
infinity. We hope to return to them in future work. 
 
\bigskip 
 
\noindent {\it (i) Solenoids.} 
The infinite tangle of bounded geodesics in $\mX_\Gamma$ can also be 
related to generalized solenoids, thus making another connection between 
Manin's theory of the closed finer at infinity and Deninger's theory 
of ``arithmetic cohomology'', for which the importance of solenoids was 
stressed in the remarks of \S 5.8 of \cite{Den2}. 
 
There are various ways in which the classical $p$--adic solenoids can 
be generalized. One standard way of producing more general solenoids is by 
considering a map $\phi : C \to C$, which is a continuous surjection 
of a Cantor set $C$, and form the quotient space 
$(C\times [0,1])/\sim$ by the identification $(x,0)\sim 
(\phi(x),1)$. The resulting space is a generalized solenoid provided 
it is an {\it indecomposable} compact connected space, where the 
indecomposable condition means that the image under the quotient map 
of any $D\times [0,1]$ where $D$ is a clopen subset of $C$ is 
connected if and only if $D=C$. It is not hard to check that this 
condition is satisfied if $\phi$ is a homeomorphism of a Cantor set 
that has a dense orbit. Thus, our model for the dual graph of the 
fiber at arithmetic infinity, given 
by the mapping torus ${\mathcal S}_T$, is a generalized solenoid. 
This can also be related to generalized solenoids for subshifts 
of finite type \cite{Wil}. 
 
\bigskip 
 
\noindent {\it (ii) Dynamical zeta function.} The action of the 
shift operator $T$ on the limit set $\Lambda_\Gamma$ provides a 
dynamical zeta function on the special fiber at arithmetic 
infinity, in terms of the Perron--Frobenius theory for the 
shift $T$, with the Ruelle transfer operator \eqref{Ruelle}. 
\begin{equation}\label{Ruelle} 
({\mathcal R}_{-s f} \,  g)(z):= \sum_{y: Ty=z} e^{-s f(y)} g(y), 
\end{equation} 
depending on a parameter $s\in \C$ and with the function $f(z)=\log | 
T'(z) |$. This is the analog of the Gauss--Kuzmin operator studied in 
\cite{ManMar} \cite{Mar} in the case of modular curves. In fact, 
for $s=\delta_H$, the Hausdorff dimension of the limit set, the
operator ${\mathcal R}:={\mathcal R}_{-\delta f}$, with $f(z)=\log | 
T'(z) |$ is the Perron--Frobenius operator of $T$, that is, the
adjoint of composition with $T$, 
\begin{equation}\label{adjoint}
\int_{\Lambda_\Gamma} h(x) {\mathcal R}f (x) \, d\mu(x) =
\int_{\Lambda_\Gamma} h(Tx)\, f(x) \, d\mu(x), 
\end{equation}
for all $h,f\in L^2(\Lambda_\Gamma,\mu)$ and $\mu$ the
Patterson--Sullivan measure \eqref{PSmeas}. In fact, 
$$ \int_{\Lambda_\Gamma} h(Tx)\, f(x) \, d\mu(x) = \sum_i
\int_{\Lambda_\Gamma(g_i^{-1})} h(g_i x) f(x) \, d\mu(x) 
 = \sum_i \int_{\Lambda_\Gamma} h(g_i x)\, \chi_{g_i^{-1}}(x)\, f(x)
\, d\mu(x), $$ 
where $\chi_\gamma(x)$ is the characteristic function of the cylinder 
$\Lambda_\Gamma(\gamma)$ of all (right) infinite reduced words 
that begin with the word $\gamma$. This then gives
$$ = \int_{\Lambda_\Gamma} h(x) \sum_{i} A_{ij} \chi_{g_i^{-1}}(x)\, 
f(g_i^{-1}(x)) \, | g_i^\prime (g_i^{-1}(x))|^{-\delta_H} \, d\mu(x), $$
where the sum is over all admissible $i$'s, namely such that $g_i\neq
g_j$, where $x=g_ja_1\ldots a_n \ldots$.  

\smallskip

The operator ${\mathcal R}$ encodes important information on the
dynamics of $T$. In fact, there exists a Banach space $({\mathbb V},\|
\cdot \|_{\mathbb V})$ of functions on $\Lambda_\Gamma$ with the
properties that $\| f\|_\infty\leq \| f \|_{\mathbb V}$ for all $f\in
{\mathbb V}$, and ${\mathbb V}\cap {\rm C}(\Lambda_\Gamma)$ is dense
in ${\rm C}(\Lambda_\Gamma)$. On this space ${\mathcal R}$ is bounded
with spectral radius $r({\mathcal R})=1$. The point $\lambda=1$ is a
simple eigenvalue, with (normalized) eigenfunction the density of the
unique $T$-invariant measure on $\Lambda_\Gamma$. 

\smallskip

Similarly, for ${\mathcal  
B}_q$ a function space of $q$--forms on ${\mathcal U}\subset 
\P^1(\C)$, with ${\mathcal U}=\cup_i D_i^{\pm}$, with analytic 
coefficients and uniformly bounded, we denote by ${\mathcal R}_{s}^{(q)}$ 
the operator \eqref{Ruelle} acting on ${\mathcal B}_q$. 
In the case $q=0$, for $s=\delta_H$ the Hausdorff dimension of 
$\Lambda_\Gamma$, the operator ${\mathcal R}_{\delta_H}^{(0)}$ has top
eigenvalue $1$. This is a simple eigenvalue, with a non--negative
eigenfunction which is the density of the invariant measure. The 
analysis is similar to the case described in \cite{ManMar}. 
For ${\rm Re}(s) >>0$ the operators ${\mathcal R}_{s}^{(q)}$ 
are nuclear. If we denote by $\{ \lambda_{i,q}(s) \}_i$ the 
eigenvalues of ${\mathcal R}_{s}^{(q)}$, we can define 
\begin{equation} \label{Fi} P_q(s) = \det( 1-{\mathcal R}_{s}^{(q)} ) = 
\prod_i (1-\lambda_{i,q}(s)), \end{equation} 
so that we have a Ruelle zeta function 
\begin{equation} \label{ZYRuelle} 
Z_\Gamma (s) := \frac{P_1(s)}{P_0(s) P_2(s)}. 
\end{equation} 
By the results of \cite{Polli} \S 4, the Ruelle zeta function 
$Z_\Gamma (s)$ 
defined above is directly related to the Selberg zeta function of the 
hyperbolic handlebody $\mX_\Gamma$ of the form 
$$ Z(s)= \prod_{\gamma \in \Gamma} (1-N(\gamma)^{-s} ), $$ 
for $N(\gamma)=\exp(\ell_\gamma)$ the length of the corresponding 
primitive closed geodesic (equivalently $N(\gamma)$ is the norm of 
the derivative of $\gamma$ at the repelling fixed point 
$z^+(\gamma)\in \P^1(\C)$). 
If we consider the set of $s\in \C$ where $1$ is an eigenvalue of 
${\mathcal R}_{s}^{(q)}$, these produce possible zeroes and poles of the 
zeta function $Z_\Gamma (s)$: among these, the Hausdorff dimension 
$s=\delta_H$. The set of complex numbers where $Z_\Gamma (s)$ has a pole 
provides this way another notion of dimension spectrum of the fractal 
$\Lambda_\Gamma$ extending the ordinary Hausdorff dimension. 
 
It should be interesting to study the properties of $Z_\Gamma 
(s)$ in relation to the properties of the arithmetic zeta function 
\eqref{ZY} we defined in \S \ref{ReiZsec}. 
 
\bigskip 
 
\noindent {\it (iii) Missing Galois theory.} 
The data of a spectral triple $({\mathcal A}, {\mathcal 
H}, D)$ associated to the special fiber at arithmetic infinity may provide 
a different approach to what Connes refers to as a {\it missing Galois 
theory at Archimedean  places}. The idea is derived from the 
general setting of Hilbert's twelfth problem, for a local or global 
field $K$, and its maximal abelian (separable) extension 
$K^{ab}$. While class field theory provides a description of 
$Gal(K^{ab}/K)$, Hilbert's twelfth problem addresses the question of 
providing explicit generators of $K^{ab}$ and an explicit action of 
the Galois group, much as in the Kronecker--Weber case over $\Q$. 
The approach to Hilbert's twelfth problem given by Stark's conjectures 
(see \eg \cite{Tate}) consists of considering a family of zeta 
functions (\eg for the Kronecker--Weber case these consist of 
arithmetic progressions $\zeta_{n,m}(s):=\sum_{k\in m+n\Z} |k|^{-s}$) 
and generate numbers, which have the form of zeta--regularized 
determinants (\eg $\exp(\zeta_{n,m}'(0))$ for the Kronecker--Weber 
case). For more general fields it is a conjecture that such numbers are 
algebraic and that they provide the desired generators with Galois 
group action. In \cite{Man3} \cite{Man5} Manin conjectured 
that, in the case of real quadratic fields, 
noncommutative geometry should 
produce Stark numbers and Galois action 
via the theory of noncommutative tori. 
 
In our setting, whenever we associate the structure of a spectral 
triple to the Archimedean  places, we obtain, in particular, 
a family of zeta functions \eqref{aDzeta} 
determined by the spectral triple. These provide numbers of the form 
considered by Stark, 
\begin{equation}\label{DSzeta-Stark} 
 \exp\left( \frac{d}{ds} \zeta_{a,D}(s) |_{s=0} \right), 
\end{equation} 
as regularized determinants. Of course, in general there is no reason 
to expect these numbers to be algebraic. However, it may be an 
interesting question whether there are spectral triples for which they 
are, and whether the original Stark numbers arise via a spectral triple. 
Notice that there is an explicit action of the algebra ${\mathcal A}$ 
on the collection of numbers \eqref{DSzeta-Stark}, which may provide 
the right framework for a Galois action. 
We hope to return to this in a future work.

\end{document}